\documentclass[final,review,11pt]{my-article}

\usepackage{etex} 
\usepackage{graphicx,subfigure,color,colordvi}
\usepackage{latexsym,amsbsy,epsfig,fancybox}
\usepackage{amstext,mathrsfs}
\usepackage{amsfonts,textgreek}
\usepackage{mathrsfs}
\usepackage{stmaryrd,dsfont}
\usepackage{bbm}
\usepackage{url}
\usepackage{wrapfig,bm}
\usepackage{tikz}
\usepackage{citesort}
\usepackage{pstricks}
\usepackage{caption}
\usepackage{psfrag}
\usepackage{upgreek}
\usepackage{isomath}
\usepackage[dvips]{geometry}
\usepackage{latexsym}
\usepackage{graphicx}
\usepackage{array}
\usepackage{multirow}
\usepackage{booktabs}
\usepackage{colortbl}
\usepackage{verbatim}
\usepackage[colorlinks=true,citecolor=blue]{hyperref}
\usepackage{epstopdf,overpic}
\usepackage{enumerate}
\usepackage{amsmath,amssymb,xspace}

\newcommand{\eref}[1]{(\ref{#1})}
\newcommand{\fref}[1]{Fig.~\ref{#1}}
\newcommand{\tref}[1]{Table~\ref{#1}}
\newcommand{\sref}[1]{Section~\ref{#1}}
\newcommand{\pref}[1]{Problem~\ref{#1}}
\newcommand{\vm}[1]{\bm{#1}} 
\newcommand{\bsym}[1]{\bm{#1}}

\renewcommand{\Re}{{\rm{I\!R}}}
\newcommand{\vx}{\vm{x}}
\newcommand{\mat}{\mathsf} 
\newcommand{\transpose}{\mathsf{T}}
\newcommand{\diffx}{\mathrm{d}\vm{x}}

\newcommand{\diffs}{\mathrm{d}s}

\newcommand\Prefixd[2]{\vphantom{#2}#1#2}
\newcommand\Prefixf[4]{\vphantom{#4}#1#2#3#4}
\DeclareMathOperator*{\trace}{tr}

\newcommand{\smat}[2][ccccccccccccccccccccccccccccccccccccccccccccccccccc]{\left
[\begin{array}{#1}#2 \\ \end{array} \right]}

\tikzstyle{nicebox}=[draw=black!100, fill=white!10, rectangle, inner sep=4pt, inner ysep=16pt]
\tikzstyle{niceboxtitle}=[draw=black!100, fill=white, text=black, rectangle]

\usetikzlibrary{arrows,decorations.markings}

\newcommand{\assembly}{\operatornamewithlimits{\bm{\mat{A}}}}

\definecolor{forestgreen}{RGB}{34, 139, 34}
\definecolor{lightgray}{gray}{0.92}

\newcommand{\alejandro}[1]{{\color{black}{#1}}}

\newtheorem{theorem}{Theorem}

\newtheorem{problem}[theorem]{Problem}

\biboptions{sort&compress}

\begin{document}

\begin{frontmatter}
\title{A volume-averaged nodal projection method for the Reissner-Mindlin plate model}

\journal{Comp. Meth. Appl. Mech. Engr.}

\author[adr1,adr2]{A.~Ortiz-Bernardin\corref{cor1}}
\ead{aortizb@uchile.cl}
\author[adr1,adr2]{P.~K\"{o}brich}
\author[adr3]{J.~S.~Hale}
\author[adr1,adr2]{E.~Olate-Sanzana}
\author[adr3,adr4]{S.~P.~A.~Bordas}
\author[adr5]{S.~Natarajan}

\cortext[cor1]{Corresponding author. Tel: +56 2 2978 4664, Fax: +56 2 2689 6057,}

\address[adr1]{Department of Mechanical Engineering, Universidad de Chile, Av. Beauchef 851, Santiago 8370456, Chile. }

\address[adr2]{Computational and Applied Mechanics Laboratory, Center for Modern Computational Engineering, Facultad de Ciencias F\'isicas y Matem\'aticas, Universidad de Chile, Av. Beauchef 851, Santiago 8370456, Chile.}

\address[adr3]{Institute of Computational Engineering, University of Luxembourg, Maison du Nombre, 6 Avenue de la Fonte, 4364 Esch-sur-Alzette, Luxembourg.}

\address[adr4]{Institute of Mechanics and Advanced Materials, School of Engineering, Cardiff University, Cardiff CF24 3AA, Wales, UK.}

\address[adr5]{Department of Mechanical Engineering, Indian Institute of Technology, Madras, Chennai - 600036, India.}

\begin{abstract}

We introduce a novel meshfree Galerkin method for the solution of Reissner-Mindlin plate problems that is written in terms of the primitive
variables only (i.e., rotations and transverse displacement) and is devoid of shear-locking. The proposed approach uses linear maximum-entropy basis functions for field variables approximation and is built variationally on a two-field potential energy functional wherein the shear strain, written in terms of the primitive variables, is computed via a volume-averaged nodal projection operator that is constructed from the Kirchhoff constraint of the three-field mixed weak form. The meshfree approximation is constructed over a set of scattered nodes that are obtained from an integration mesh of three-node triangles on which the meshfree stiffness matrix and nodal force vector are numerically integrated. The stability of the method is rendered by adding bubble-like enrichment to the rotation degrees of freedom. Some benchmark problems are presented to demonstrate the accuracy and performance of the proposed method for a wide range of plate thicknesses.

\end{abstract}
\begin{keyword}
meshfree methods \sep maximum-entropy approximation \sep Reissner-Mindlin plate \sep shear-locking \sep \texttt{VANP} method
\end{keyword}

\end{frontmatter}

\section{Introduction}
\label{sec:intro}
Shear-deformable thin-structural theories such as the Timoshenko beam and
Reissner-Mindlin plate theories are widely used throughout engineering practice to simulate the mechanical response of structures with planar dimensions far greater than their thickness. The shear deformable theories' popularity over the classical thin-structural theories, Euler-Bernoulli beam and Kirchhoff-Love plate theory, is primarily due to the following two factors:

\begin{itemize}
\item They capture the shear-deformable behavior
inherent to thicker structures~\cite{arnold_asymptotic_1996,dauge_plates_2004}.
Capturing this behavior is necessary for simulating modern engineering
structures constructed from e.g.  functionally graded materials~\cite{yin_isogeometric_2014}.
\item They are second-order PDEs giving rise to weak formulations with
$H^1(\Omega)$ regularity, whereas the classical theories are fourth-order PDEs with weak formulations demanding $H^2(\Omega)$ regularity.  The
difficulty in the construction of an efficient $H^2(\Omega)$-conforming finite element method (FEM) is well-known. In contrast, $H^1(\Omega)$-conforming FEMs are straightforward.
\end{itemize}

Unfortunately, it is also the case that na{\"i}vely constructed low-order
polynomial $H^1(\Omega)$-conforming numerical methods suffer from
shear-locking.  Shear-locking is a numerical issue caused by the inability of a numerical method to represent the Kirchhoff limit as the plate thickness
parameter tends to zero.  This usually results in a nonconvergent numerical
method, or at best, very poor convergence.

The majority of solutions to the shear-locking issue in the finite element
literature resort to a mixed variational method where the transverse shear
stress is treated as an independent variational quantity in the weak form. There is a huge amount of mathematical and engineering literature related to constructing plate elements for FEM analysis, which we cannot comprehensively review here. Particularly popular examples of this approach include the
Mixed Interpolation of Tensorial Components (MITC) or Assumed Natural Strain
(ANS) approach~\cite{bathe_mitc7_1989,lyly_stable_1993,bathe_four-node_1985,duran_mixed_1992},
and the Discrete Shear Gap (DSG) method~\cite{bletzinger_unified_2000,donning_meshless_1998,kanok-nukulchai_elimination_2001}.

One desirable aspect of both MITC and DSG is that even though they can be
mathematically analyzed using (and are based on) an underlying mixed formulation, the final systems of equations are expressed in terms of the primal unknowns of the standard Reissner-Mindlin problem only.  This is achieved by the use of an operator that reduces or projects the shear stresses expressed in terms of the primal variables onto an underlying mixed finite element space. This ``unlocks'' the formulation.

Given the success of using mixed variational methods in constructing
shear-locking free finite elements, it should be no surprise that many authors have taken this route to construct convergent methods for more modern numerical techniques, such as in Isogeometric Analysis (IGA)~\cite{beirao_da_veiga_isogeometric_2012,echter_hierarchic_2013,beirao_da_veiga_avoiding_2012}
and meshfree methods~\cite{hale_locking-free_2012,tiago_eliminating_2007,kanok-nukulchai_elimination_2001,donning_meshless_1998,wang_locking-free_2004,cho_analysis_2001}.

The work described in this paper is a continuation of the line of research
presented in the PhD thesis of Hale~\cite{hale:2013:MMSDBP} on developing
meshfree methods for shear-deformable structures using mixed formulations.
There the volume-averaged nodal pressure technique that was proposed in Refs.~\cite{ortiz:2010:MEM,ortiz:2011:MEI} for the Stokes and
nearly-incompressible elasticity problems and later generalized as the volume-averaged nodal projection (\texttt{VANP}) method~\cite{ortiz:VAN:2015}, was adapted to the Reissner-Mindlin problem. In the \texttt{VANP} approach, bubble-like enrichment is used to ensure inf--sup stability~\cite{brezzi_existence_1974} mimicking the MINI element~\cite{arnold_stable_1984} and the volume-averaging procedure is closely related to the average-nodal strain finite element formulation proposed in Ref.~\cite{lamichhane_hu-washizu_2009}. The adaptation of the \texttt{VANP} approach from the Stokes problem to the Reissner-Mindlin in Ref.~\cite{hale:2013:MMSDBP} was achieved using a stabilized mixed variational formulation developed in Ref.~\cite{arnold_new_1993}. Using this stabilization it is possible to bypass the coercivity on the kernel condition in the LBB theorem~\cite{brezzi_existence_1974}, opening up the possibility of using inf--sup stable designs for the Stokes problem (e.g. the MINI element~\cite{arnold_stable_1984} or the \texttt{VANP} operator~\cite{ortiz:2010:MEM,ortiz:2011:MEI}) almost directly for the Reissner-Mindlin problem. This stabilization comes at the expense of the introduction of a stabilization constant. A more detailed analysis in Refs.~\cite{boffi_analysis_1997,lovadina_new_1996} shows that in a finite element
context it is possible to quite precisely relate this constant to the element
size and obtain good convergence rates. Unfortunately, numerical experiments to choose a good scheme for the stabilization constant in the meshfree context of Ref.~\cite{hale:2013:MMSDBP} were less successful.

In this paper, we develop a new meshfree scheme for the Reissner-Mindlin plate model with many of the best aspects of the formulation in Ref.~\cite{hale:2013:MMSDBP}, but with none of its drawbacks, such as reliance on a stabilization scheme with an a priori unknown constant. The scheme uses linear maximum-entropy basis functions for field variables approximation and is built variationally on a two-field potential energy functional wherein the shear strain, written in terms of primitive variables (i.e., rotations and transverse displacement), is computed via a volume-averaged nodal projection operator that is constructed from the Kirchhoff constraint of the three-field mixed weak form, which is an idea adapted from the \texttt{VANP} formulation of Ref.~\cite{ortiz:VAN:2015} and that leads to a symmetric stiffness matrix. The meshfree approximation is constructed over a set of scattered nodes that are obtained from an integration mesh of three-node triangles on which the meshfree stiffness matrix and nodal force vector are numerically integrated. We use recent advances in integration techniques for meshfree methods, e.g. the work of Duan et. al~\cite{duan:2012:SOI}, to ensure efficient and accurate integration of the weak form. Bubble-like enrichment of the rotation degrees of freedom is added to ensure inf--sup stability, similarly to the recent finite element of Song and Niu~\cite{Song:2016:MFERM} and many others. No further stabilization  is required to ensure the stability of the discrete problem. The final system of equations is expressed in terms of the primal unknowns only, an improvement over the work of Hale~\cite{hale_locking-free_2012}. Our numerical experiments show that the proposed method is optimally convergent for a wide range of thicknesses (shear-locking-free).

The remainder of the paper is given as follows. Section~\ref{sec:notation}
provides a summary of the notation used in this paper. The main ingredients for the computation of the maximum-entropy basis functions are given in
Section~\ref{sec:maxent}. In Section~\ref{sec:goveqns}, the classical
three-field formulation for the Reissner-Mindlin plate model is summarized
along with its discretization using meshfree basis functions. The~\texttt{VANP} method for the Reissner-Mindlin plate model is developed in Section~\ref{sec:vanp_method}. Section~\ref{sec:numexamples} presents some numerical examples that are solved using the proposed \texttt{VANP} approach. We end with some concluding remarks in Section~\ref{sec:conclusions}.

\section{Notation}
\label{sec:notation}
The following is a summary of the main notation used in this paper. Slanted bold symbols such as $\vm{v}$ are used to represent vectors and tensors. In particular, the following notation is used to represent vectors in components form: $\vm{v}=(v_1,\ldots, v_n)$ in an $n$-dimensional space and $\vm{v}=(v_x,v_y)$ in the two-dimensional Cartesian coordinate system.

Lowercase (nonbold) upright symbols are used to represent row and column vectors. Their entries are written between square brackets. For instance, $\mat{r}=[r_1 \,\, \cdots \,\, r_n]$ is a row vector and $\mat{c}=[c_1 \,\, \cdots \,\, c_n]^\transpose$ is a column vector.

Uppercase (nonbold) upright symbols are used to represent matrices. Their entries are written between square brackets. An example of a matrix representation is given as follows:
\begin{equation*}
\mat{M} = \smat{m_{11} & m_{12} & \cdots & m_{1n} \\ m_{21} & m_{22} & \cdots & m_{2n} \\ \vdots & \vdots & \ddots & \vdots \\ m_{n1} & m_{n2} & \cdots & m_{nn}}.
\end{equation*}

The gradient operator is denoted as $\bsym{\nabla}$ and the trace operator as $\trace(\cdot)$.

\section{Maximum-entropy basis functions}
\label{sec:maxent}

Consider a convex domain represented by a set of $n$
scattered nodes and a prior (weight) function $w_{a}(\vm{x})$
associated with each node $a$.  The
approximation for a scalar-valued function $u(\vx)$ is given in the form:
\begin{equation}\label{eq:trial}
u(\vm{x})=\sum_{a=1}^m\phi_a(\vm{x}) u_a,
\end{equation}
where $u_a$ are nodal coefficients.
On using the Shannon--Jaynes (or relative) entropy
functional, the maxent basis functions $\{\phi_{a}(\vm{x})\geq 0\}_{a=1}^m$
are obtained via the solution of the following convex optimization
problem~\cite{sukumar:2007:OAC}:
\begin{problem}\label{pr:maxent_problem}
\begin{align*}
\min_{\bsym{\phi} \in {\Re}_{+}^{m}} \sum_{a=1}^{m} &
\phi_{a}(\vm{x})\ln\left(\frac{\phi_{a}(\vm{x})}{w_{a}(\vm{x})}\right)
\\ \intertext{subject to the linear reproducing conditions:}
\sum_{a=1}^{m} \phi_{a}(\vm{x}) = 1, & \quad \sum_{a=1}^{m}
\phi_{a}(\vm{x})\,\vm{c}_a  = \vm{0}.
\end{align*}
\end{problem}
In~\pref{pr:maxent_problem}, $\vm{c}_a=\vx_a - \vx$ are shifted nodal coordinates and $\Re_{+}^{m}$ is the nonnegative orthant.

In this paper, we use as the prior function the Gaussian radial basis function given by~\cite{arroyo:2006:LME}
\begin{equation}\label{priors}
w_{a}(\vm{x}) = \exp\left(-\frac{\gamma}{h_{a}^2}\|\vm{c}_a\|^2\right),
\end{equation}
where $\gamma$ is a parameter that controls the support size of the basis function and $h_a$ is a characteristic nodal spacing associated with node $a$.

On using the method of Lagrange multipliers, the
solution to~\pref{pr:maxent_problem} is given
by~\cite{sukumar:2007:OAC}
\begin{equation}\label{eq:maxent_bfun}
\phi_{a}(\vm{x},\bsym{\lambda}) = \frac{Z_a(\vm{x},\bsym{\lambda}(\vm{x}))}{\sum_b Z_b(\vm{x},\bsym{\lambda}(\vm{x}))}, \quad Z_a = w_{a}(\vm{x}) \exp(- \bsym{\lambda}(\vm{x})\cdot
\vm{c}_a(\vm{x})),
\end{equation}
where the Lagrange multiplier vector $\bsym{\lambda}(\vm{x})$ is
obtained as the minimizer of the following dual optimization problem
($\vm{x}$ is fixed):
\begin{problem}\label{pr:maxent_dual}
\begin{equation*}
\bsym{\lambda}^*(\vm{x}) = \arg \min_{\bsym{\lambda}
\in {\Re}^{d}} \ln Z(\vm{x},\bsym{\lambda}).
\end{equation*}
\end{problem}

\pref{pr:maxent_dual} leads to a system of two nonlinear equations given by
\begin{equation}\label{eq:maxent_nlsys}
\vm{f}(\boldsymbol{\lambda})=\bsym{\nabla}_{\bsym{\lambda}} \ln
Z(\bsym{\lambda})=-\sum_{a}^{n}
\phi_{a}(\vm{x})\tilde{\vm{x}}_a=\vm{0},
\end{equation}
where $\bsym{\nabla}_{\bsym{\lambda}}$ stands for the gradient with respect to $\bsym{\lambda}$. Once the converged solution for the Lagrange multiplier vector $\bsym{\lambda}^*$ is found through~\eref{eq:maxent_nlsys}, the basis functions $\phi_a(\vm{x})$ are obtained by setting
$\bsym{\lambda}=\bsym{\lambda}^*$ in~\eref{eq:maxent_bfun}.

Finally, the gradient of the basis function is~\cite{arroyo:2006:LME}:
\begin{equation}\label{maxent_bfgrad}
\bsym{\nabla}\phi_{a}(\vm{x}) = \phi_{a}(\vm{x},\bsym{\lambda}^*)\left(\vm{J}(\vm{x},\bsym{\lambda}^*)\right)^{-1}\vm{c}_a(\vm{x}),
\end{equation}
where
\begin{equation}
\vm{J}(\vm{x},\bsym{\lambda})=
\sum_{a=1}^m\phi_a(\vm{x},\bsym{\lambda})\,\vm{c}_a(\vm{x})\otimes\vm{c}_a(\vm{x})-
\vm{r}(\vm{x},\bsym{\lambda})\otimes\vm{r}(\vm{x},\bsym{\lambda}),\quad
\vm{r}(\vm{x},\bsym{\lambda})= - \sum_{a=1}^m\phi_a(\vm{x},\bsym{\lambda})\,\vm{c}_a(\vm{x}).
\end{equation}

\section{Governing equations for the three-field formulation}
\label{sec:goveqns}

The method that is proposed in the next section relies on the classical three-field Reissner-Mindlin plate problem formulation. It is instructive to review this formulation as many of its aspects are preserved in the new method. Therefore, in this section we provide a summary of the classical three-field formulation for the Reissner-Mindlin plate model and its discretization procedure using the maxent basis functions.

\subsection{Strong form}
\label{sec:strongform}

\begin{figure}[htbp]
  \centering
  \epsfig{file = 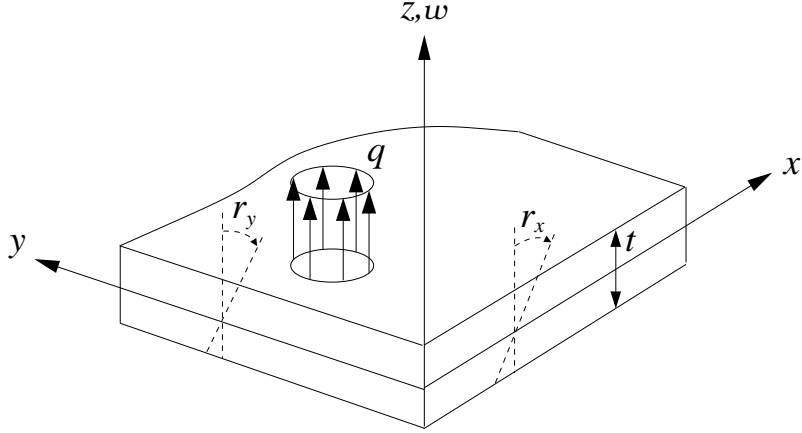, width = 0.7\textwidth}
  \caption{Schematic representation of the plate.}
  \label{fig:platedef}
\end{figure}

Consider the midplane of an elastic plate of uniform thickness $t$ that occupies the open domain $\Omega \subset \Re^2$ and is bounded by the one-dimensional surface $\Gamma$. The coordinates of a point in this domain are denoted by $\vm{x}=(x,y)$. The rotations of fibers normal to the midplane, the transverse displacement of the midplane, and the scaled transverse shear stresses are denoted by $\vm{r}(\vm{x})=(r_x,r_y)$, $w(\vm{x})$, and $\vm{s}(\vm{x})=(s_x,s_y)$, respectively. A transverse load $q(\vm{x}) \in L^2(\Omega)$ is also acting on the plate. A schematic representation of the plate is depicted in~\fref{fig:platedef}.

The boundary of the plate is assumed to be entirely subjected to the essential (Dirichlet) boundary conditions $\hat{\vm{r}}(\vm{x}):\Gamma_D\rightarrow \Re^2$ and $\hat{w}(\vm{x}):\Gamma_D\rightarrow \Re$, which implies that $\Gamma=\Gamma_D$.

The boundary-value problem for this Reissner-Mindlin plate configuration reads~\cite{boffi:2008:MFE}:
\begin{problem}\label{pr:strong_form}
Find $\vm{r}(\vm{x}):\Omega \rightarrow \Re^2$, $w(\vm{x}): \Omega \rightarrow \Re$ and $\vm{s}(\vm{x}):\Omega \rightarrow \Re^2$
such that
\begin{align*}
-\bsym{\nabla}\cdot\left(\bsym{\mathcal{C}}\bsym{\varepsilon}(\vm{r})\right)-\vm{s} & =\vm{0} \quad \forall \vm{x} \in \Omega,\\
-\bsym{\nabla}\cdot\vm{s}-q &= 0 \quad \forall \vm{x} \in \Omega, \\
(\bsym{\nabla} w - \vm{r}) - \frac{1}{\lambda t^{-2}}\vm{s} &= \vm{0} \quad \forall \vm{x} \in \Omega,\\
\vm{r}=\hat{\vm{r}}, \quad w &= \hat{w} \quad \forall \vm{x} \in \Gamma_D.
\end{align*}
\end{problem}
In~\pref{pr:strong_form}, $\bsym{\varepsilon}(\vm{r})=\frac{1}{2}\left(\bsym{\nabla}\vm{r}+(\bsym{\nabla}\vm{r})^\transpose\right)$ is the strain tensor, $\bsym{\mathcal{C}}=\frac{E_\mathrm{Y}}{12(1-\nu^2)}\left((1-\nu)\bsym{\varepsilon}+\nu\trace(\bsym{\varepsilon})\vm{I}\right)$ is the tensor of bending moduli, and $\lambda=\kappa\frac{E_\mathrm{Y}}{2(1+\nu)}$, where $\kappa=5/6$ is the shear correction factor; $E_\mathrm{Y}$ and $\nu$ are the Young's modulus and the Poisson's ratio of the plate material, respectively.

\subsection{Three-field mixed weak form}
\label{sec:weakform}

Let the spaces
\begin{align*}
R :=& \left\{\vm{r}:\vm{r} \in [H^{1}(\Omega)]^2, \ \vm{r} = \hat{\vm{r}} \ \textrm{on} \ \Gamma_D \right\}\\
R_0 :=& \left\{\delta\vm{r}:\delta\vm{r} \in [H^{1}(\Omega)]^2, \ \delta\vm{r} = \vm{0} \ \textrm{on} \ \Gamma_D \right\}\\
\intertext{be the trial and virtual spaces for the rotation field, respectively,}
W :=& \left\{w:w \in H^{1}(\Omega), \ w = \hat{w} \ \textrm{on} \ \Gamma_D \right\}\\
W_0 :=& \left\{\delta w:\delta w \in H^{1}(\Omega), \ \delta w = 0 \ \textrm{on} \ \Gamma_D \right\}\\
\intertext{be the trial and virtual spaces for the transverse displacement field, respectively, and}
S :=& \left\{\vm{z}:\vm{z} \in [L^{2}(\Omega)]^2 \right\}
\end{align*}
be the space for the trial and virtual scaled transverse shear stresses.

On using the preceding space definitions, the three-field mixed weak form reads~\cite{falk:LFFEMRM:2000,boffi:2008:MFE}:

\begin{problem}\label{pr:weakform_RM_plate}
Find $(\vm{r},w,\vm{s}) \in R\times W\times S$ such that
\begin{align*}
\int_{\Omega} (\upvarepsilon(\delta\vm{r}))^\transpose\mat{C}\,\upvarepsilon(\vm{r})\,\diffx -\int_{\Omega} \vm{s}\cdot\delta\vm{r}\,\diffx &=0 \quad \forall \delta\vm{r} \in R_0,\\
\int_{\Omega} \vm{s}\cdot\bsym{\nabla}\delta w\,\diffx - \int_{\Omega} q\,\delta w\,\diffx &= 0 \quad \forall \delta w \in W_0,\\
\int_{\Omega}\left(\left(\bsym{\nabla} w - \vm{r}\right) - \frac{\vm{s}}{\lambda t^{-2}}\right)\cdot\delta\vm{s}\,\diffx  &=0 \quad \forall \delta\vm{s} \in S.
\end{align*}
\end{problem}
In~\pref{pr:weakform_RM_plate}, $\upvarepsilon(\vm{r})=[\varepsilon_{xx} \,\, \varepsilon_{yy} \,\, 2\varepsilon_{xy}]^\transpose$ is the strain tensor in Voigt notation, and $\mat{C}$ is the Voigt representation of the tensor of bending moduli and is given by
\begin{equation}\label{eq:bendingmoduli}
\mat{C} =
\frac{E_Y}{12(1-\nu^2)}\smat{1 & \nu & 0\\ \nu & 1 & 0\\ 0 & 0 & \frac{1-\nu}{2}}.
\end{equation}

\subsection{Discrete equations using maxent basis functions}
\label{sec:weakform_threefield_discrete}

The discrete version of~\pref{pr:weakform_RM_plate} is obtained by approximating the field variables using the maxent basis functions over a set of scattered nodes that discretely represent the continuous plate. Due to the nonpolynomial nature of the maxent basis functions, the weak form integrals cannot be computed exactly. Thus, numerical quadrature is used to evaluate them. For this purpose, we construct a finite element mesh whose elements are used to define integration points and its nodes to discretize the field variables. The whole procedure can be thought as a finite element method with basis functions having a radial support. The support is controlled by the maxent parameters and its size is typically larger than the support of a finite element basis function. The construction of maxent basis functions depends only on the nodal coordinates (see Section~\ref{sec:maxent}) for which they are regarded as ``meshfree." An advantage of using meshfree basis functions is that since the element is not involved in the computation of them, the resulting method is less sensitive to mesh distortion than the finite element method.

\begin{figure}[htbp]
  \centering
  \epsfig{file = 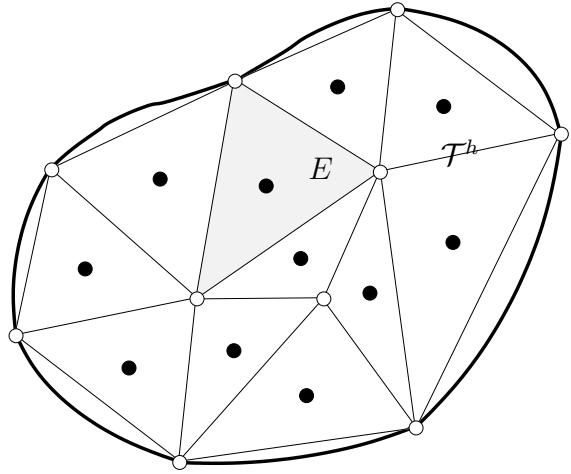, width = 0.5\textwidth}
  \caption{Schematic representation of the domain partition into cells and nodes for the discretization of the three-field mixed weak form. The shaded triangle is a typical cell of the partition. The white circles represent the standard node set $\mathcal{N}^s$ and the black ones the barycenter node set $\mathcal{N}^{b}$.}
  \label{fig:nodesets}
\end{figure}

For the construction of the integration mesh, we follow the standard practice in finite elements. That is, we consider a mesh of so-called mixed finite elements that would produce convergent and stable finite element solutions for the three-field mixed weak form that models the Reissner-Mindlin plate problem. Several mixed finite elements are available in the finite element literature (see for instance Ref.~\cite{zienkiewicz:2005:FEM2}). In this paper, we construct the integration mesh inspired by the recent work of Song and Niu~\cite{Song:2016:MFERM}, as follows: let the domain be partitioned into disjoint nonoverlapping three-node triangular cells. We denote an integration cell as $E$. The partition formed by these cells is denoted as $\mathcal{T}^h$, where $h$ is the maximum cell diameter among the cells in the partition. The standard set of nodes, denoted by $\mathcal{N}^s$, is formed by the vertices of the triangular mesh. In addition to the standard node set, we define a barycenter node set as $\mathcal{N}^{b}$ with nodes located at the barycenter of each cell. So, the enhanced node set is defined as $\mathcal{N}^{+}=\mathcal{N}^s \cup \mathcal{N}^b$. The degrees of freedom associated with this partition is summarized as follows:
\begin{itemize}
\item each node in the standard node set carries two rotations, one transverse displacement, and two transverse shear stresses.
\item each node in the barycenter node set carries two rotations.
\end{itemize}
\fref{fig:nodesets} presents a schematic representation of the integration mesh.

The partition $\mathcal{T}^h$ is constructed using a triangular mesh generator and the location of quadrature points is computed based on this partition. The enhanced node set $\mathcal{N}^{+}$ is constructed when needed by adding the barycenter node set $\mathcal{N}^{b}$ to the standard node set $\mathcal{N}^{s}$. This poses no problem or additional complexity in the method since the absence of the finite element structure in the computation of meshfree basis functions permits the addition of nodes to the mesh very easily.

\begin{figure}[!tbh]
  \centering
  \mbox{
  \subfigure[]{\label{fig:cellcont_a}
  \epsfig{file = 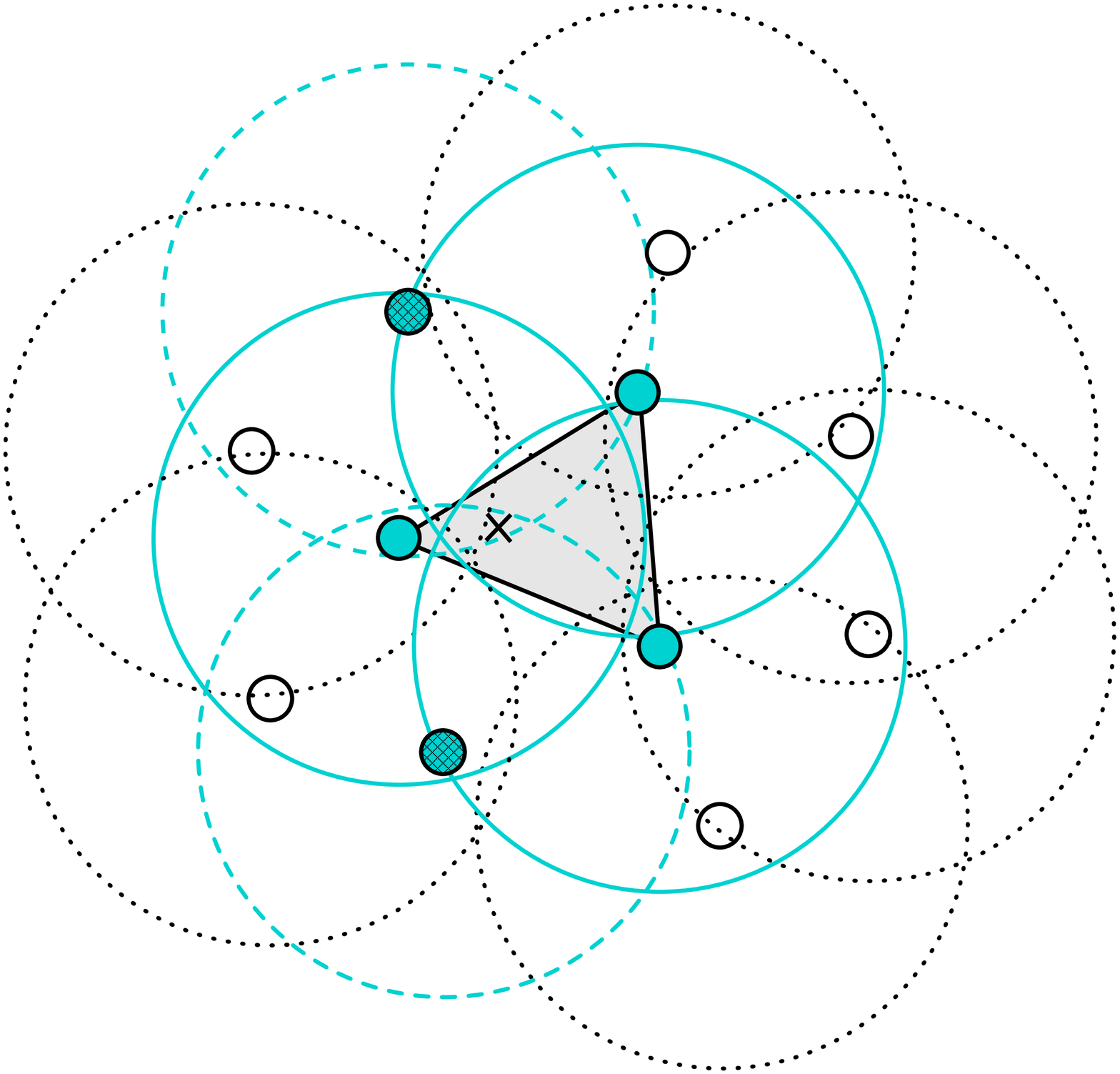, width = 0.48\textwidth}}
  \subfigure[]{\label{fig:cellcont_b} \epsfig{file = 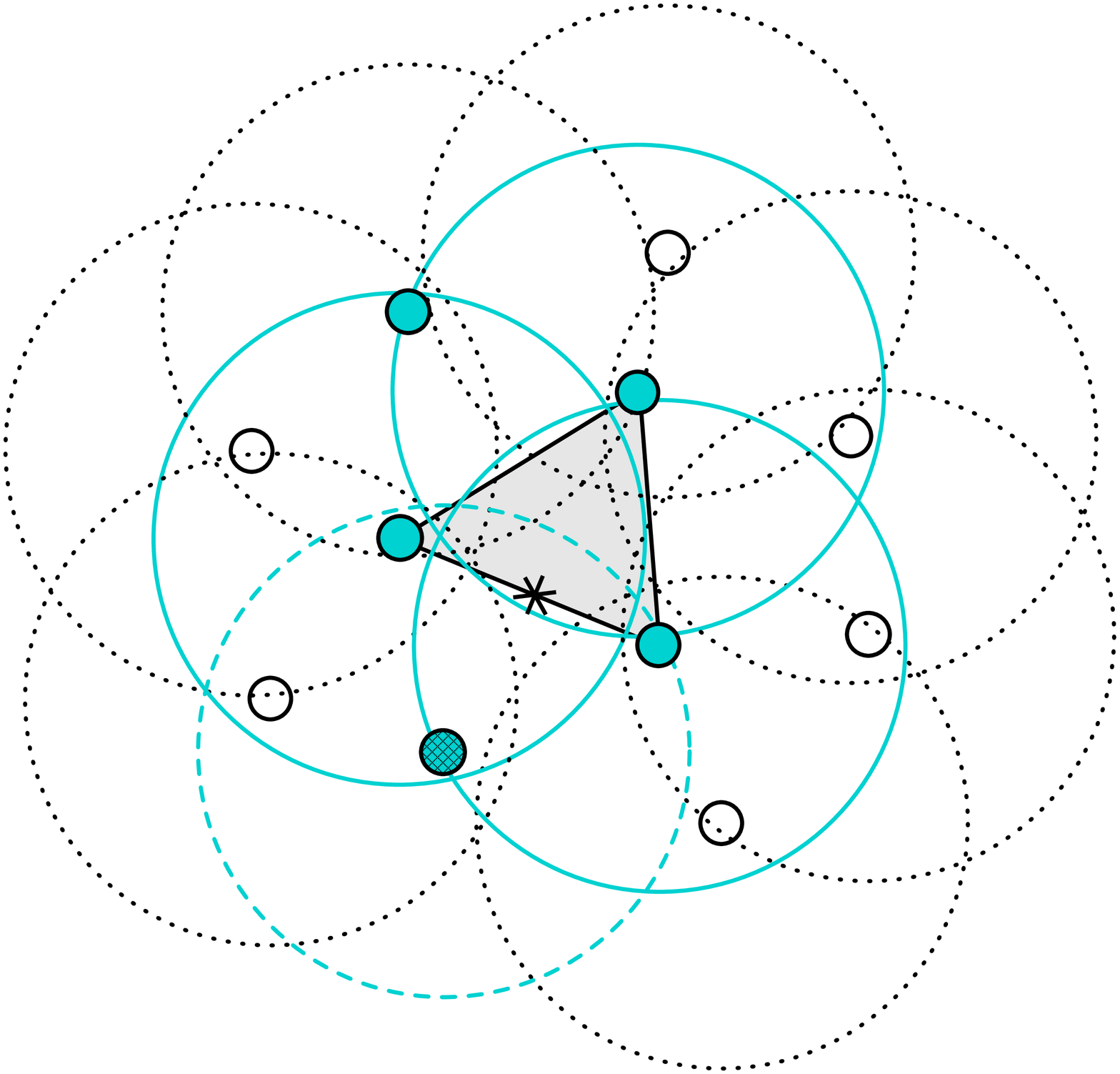, width = 0.48\textwidth}}
  }
  \caption{Graphical explanation of the nodal contribution concept for the evaluation of nodal basis functions at (a) an integration point (depicted as $\times$) located in the interior of the cell and (b) an integration point (depicted as $\ast$) located on an edge of the cell. The support of nodal basis functions are represented by circles centered at nodes. The indices of the nodes whose basis functions have a nonzero value at the integration point (i.e., the indices of the nodes whose associated circles contain the integration point) define the nodal contribution. In this example, the nodal contribution contains the indices of the nodes that define the integration cell (filled nodes) and some nodes that lie outside the integration cell (filled dashed nodes).}\label{fig:cellcont}
\end{figure}

In the process of evaluating the discrete weak form, maxent basis functions need to be computed at integration points. To make clear the implications of the numerical integration procedure using meshfree basis functions, the concept of nodal contribution is introduced as follows: the nodal contribution at a given integration point with coordinates $\vm{x}$ is defined as the indices of the nodes whose
basis functions have a nonzero value at $\vm{x}$. It should be noted that due to the radial support of the maxent basis functions, their evaluation at an integration point is likely to have a nodal contribution stemming not only from the nodes that define the integration cell, but also from nodes located outside the integration cell. A graphical explanation of the nodal contribution at an integration point located in the interior of the cell and an integration point located on an edge of the cell is provided in~\fref{fig:cellcont}, where the support of nodal basis functions are represented by circles centered at nodes. The circles drawn with continuous line and centered at filled nodes represent basis functions of nodes defining the integration cell and having a nonzero value at the integration point. Hence, the indices of filled nodes are part of the nodal contribution. The circles drawn with dashed line and centered at filled dashed nodes represent basis functions of nodes lying outside the integration cell and having a nonzero value at the integration point. Thus, the indices of filled dashed nodes are also part of the nodal contribution. The circles drawn with dotted line and centered at empty nodes represent nodal basis functions having a zero value at the integration point. The indices of empty nodes are then not part of the nodal contribution.

The nodal contribution concept is not restricted to the numerical integration procedure only, it is in general applicable to any evaluation of basis functions within $\Omega$ or on $\Gamma$.

On using the maxent basis functions, the discrete trial and virtual field variables are obtained as follows:
\begin{subequations}\label{eq:fieldapprox}
\begin{align}
\vm{r}^h(\vm{x})=\sum_{a=1}^{nenh}\phi_a(\vm{x}) \vm{r}_a, & \quad \delta\vm{r}^h(\vm{x})=\sum_{b=1}^{nenh}\phi_b(\vm{x}) \delta\vm{r}_b, \label{eq:fieldapprox_a} \\
w^h(\vm{x})=\sum_{a=1}^{nstd}\phi_a(\vm{x}) w_a, & \quad \delta w^h(\vm{x})=\sum_{b=1}^{nstd}\phi_b(\vm{x}) \delta w_b,\label{eq:fieldapprox_b}\\
\vm{s}^h(\vm{x})=\sum_{a=1}^{nstd}\phi_a(\vm{x}) \vm{s}_a, & \quad \delta\vm{s}^h(\vm{x})=\sum_{b=1}^{nstd}\phi_b(\vm{x}) \delta\vm{s}_b,\label{eq:fieldapprox_c}
\end{align}
\end{subequations}
where $nenh$ and $nstd$ are the number of nodes that define the nodal contributions at $\vm{x}$ in the enhanced node set ($\mathcal{N}^{+}$) and the standard node set ($\mathcal{N}^s$), respectively, and $\vm{r}_a=\vm{r}(\vm{x}_a)$, $w_a=w(\vm{x}_a)$ and $\vm{s}_a=\vm{s}(\vm{x}_a)$ are the unknown nodal coefficients.

Thus, the discrete three-field mixed weak form at the integration cell level reads:
\begin{problem}\label{pr:disc_weakform_RM_plate}
Find $(\vm{r}^h,w^h,\vm{s}^h) \in (R^h \subset R) \times (W^h \subset W) \times (S^h \subset S)$ such that
\begin{align*}
\int_E (\upvarepsilon^h(\delta\vm{r}^h))^\transpose\mat{C}\,\upvarepsilon^h(\vm{r}^h)\,\diffx -\int_E \vm{s}^h\cdot\delta\vm{r}^h\,\diffx &=0 \quad \forall \delta\vm{r}^h \in R_0^h \subset R_0,\\
\int_E \vm{s}^h\cdot\bsym{\nabla}\delta w^h\,\diffx - \int_E q\,\delta w^h\,\diffx &= 0 \quad \forall \delta w^h \in W_0^h \subset W_0,\\
\int_E\left(\left(\bsym{\nabla} w^h - \vm{r}^h\right)-\frac{\vm{s}^h}{\lambda t^{-2}}\right)\cdot\delta\vm{s}^h\,\diffx  &=0 \quad \forall \delta\vm{s}^h \in S_0^h \subset S_0.
\end{align*}
\end{problem}

\section{The volume-averaged nodal projection method}
\label{sec:vanp_method}

In analogy to the \texttt{VANP} method for nearly-incompressible elasticity~\cite{ortiz:VAN:2015}, this method when applied to the Reissner-Mindlin plate model allows the elimination of the scaled shear stresses from the analysis, which leads to a method written in terms of the primitive variables $\vm{r}$ and $w$. In this section, the \texttt{VANP} method for the Reissner-Mindlin plate model is formulated.

\subsection{Projection operator}
\label{sec:vanp_operator}

Consider the discrete two-field scaled variational formulation for the Reissner-Mindlin plate model~\cite{chinosi:NAMFERM:1995,boffi:2008:MFE}:

\begin{problem}\label{pr:discr_potential_RM}
The field variables $(\vm{r}^h,w^h) \in (R^h \subset R) \times (W^h \subset W)$ can be found as the unique minimum point of the following potential energy functional:
\begin{equation*}
\Psi(\vm{r}^h,w^h) =  \inf_{\substack{\vm{r}^h,w^h}} \frac{1}{2}\int_E (\upvarepsilon^h(\vm{r}^h))^\transpose\mat{C}\upvarepsilon^h(\vm{r}^h)\,\diffx +\frac{\lambda t^{-2}}{2}\int_E \left(\bsym{\nabla} w^h - \vm{r}^h\right)^\transpose\left(\bsym{\nabla} w^h - \vm{r}^h\right)\,\diffx - \int_E q w^h\,\diffx.
\end{equation*}
\end{problem}

\pref{pr:discr_potential_RM} requires the minimizing pair ${\vm{r}^h,w^h}$ to satisfy the Kirchhoff constraint $\bsym{\nabla} w^h - \vm{r}^h = \vm{0}$ as the thickness of the plate becomes very small, which at the element level is a severe condition that leads to shear-locking. This issue is purely numerical and manifests itself as the impossibility for a displacement-based formulation (i.e., a formulation in primitive variables $\vm{r}$ and $w$) to undergo deformations as the thickness of the plate becomes too small.

As a remedy for shear-locking, the following modified version of~\pref{pr:discr_potential_RM} is considered:

\begin{problem}\label{pr:discr_modified_potential_RM}
The field variables $(\vm{r}^h,w^h) \in (R^h \subset R) \times (W^h \subset W)$ can be found as the unique minimum point of the following modified potential energy functional:
\begin{equation*}
\overline{\Psi}(\vm{r}^h,w^h) =  \inf_{\substack{\vm{r}^h,w^h}} \frac{1}{2}\int_E (\upvarepsilon^h(\vm{r}^h))^\transpose\mat{C}\upvarepsilon^h(\vm{r}^h)\,\diffx +\frac{\lambda t^{-2}}{2}\int_E (\overline{\bsym{\nabla} w^h - \vm{r}^h})^\transpose(\overline{\bsym{\nabla} w^h - \vm{r}^h})\,\diffx - \int_E q w^h\,\diffx.
\end{equation*}
\end{problem}
The ``bar'' symbol in~\pref{pr:discr_modified_potential_RM} is intended to define a modified shear strain that alleviates shear-locking.

On taking the first variation of the modified potential energy functional in~\pref{pr:discr_modified_potential_RM}, leads to the following discrete modified two-field weak form for the Reissner-Mindlin model:
\begin{problem}\label{pr:disc_weakform_VANP_RM_plate}
Find $(\vm{r}^h,w^h) \in (R^h \subset R) \times (W^h \subset W)$ such that
\begin{align*}
&\lambda t^{-2}\int_E(\overline{\bsym{\nabla}\delta w^h})^\transpose\,\overline{\bsym{\nabla} w^h}\,\diffx
- \lambda t^{-2}\int_E(\overline{\bsym{\nabla}\delta w^h})^\transpose\,\overline{\vm{r}^h}\,\diffx - \int_E \delta w^h\,q\,\diffx = 0\\
&\forall \delta w^h \in W_0^h \subset W_0,\\
&- \lambda t^{-2}\int_E(\overline{\delta \vm{r}^h})^\transpose\,\overline{\bsym{\nabla} w^h}\,\diffx
+ \int_E (\upvarepsilon^h(\delta\vm{r}^h))^\transpose\mat{C}\,\upvarepsilon^h(\vm{r}^h)\,\diffx
+ \lambda t^{-2}\int_E(\overline{\delta \vm{r}^h})^\transpose\,\overline{\vm{r}^h}\,\diffx = 0\\
&\forall \delta\vm{r}^h \in R_0^h \subset R_0.
\end{align*}
\end{problem}

In contrast to the system that gives form to the local patch projection method~\cite{hale:2013:MMSDBP}, \pref{pr:disc_weakform_VANP_RM_plate} is a symmetric system. This is a consequence of the modified shear strain being applied to the potential energy functional.

In order to develop the stiffness matrix from~\pref{pr:disc_weakform_VANP_RM_plate}, an appropriate construction for the ``barred'' quantities that appear therein is needed --- here ``appropriate'' means that shear-locking is precluded. An effective procedure to achieve this aim is offered by the Kirchhoff constraint given in the last equality of~\pref{pr:disc_weakform_RM_plate}. The procedure consists in rearranging the Kirchhoff constraint such that $\vm{s}^h$ can be computed in terms of the primitive variables, as follows:
\begin{equation}\label{eq:projected_shear_stress}
\vm{s}^h = \lambda t^{-2} \pi^h \left[\bsym{\nabla} w^h - \vm{r}^h\right]=\lambda t^{-2}\left(\pi^h \left[\bsym{\nabla} w^h\right]-\pi^h \left[\vm{r}^h\right]\right),
\end{equation}
where $\pi^h$ is a projection operator that adopts the form of an $L^2$ projection. By comparing the second equation in~\pref{pr:disc_weakform_VANP_RM_plate} with~\eref{eq:projected_shear_stress}, the following equalities are proposed:
\begin{equation}\label{eq:projected_strain_components}
\overline{\bsym{\nabla}w^h} = \pi^h\left[\bsym{\nabla}w^h\right], \quad \overline{\vm{r}^h} = \pi^h\left[\vm{r}^h\right],
\end{equation}
which give the definition of the ``bar'' operator as $\overline{(\,\cdot\,)} := \pi^h[\,\cdot\,]$.

We are left with the explicit expression for the projection operator. It is derived as follows: the discrete scaled transverse shear stresses given in~\eref{eq:fieldapprox_c} are replaced in the last equation in~\pref{pr:disc_weakform_RM_plate} (the Kirchhoff constraint), which after relying on the arbitrariness of nodal variations yields in nodal form
\begin{equation}\label{eq:volume_avg}
\int_E\phi_c(\vm{x})\left[\bsym{\nabla} w^h - \vm{r}^h\right]\,\diffx -\frac{1}{\lambda t^{-2}}\int_E\phi_c(\vm{x})\phi_b(\vm{x})\vm{s}_b\,\diffx=0.
\end{equation}
\alejandro{The integral that accompanies the nodal transverse shear stress on the left-hand side of~\eref{eq:volume_avg} defines a matrix $\mat{H}$ whose entries are given by
\begin{equation}\label{eq:mass_matrix}
\mat{H}_{cb}=\int_E\phi_c(\vm{x})\phi_b(\vm{x})\,\diffx.
\end{equation}

Eq.~\eref{eq:volume_avg} can be solved for the nodal scaled transverse shear stress as follows:
\begin{equation}\label{eq:nodal_shearstress}
\vm{s}_b = \lambda t^{-2}\,\mat{H}_{cb}^{-1}\int_E\phi_c(\vm{x})\left[\bsym{\nabla} w^h - \vm{r}^h\right]\,\diffx,
\end{equation}
where $\mat{H}_{cb}^{-1}$ is read as the nodal component of the inverse of the matrix $\mat{H}$. Since the maxent basis functions are positive functions, all the entries in the matrix $\mat{H}$ are nonnegative. Hence, Eq.~\eref{eq:nodal_shearstress} can be safely simplified by performing row-sum on $\mat{H}$ (i.e., by lumping) leading to the following
\textit{volume-averaged nodal scaled transverse shear stress} vector:}
\begin{equation}\label{eq:nodal_shearstress_lumped}
\vm{s}_c = \lambda t^{-2}\frac{\int_{E_c}\phi_c(\vm{x})\left[\bsym{\nabla} w^h - \vm{r}^h\right]\,\diffx}{\int_{E_c}\phi_c(\vm{x})\,\diffx},
\end{equation}
which is used to project the scaled transverse shear stress field, as follows:
\begin{equation}\label{eq:vanp_projected_shear_stress}
    \vm{s}^h = \sum_{c=1}^{nstd}\phi_c(\vm{x})\vm{s}_c=\lambda t^{-2}\sum_{c=1}^{nstd}\phi_c(\vm{x})\left\lbrace \frac{\int_{E_c}\phi_c(\vm{x})\left[\bsym{\nabla} w^h - \vm{r}^h\right]\,\diffx}{\int_{E_c}\phi_c(\vm{x})\,\diffx}\right\rbrace.
\end{equation}
The projection operator that defines the ``bar'' operator is realized by comparing~\eref{eq:projected_shear_stress} with~\eref{eq:vanp_projected_shear_stress} and is given by
\begin{equation}\label{eq:projection_operator}
    \overline{(\,\cdot\,)} = \pi^h[\,\cdot\,] = \sum_{c=1}^{nstd} \phi_c(\vm{x})\pi_c[\,\cdot\,],
\end{equation}
where $\pi_c[\,\cdot\,]$ is the \textit{volume-averaged nodal projection} (\texttt{VANP}) operator given by
\begin{equation}\label{eq:vanp_operator}
    \pi_c[\,\cdot\,] = \frac{\int_{E_c} \phi_c(\vm{x}) [\,\cdot\,] \, \diffx}{\int_{E_c} \phi_c(\vm{x}) \, \diffx}.
\end{equation}

\begin{figure}[!tbhp]
\centering
\subfigure[]{\label{fig:volumedef_a} \epsfig{file = 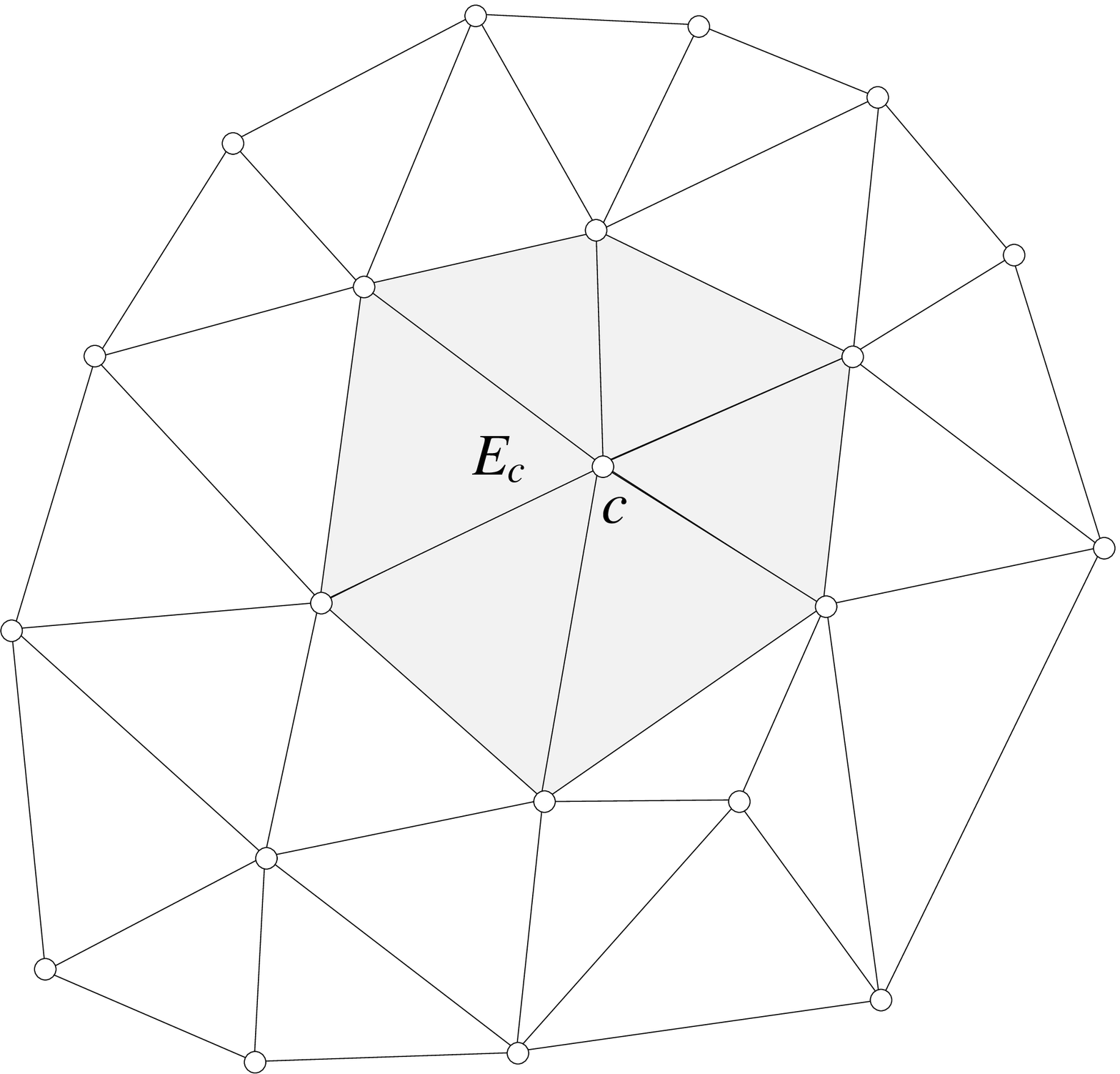, width = 0.4\textwidth}}
\subfigure[]{\label{fig:volumedef_b} \epsfig{file = 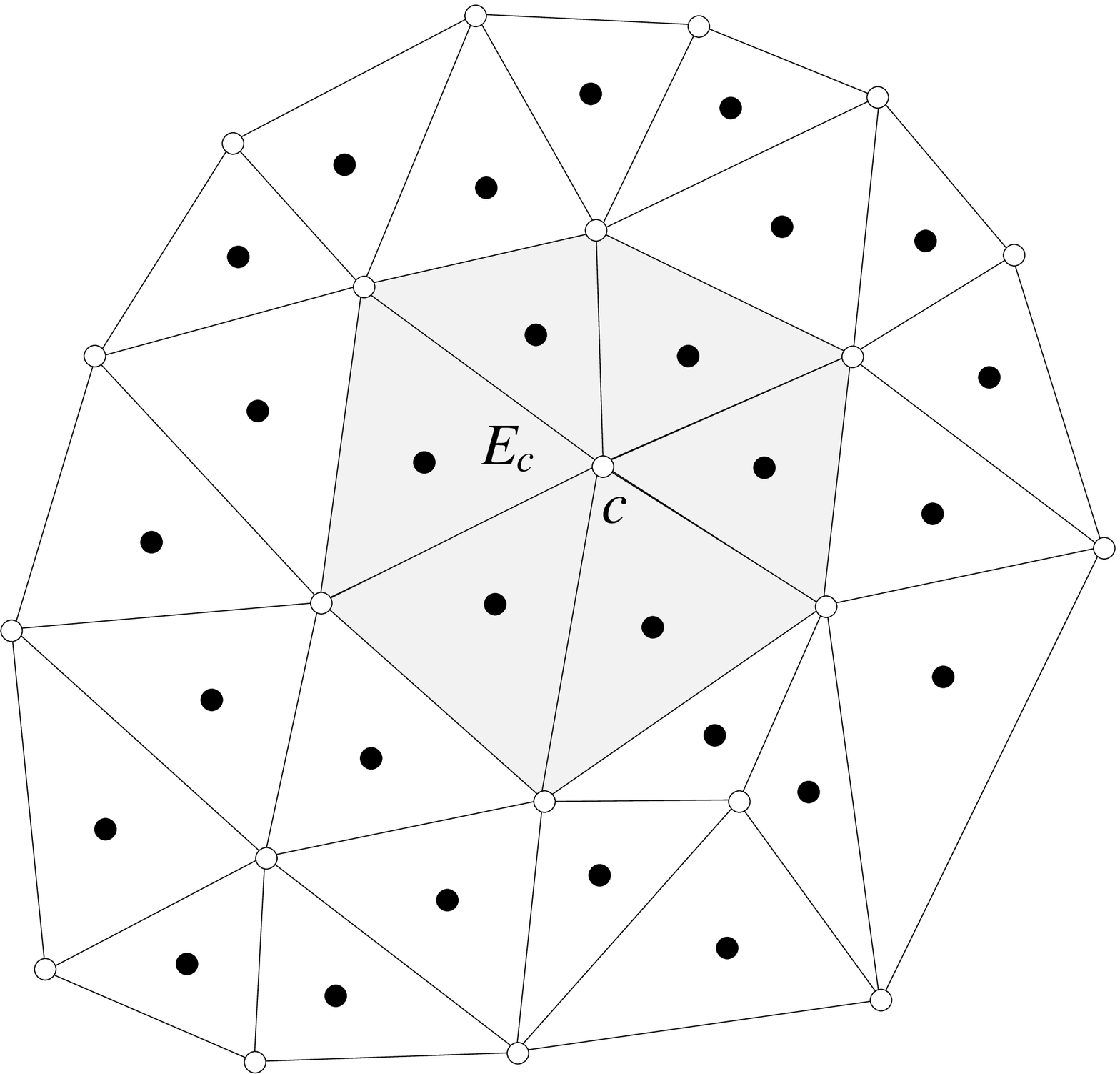, width = 0.4\textwidth}}
\caption{Definition of the representative nodal volume $E_c$ (shaded area) for the evaluation of the integrals that appear in the volume-averaged nodal projection operator. (a) nodal volume based on the standard node set $\mathcal{N}^{s}$, and (b) nodal volume based on the enhanced node set $\mathcal{N}^{+}$.}
\label{fig:volumedef}
\end{figure}

In~\eref{eq:vanp_operator}, $E_c$ is a representative nodal volume defined as the union of all the elements attached to node $c$. \fref{fig:volumedef_a} depicts the nodal volume $E_c$ when the standard node set $\mathcal{N}^{s}$ is used, and \fref{fig:volumedef_b} when the enhanced node set $\mathcal{N}^{+}$ is used.

The following precautions must be taken into account when computing the \texttt{VANP} operator:
\begin{itemize}
\item It should be noted that since $\phi_c(\vm{x})$ in~\eref{eq:vanp_operator} stems from the nodal shear stresses variations its evaluation must always be performed in $\mathcal{N}^{s}$, which requires $E_c$ to be defined as in~\fref{fig:volumedef_a}.
\item Observing the ``barred'' terms in~\pref{pr:disc_weakform_VANP_RM_plate}, it must be realized that between the square brackets in the \texttt{VANP} operator we will have either $\vm{r}^h$ or $\bsym{\nabla} w^h$. The computation of the former requires the enhanced node set and thus $E_c$ must be defined as in~\fref{fig:volumedef_b}, and the computation of the latter requires the standard node set and thus $E_c$ must be defined as in~\fref{fig:volumedef_a}.
\item The evaluation of the \texttt{VANP} operator requires numerical integration at quadrature points and thus the nodal contribution concept (see~\fref{fig:cellcont} to recall this concept) must also be considered.
\end{itemize}

\subsection{Stiffness matrix and nodal force vector}
\label{sec:vanp_stiffnessmatrix}

The stiffness matrix and nodal force vector are developed by discretizing~\pref{pr:disc_weakform_VANP_RM_plate} with the rotation and transverse displacement fields approximations given in~\eref{eq:fieldapprox_a}  and~\eref{eq:fieldapprox_b}, respectively. On using these approximations, we write
\begin{equation} \label{eq:std_strain}
\bsym{\varepsilon}^h(\vm{r}^h)=\sum_{a=1}^{nenh}\mat{B}_a(\vm{x})\vm{r}_a, \quad \bsym{\varepsilon}^h(\delta\vm{r}^h)=\sum_{b=1}^{nenh}\mat{B}_b(\vm{x})\delta\vm{r}_b, \quad
\mat{B}_a(\vm{x})=\smat{\phi_{a,x} & 0 \\ 0 & \phi_{a,y} \\ \phi_{a,y} & \phi_{a,x}},
\end{equation}
and
\begin{equation}\label{eq:shear_strain_wcomp}
\bsym{\nabla}w^h = \sum_{a=1}^{nstd}\mat{G}_a(\vm{x})w_a, \quad \bsym{\nabla}\delta w^h = \sum_{b=1}^{nstd}\mat{G}_b(\vm{x})\delta w_b, \quad \mat{G}_a(\vm{x})=\smat{\phi_{a,x} \\ \phi_{a,y}}.
\end{equation}
In addition, the discrete rotation field is conveniently rewritten as
\begin{equation}\label{eq:shear_strain_rcomp}
\vm{r}^h = \sum_{a=1}^{nenh}\mat{N}_a(\vm{x})\vm{r}_a, \quad \delta\vm{r}^h = \sum_{b=1}^{nenh}\mat{N}_b(\vm{x})\delta\vm{r}_b, \quad \mat{N}_a(\vm{x})=\smat{\phi_a & 0 \\ 0 & \phi_a}.
\end{equation}
And on using~\eref{eq:projection_operator}, \eref{eq:shear_strain_wcomp} and~\eref{eq:shear_strain_rcomp}, the following projected terms are obtained:
\begin{equation} \label{eq:projected_shear_strain_wcomp}
\overline{\bsym{\nabla}w^h} =\sum_{a=1}^{nstd}\left\lbrace\sum_{c=1}^{nstd} \phi_c(\vm{x})\pi_c\left[\mat{G}_a(\vm{x})\right]\right\rbrace w_a, \quad
\overline{\bsym{\nabla}\delta w^h} =\sum_{b=1}^{nstd}\left\lbrace\sum_{c=1}^{nstd} \phi_c(\vm{x})\pi_c\left[\mat{G}_b(\vm{x})\right]\right\rbrace\delta w_b,
\end{equation}
and
\begin{equation} \label{eq:projected_shear_strain_rcomp}
\overline{\vm{r}^h} =\sum_{a=1}^{nenh}\left\lbrace\sum_{c=1}^{nstd} \phi_c(\vm{x})\pi_c\left[\mat{N}_a(\vm{x})\right]\right\rbrace \vm{r}_a, \quad
\overline{\delta\vm{r}^h} =\sum_{b=1}^{nenh}\left\lbrace\sum_{c=1}^{nstd} \phi_c(\vm{x})\pi_c\left[\mat{N}_b(\vm{x})\right]\right\rbrace \delta\vm{r}_b.
\end{equation}

Finally, by collecting all the discrete quantities and replacing them into the modified two-field weak form (\pref{pr:disc_weakform_VANP_RM_plate}), and appealing to the arbitrariness of nodal variations, the following global system of equations is obtained after assembling the local contributions:
\begin{subequations} \label{eq:discrete_system}
\begin{equation}
\smat{\mat{K}_{ww} & -\mat{K}_{wr} \\ -\mat{K}_{wr}^\transpose & \mat{K}_{ee}+\mat{K}_{rr}}\smat{\mat{w} \\ \mat{r}}=\smat{\mat{f}_w \\ \mat{0}}, \label{eq:discrete_system_a}
\end{equation}
where $\mat{w}$ and $\mat{r}$  are the global column vectors of nodal coefficients for the transverse displacement and rotations, respectively; $\mat{K}_{ww}$, $\mat{K}_{wr}$, $\mat{K}_{ee}$ and $\mat{K}_{rr}$ are the assembled stiffness matrices, and $\mat{f}_w$ is the assembled nodal force vector.

On defining the assembly operator as $\bm{\mat{A}}$, the assembled stiffness matrices for the partition of the domain into $nel$ integration cells are obtained as
\begin{equation}
\mat{K}_{ww}=\assembly_{E=1}^{nel}\mat{K}_{ww}^E,\quad
\mat{K}_{wr}=\assembly_{E=1}^{nel}\mat{K}_{wr}^E,\quad
\mat{K}_{ee}=\assembly_{E=1}^{nel}\mat{K}_{ee}^E,\quad
\mat{K}_{rr}=\assembly_{E=1}^{nel}\mat{K}_{rr}^E,\quad
\mat{f}_w=\assembly_{E=1}^{nel}\mat{f}_w^E,
\label{eq:discrete_system_b}
\end{equation}
where the local stiffness matrices evaluated on the integration cell $E$ are
\begin{align}
\mat{K}_{ww}^E  &=  \sum_{a=1}^{nstd}\sum_{b=1}^{nstd}\left(\lambda t^{-2}\int_E\left(\sum_{c=1}^{nstd} \phi_c(\vm{x})\pi_c\left[\mat{G}_a(\vm{x})\right]\right)^\transpose\sum_{c=1}^{nstd} \phi_c(\vm{x})\pi_c\left[\mat{G}_b(\vm{x})\right]\,\diffx\right),\label{eq:discrete_system_c}\\
\mat{K}_{wr}^E  &=  \sum_{a=1}^{nstd}\sum_{b=1}^{nenh}\left(\lambda t^{-2}\int_E\left(\sum_{c=1}^{nstd} \phi_c(\vm{x})\pi_c\left[\mat{G}_a(\vm{x})\right]\right)^\transpose\sum_{c=1}^{nstd} \phi_c(\vm{x})\pi_c\left[\mat{N}_b(\vm{x})\right]\,\diffx\right),\label{eq:discrete_system_d}\\
\mat{K}_{ee}^E  &=  \sum_{a=1}^{nenh}\sum_{b=1}^{nenh}\left(\int_E\mat{B}_a^\transpose(\vm{x})\mat{C}\mat{B}_b(\vm{x})\,\diffx\right),\label{eq:discrete_system_e}\\
\mat{K}_{rr}^E  &=  \sum_{a=1}^{nenh}\sum_{b=1}^{nenh}\left(\lambda t^{-2}\int_E\left(\sum_{c=1}^{nstd} \phi_c(\vm{x})\pi_c\left[\mat{N}_a(\vm{x})\right]\right)^\transpose\sum_{c=1}^{nstd} \phi_c(\vm{x})\pi_c\left[\mat{N}_b(\vm{x})\right]\,\diffx\right),\label{eq:discrete_system_f}\\
\intertext{and the nodal force vector is}
\mat{f}_w^E  &= \sum_{a=1}^{nstd}\left(\int_E\phi_{a}(\vm{x})q(\vm{x})\,\diffx\right).\label{eq:discrete_system_g}
\end{align}
\end{subequations}

It is recalled that in the implementation of these stiffness matrices and nodal force vector, $nstd$ and $nenh$ are the number of nodes that define the nodal contributions in the standard node set ($\mathcal{N}^{s}$) and the enhanced node set ($\mathcal{N}^{+}$), respectively, that result from the numerical integration process. The numerical integration of the nodal force vector is done using standard Gauss integration on triangles, but the numerical integration of the stiffness matrices requires a special treatment which is elaborated in the next subsection.

\subsection{Numerical integration}
\label{sec:numericalintegration}

The cell-based integration of the stiffness matrices that depend on basis functions derivatives (i.e., Eqs.~\eref{eq:discrete_system_c}-\eref{eq:discrete_system_e}) introduces integration errors when standard Gauss integration is used, which results in convergence issues and the patch test is not satisfied. To alleviate these integration errors in the \texttt{VANP} method, a smoothing procedure known as quadratically consistent 3-point integration scheme~\cite{duan:2012:SOI} is performed to correct the values of the basis functions derivatives at the integration points. This smoothing procedure was already used in the linear~\cite{ortiz:VAN:2015} and nonlinear~\cite{ortiz:IRN:2015} \texttt{VANP} formulations for nearly-incompressible solids. In this paper, we follow the same approach.

\begin{figure}
  \centering
  \mbox{
  \subfigure[]{\label{fig:modint_a}
  \epsfig{file = 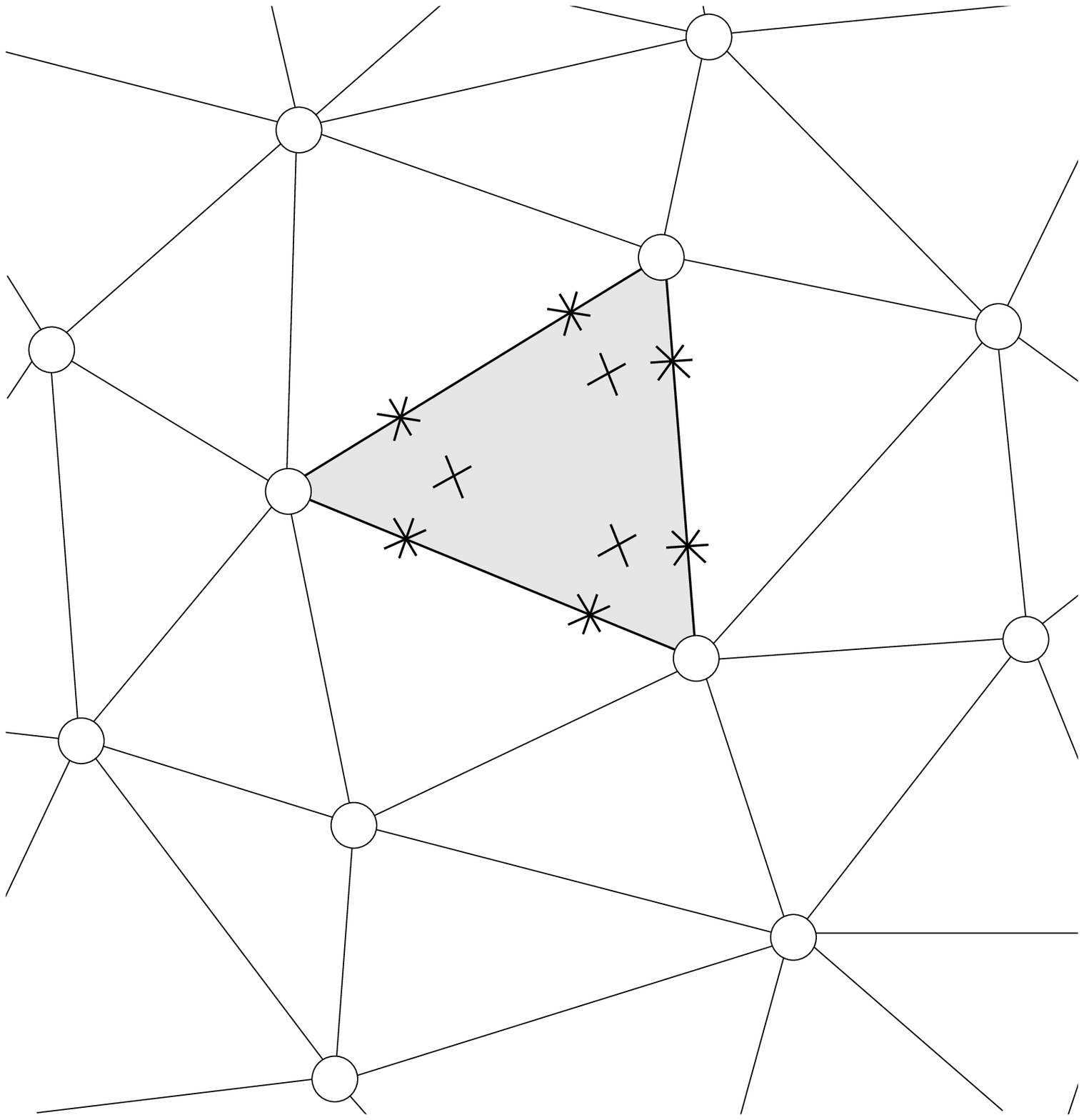, width = 0.45\textwidth}}
  \subfigure[]{\label{fig:modint_b} \epsfig{file = 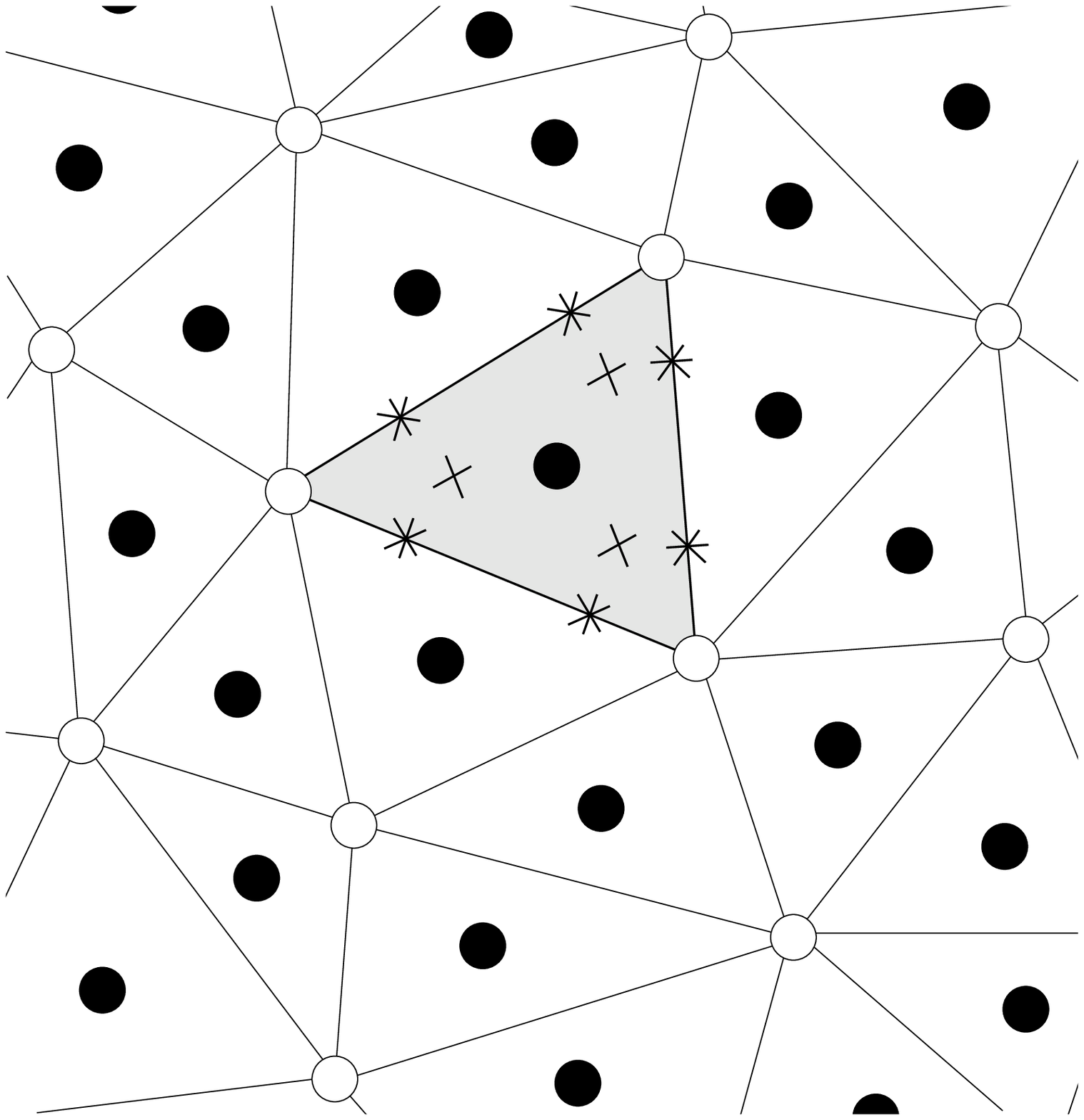, width = 0.45\textwidth}}
  }
  \caption{Schematic representation of integration cells and their integration points for correction of basis functions derivatives in the \texttt{VANP} approach. The shaded region is the integration cell $E$. The symbol $\times$ represents integration points located in the interior of the integration cell and the symbol $\ast$ represents integration points located on the cell boundary. (a) Integration cell when the standard node set ($\mathcal{N}^{s}$) is used and (b) integration cell when the enhanced node set ($\mathcal{N}^{+}$) is used.}
  \label{fig:modint}
\end{figure}

A representative integration cell $E$ along with its integration points is depicted in~\fref{fig:modint_a} when the standard node set ($\mathcal{N}^{s}$) is used and in~\fref{fig:modint_b} when the enhanced node set ($\mathcal{N}^{+}$) is used. Basically, the situation shown in~\fref{fig:modint_a} is used to evaluate the derivatives appearing in~\eref{eq:discrete_system_c} and~\eref{eq:discrete_system_d} through the nodal matrix $\mat{G}_a$ and the situation depicted in~\fref{fig:modint_b} to evaluate the derivatives appearing in~\eref{eq:discrete_system_e} through the nodal matrix $\mat{B}_a$.

A summary of the basis functions derivatives correction procedure follows. The Cartesian coordinate system is chosen, where for convenience $x\equiv x_1$ and $y\equiv x_2$. In addition, $n_j$ ($j=1,2$) is the $j$-th component of the unit outward normal to a cell edge in the Cartesian coordinate system. The integration accuracy of the smoothing procedure is of second-order, which is obtained by requiring the basis functions derivatives to satisfy the divergence constraint
\begin{equation}\label{eq:divconsistency}
\int_E\phi_{a,j}f(\vm{x})\,\diffx  =  \int_{\partial E}\phi_a f(\vm{x})n_j\,\diffs - \int_E\phi_a f_{,j}(\vm{x})\,\diffx,
\end{equation}
where $f(\vm{x})$ is the first-order polynomial base
\begin{equation}
f(\vm{x})=[1 \,\,\, x_1 \,\,\, x_2]^\transpose,
\end{equation}
whose derivative ($\delta_{ij}$ is the Kronecker delta symbol) is
\begin{equation}
f_{,j}(\vm{x})=[0 \,\,\, \delta_{1j} \,\,\, \delta_{2j}]^\transpose.
\end{equation}
The expanded version of~\eref{eq:divconsistency} is:
\begin{subequations}\label{eq:intconstraints}
\begin{align}
\int_E\phi_{a,1}\,\diffx =&  \int_{\partial E}\phi_a n_1 \,\diffs,\\
\int_E\phi_{a,1}x_1\,\diffx =&  \int_{\partial E}\phi_a x_1 n_1\,\diffs-\int_E\phi_a\,\diffx,\\
\int_E\phi_{a,1}x_2\,\diffx =&  \int_{\partial E}\phi_a x_2 n_1 \,\diffs,\\
\intertext{for $\phi_{a,1}$, and}
\int_E\phi_{a,2}\,\diffx =&  \int_{\partial E}\phi_a n_2 \,\diffs,\\
\int_E\phi_{a,2}x_1\,\diffx =&  \int_{\partial E}\phi_a x_1 n_2\,\diffs,\\
\int_E\phi_{a,2}x_2\,\diffx =&  \int_{\partial E}\phi_a x_2 n_2 \,\diffs-\int_E\phi_a\,\diffx
\end{align}
\end{subequations}
for $\phi_{a,2}$.

The integration constraints~\eref{eq:intconstraints} are solved using Gauss integration on the integration cell $E$ shown in~\fref{fig:modint}. Consider the following notations:
\begin{itemize}
\item $\Prefixd^{h}{\vm{p}}=(\Prefixd^{h}{p_1},\Prefixd^{h}{p_2})$ as the Cartesian coordinates of the $h$-th interior integration point with an associated Gauss weight $\Prefixd^{h}{w}$.
\item $\Prefixf_{k}^{g}{\vm{e}}=(\Prefixf_{k}^{g}{e_1},\Prefixf_{k}^{g}{e_2})$ as the Cartesian coordinates of the $g$-th integration point that is located on the $k$-th edge of the cell with an associated Gauss weight $\Prefixf_{k}^{g}{v}$.
\item $\Prefixd_{k}{\vm{n}}=(\Prefixd_{k}{n_1},\Prefixd_{k}{n_2})$ as the unit outward normal to the $k$-th edge of the cell.
\end{itemize}

On using the preceding notations, the discrete version of the integration constraints~\eref{eq:intconstraints} is:
\begin{subequations}\label{eq:discreteintconstraints}
\begin{equation}
\mat{W}\mat{d}_j =  \mat{f}_j, \,\,\, j={1,2}
\end{equation}
where
\begin{equation}
\mat{W} = \left[
\begin{array}{c c c}
\Prefixd^{1}{w} & \Prefixd^{2}{w} & \Prefixd^{3}{w}  \\
\Prefixd^{1}{w}\,\Prefixd^{1}{p_1} & \Prefixd^{2}{w}\,\Prefixd^{2}{p_1} & \Prefixd^{3}{w}\,\Prefixd^{3}{p_1}  \\
\Prefixd^{1}{w}\,\Prefixd^{1}{p_2} & \Prefixd^{2}{w}\,\Prefixd^{2}{p_2} & \Prefixd^{3}{w}\,\Prefixd^{3}{p_2}
\end{array}
\right],
\end{equation}
\begin{equation}
\mat{f}_1 = \left[
\begin{array}{l}
\sum\limits_{k=1}^3\sum\limits_{g=1}^2 \phi_a(\Prefixf_{k}^{g}{\vm{e}})\,\Prefixd_{k}{n_1}\,\Prefixf_{k}^{g}{v}\\
\sum\limits_{k=1}^3\sum\limits_{g=1}^2 \phi_a(\Prefixf_{k}^{g}{\vm{e}})\,\Prefixf_{k}^{g}{e_1}\,\Prefixd_{k}{n_1}\,\Prefixf_{k}^{g}{v}
-\sum\limits_{h=1}^3\phi_a(\Prefixd^{h}{\vm{p}})\,\Prefixd^{h}{w}\\
\sum\limits_{k=1}^3\sum\limits_{g=1}^2 \phi_a(\Prefixf_{k}^{g}{\vm{e}})\,\Prefixf_{k}^{g}{e_2}\,\Prefixd_{k}{n_1}\,\Prefixf_{k}^{g}{v}
\end{array}
\right],
\end{equation}
\begin{equation}
\mat{f}_2 = \left[
\begin{array}{l}
\sum\limits_{k=1}^3\sum\limits_{g=1}^2 \phi_a(\Prefixf_{k}^{g}{\vm{e}})\,\Prefixd_{k}{n_2}\,\Prefixf_{k}^{g}{v}\\
\sum\limits_{k=1}^3\sum\limits_{g=1}^2 \phi_a(\Prefixf_{k}^{g}{\vm{e}})\,\Prefixf_{k}^{g}{e_1}\,\Prefixd_{k}{n_2}\,\Prefixf_{k}^{g}{v}\\
\sum\limits_{k=1}^3\sum\limits_{g=1}^2 \phi_a(\Prefixf_{k}^{g}{\vm{e}})\,\Prefixf_{k}^{g}{e_2}\,\Prefixd_{k}{n_2}\,\Prefixf_{k}^{g}{v}
-\sum\limits_{h=1}^3\phi_a(\Prefixd^{h}{\vm{p}})\,\Prefixd^{h}{w}
\end{array}
\right],
\end{equation}
and the solution vector of the $j$-th basis function derivative evaluated at the three interior integration points is
\begin{equation} \label{eq:duanderivatives}
\mat{d}_j = \left[
\begin{array}{c c c }
\phi_{a,j}(\Prefixd^{1}{\vm{p}}) &
\phi_{a,j}(\Prefixd^{2}{\vm{p}}) &
\phi_{a,j}(\Prefixd^{3}{\vm{p}})
\end{array}
\right]^\transpose.
\end{equation}
\end{subequations}

In the foregoing equations, index $a$ runs through the nodes that define the nodal contribution either in $\mathcal{N}^{s}$ or $\mathcal{N}^{+}$.

Finally, the numerical integration of the stiffness matrices is performed using 3-point Gauss rule on the triangular cells, but the basis functions derivatives at these three integration points are replaced by the corrected derivatives given in~\eref{eq:duanderivatives}.

\section{Numerical experiments}
\label{sec:numexamples}

In this section, several numerical experiments are performed to assess the accuracy of the proposed \texttt{VANP} formulation for the Reissner-Mindlin plate model. Unless stated otherwise, the default numerical integration procedure for the \texttt{VANP} formulation is the quadratically consistent 3-point integration scheme (QC3). For assessing Reissner-Mindlin plate problems with known global solution, we use the relative $L^2$-norm of error and the relative $H^1$-seminorm of the error, which are defined, respectively, as follows:
\begin{equation*}
\frac{\|\mat{u}-\mat{u}^h\|_{L^2(\Omega)}}{\|\mat{u}\|_{L^2(\Omega)}}
=\sqrt{\frac{\sum_E\int_E\left(\mat{u}-\mat{u}^h\right)^\transpose
       \left(\mat{u}-\mat{u}^h\right)\,\diffx}
       {\sum_E\int_E\mat{u}^\transpose\mat{u}\,\diffx}},
\end{equation*}
\begin{equation*}
\frac{\|\mat{u}-\mat{u}^h\|_{H^1(\Omega)}}{\|\mat{u}\|_{H^1(\Omega)}}
=\sqrt{\frac{\sum_E\int_E\left(\mat{u}'-{\mat{u}'}^h\right)^\transpose
       \left(\mat{u}'-{\mat{u}'}^h\right)\,\diffx}
       {\sum_E\int_E{\mat{u}'}^\transpose\mat{u}'\,\diffx}},
\end{equation*}
where $\mat{u}=[w \,\, r_x \,\, r_y]^\transpose$ and $\mat{u}'=[w_{,x} \,\, w_{,y} \,\, r_{x,x} \,\, r_{x,y} \,\, r_{y,x} \,\, r_{y,y}]^\transpose$ are the exact solutions, and $\mat{u}^h$ and ${\mat{u}'}^h$ are their corresponding approximations.

\alejandro{In addition, on using the exact nodal scaled transverse shear stress solution $\vm{s}_a$ and its approximation $\vm{s}_a^h$, the following relative $L^2$-norm of the nodal error is defined to assess the convergence of the \texttt{VANP} method in the scaled transverse shear stress variable:
\begin{equation*}
\frac{\|\vm{s}-\vm{s}^h\|_{L^2(\Omega)}}{\|\vm{s}\|_{L^2(\Omega)}}
=\sqrt{\frac{\sum_a\left(\vm{s}_a-\vm{s}_a^h\right)^\transpose
       \left(\vm{s}_a-\vm{s}_a^h\right)}
       {\sum_a\vm{s}_a^\transpose\vm{s}_a}},
\end{equation*}
since in the \texttt{VANP} approach the scaled transverse shear stress is a nodal quantity that can be computed a \textit{posteriori} from the primitive variables using~\eref{eq:nodal_shearstress_lumped}.}

\subsection{Zero shear deformation patch test}
\label{sec:patch_test}
\alejandro{We start by performing a patch test to evaluate whether our \texttt{VANP} formulation, which uses linear approximations, can reproduce a linear solution within machine precision. We also want to check that our method is devoid of the shear-locking phenomenon when the transverse shear deformation approaches zero. For the Reissner-Mindlin problem, the condition of zero shear deformation requires that the transverse displacement $w$ is one order higher than the order of the rotations $\vm{r}$. To obtain an exact linear solution for $w$ and be able to test our method, we use the zero shear deformation patch test that is provided in Ref.~\cite{chen_wang_zhao_2009}. The lowest order solution given therein is quadratic in $w$ and linear in $\vm{r}$. To obtain a linear solution in $w$, the exact solution provided in Ref.~\cite{chen_wang_zhao_2009} is managed by appropriately choosing the arbitrary constants so that the following particular exact solution is obtained:
\begin{equation*}
w = 1 + x + y,\quad r_x=1, \quad r_y=1.
\end{equation*}

The linear patch test is built by imposing the above exact solution along the entire boundary of a unit square domain. Three integration meshes are used as shown in~\fref{fig:patch_test_meshes}.
\begin{figure}[!htbp]
\centering
\mbox{
\subfigure[]{\label{fig:patch_test_meshes_a} \epsfig{file = 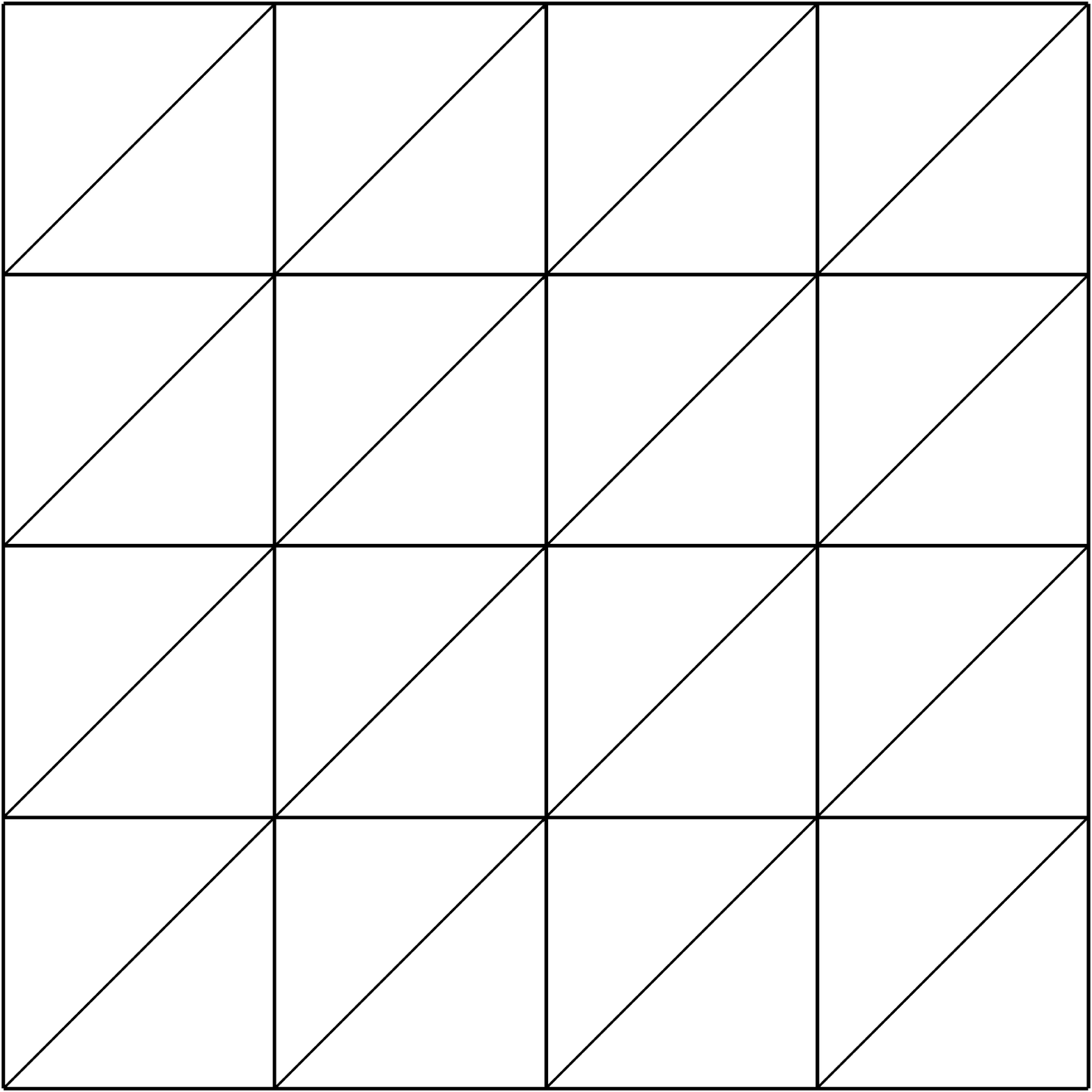, width = 0.3\textwidth}}
\subfigure[]{\label{fig:patch_test_meshes_b} \epsfig{file = 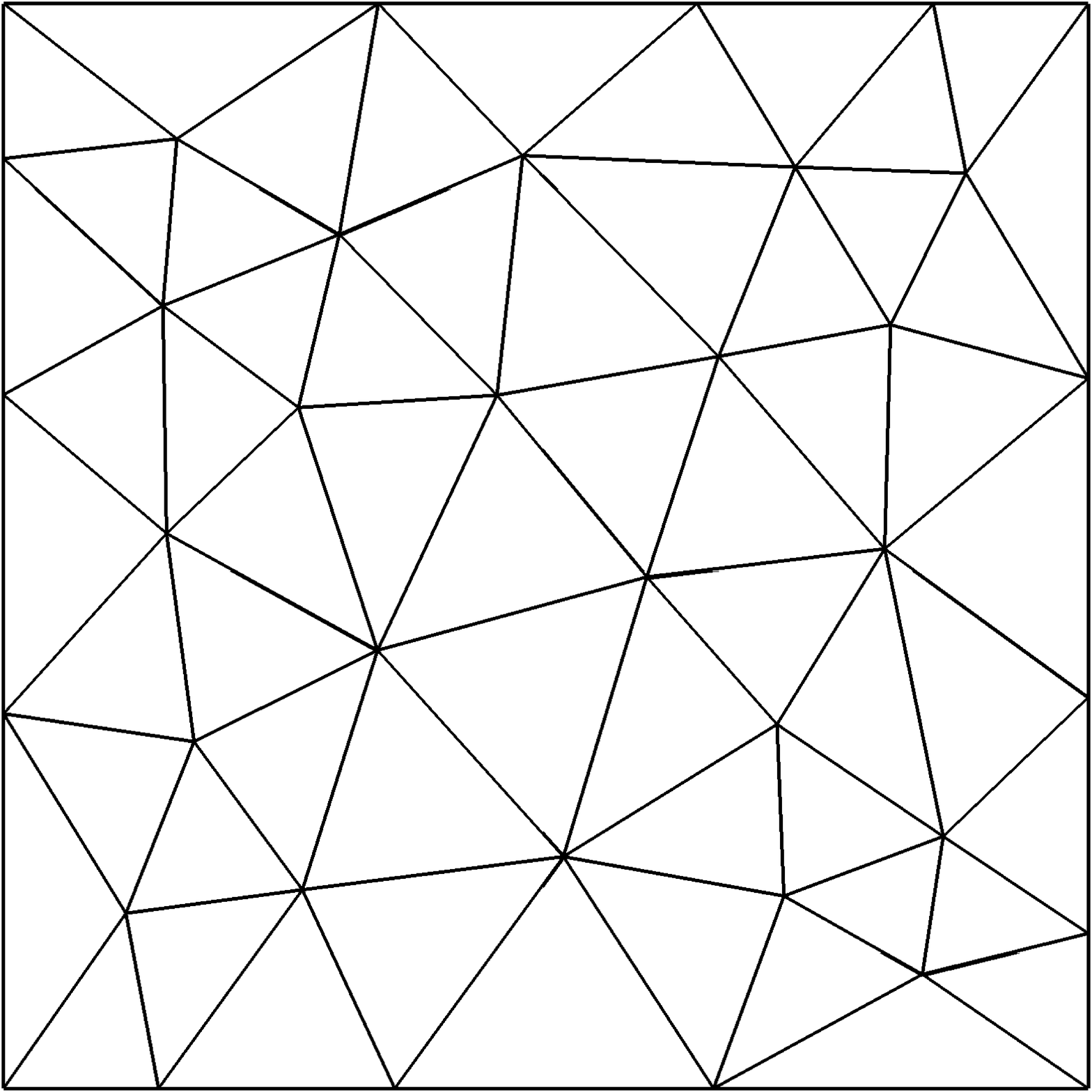, width = 0.3\textwidth}}
\subfigure[]{\label{fig:patch_test_meshes_c} \epsfig{file = 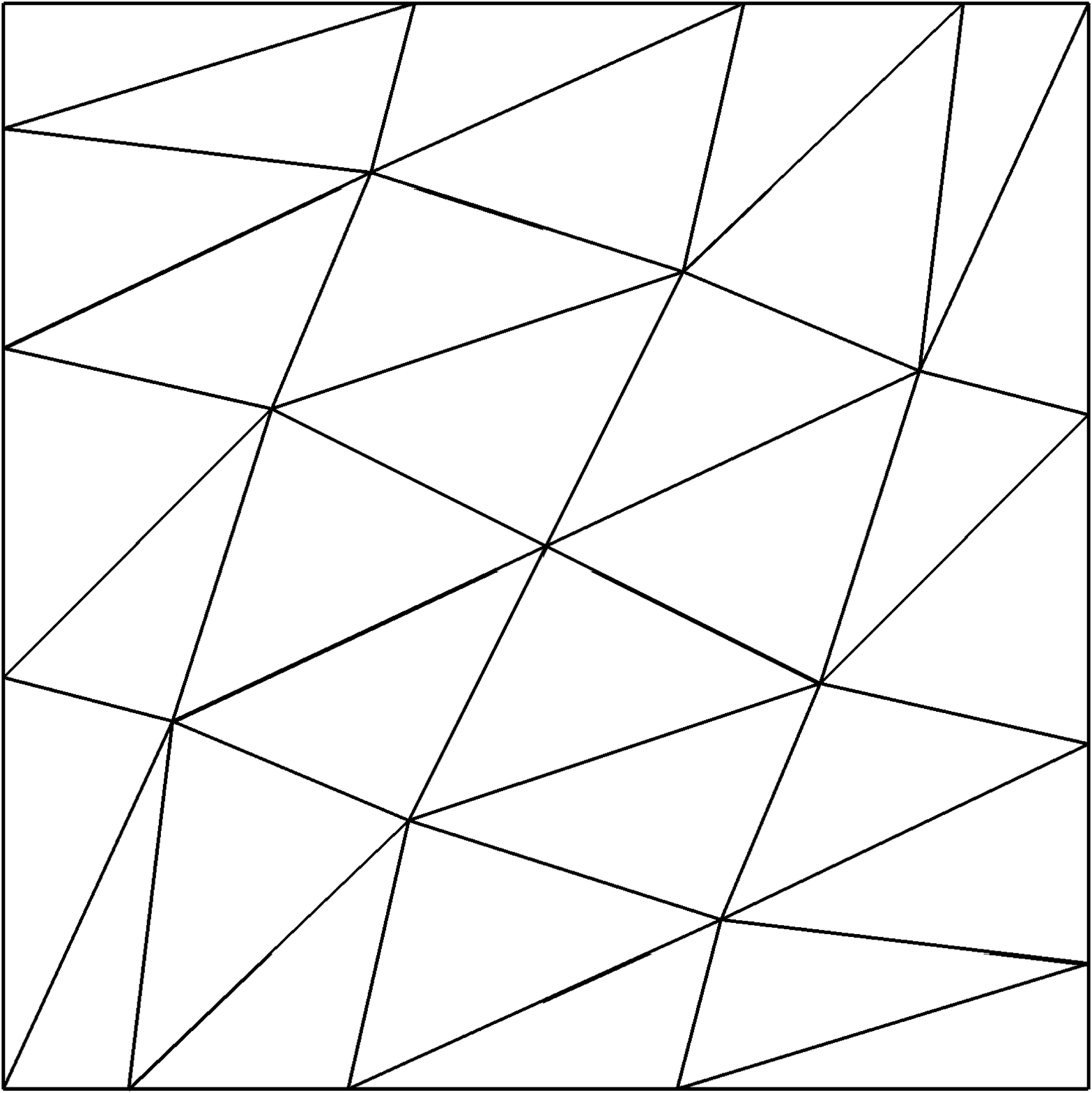, width = 0.3\textwidth}}
}
\caption{Integration meshes for the zero shear deformation patch test.}
\label{fig:patch_test_meshes}
\end{figure}
The elastic parameters for the material of the plate are set to $E_\textrm{Y}=10.92\times 10^6$ psi and $\nu=0.3$. The side of the plate $L$, which in this case is the unit, is taken as the characteristic length for defining the normalized thickness of the plate as $t/L$. The relative $L^{2}$-norm and $H^1$-seminorm of the error are shown in~\tref{tab:patchtest_L2_errors} and~\tref{tab:patchtest_H1_errors}, respectively, for the three integration meshes and various normalized thicknesses. These results reveal that for the plates with normalized thicknesses $t/L=0.1$, $t/L=0.01$ and $t/L=0.001$, the errors are extremely small and approaching machine precision; thus, it can be said that they pass the zero shear deformation patch test in the numerical sense. Even though the errors for the plate with normalized thickness $t/L=0.0001$ are not within machine precision, they are sufficiently small to nearly pass the zero shear deformation patch test in the numerical sense. Also, the absence of the shear-locking phenomenon in the \texttt{VANP} formulation is made evident by these small errors.}

\begin{table}[!htbp]
\centering
\caption{Relative $L^{2}$-norm of the error for the zero shear deformation patch test.}
\begin{tabular}{|c|c|c|c|c|}
\hline
Mesh & $t/L=0.1$            & $t/L=0.01$          & $t/L=0.001$          & $t/L=0.0001$        \\ \hline
(a)  & $2.3\times 10^{-14}$ & $2.6\times 10^{-12}$ & $1.5\times 10^{-10}$ & $2.0\times 10^{-8}$ \\
(b)  & $1.5\times 10^{-14}$ & $6.3\times 10^{-13}$ & $6.0\times 10^{-11}$ & $6.2\times 10^{-9}$ \\
(c)  & $1.7\times 10^{-14}$ & $7.5\times 10^{-13}$ & $4.2\times 10^{-11}$ & $5.7\times 10^{-9}$ \\ \hline
\end{tabular}
\label{tab:patchtest_L2_errors}
\end{table}

\begin{table}[!htbp]
\centering
\caption{Relative $H^1$-seminorm of the error for the zero shear deformation patch test.}
\begin{tabular}{|c|c|c|c|c|}
\hline
Mesh & $t/L=0.1$            & $t/L=0.01$          & $t/L=0.001$          & $t/L=0.0001$       \\ \hline
(a)  & $4.7\times 10^{-14}$ & $9.7\times 10^{-13}$ & $5.3\times 10^{-11}$ & $7.5\times 10^{-9}$ \\
(b)  & $3.4\times 10^{-13}$ & $7.0\times 10^{-13}$ & $5.0\times 10^{-11}$ & $6.0\times 10^{-9}$ \\
(c)  & $2.1\times 10^{-13}$ & $3.9\times 10^{-13}$ & $2.5\times 10^{-11}$ & $5.6\times 10^{-9}$ \\ \hline
\end{tabular}
\label{tab:patchtest_H1_errors}
\end{table}

\subsection{Circular plate subjected to a uniform load}
\label{sec:numexamples_falk}

\begin{figure}[!htbp]
  \centering
  \epsfig{file = 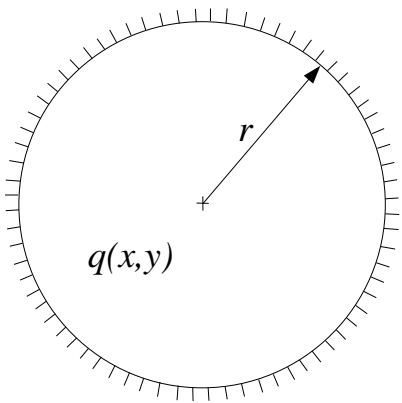, width = 0.35\textwidth}
  \caption{Circular plate subjected to a uniform load.}
  \label{fig:falk_prob}
\end{figure}

\fref{fig:falk_prob} depicts a circular plate of radius $r$ that is subjected to a uniform load $q$ and is clamped along its entire boundary. The normalized thickness of the plate is $t/L$, where $t$ is the thickness of the plate and $L$ is a characteristic length of the physical domain, which in this case is taken as the radius of the plate. The radius of the plate is set to $r=1$ in so that $L=1$ in, and the uniform load is set to $q=1$ psi. The following elastic parameters are considered for the material of the plate: $E_\textrm{Y}=10.92\times 10^6$ psi and $\nu=0.3$. The exact solution for this problem is given by~\cite{falk:LFFEMRM:2000}
\begin{align*}
r_x &= \frac{x(x^2+y^2-1)}{16D},\quad r_y=\frac{y(x^2+y^2-1)}{16D},\\
w&=\frac{(x^2+y^2)^2}{64D}-(x^2+y^2)\left(\frac{\lambda^{-1}t^2}{4}+\frac{1}{32D}\right)+\frac{1}{4}\lambda^{-1}t^2+\frac{1}{64D},
\end{align*}
where $D=E_\textrm{Y}/(12(1-\nu^2))$.

\begin{figure}[!tbhp]
\centering
\mbox{
\subfigure[]{\label{fig:meshes_falk_prob_a} \epsfig{file = 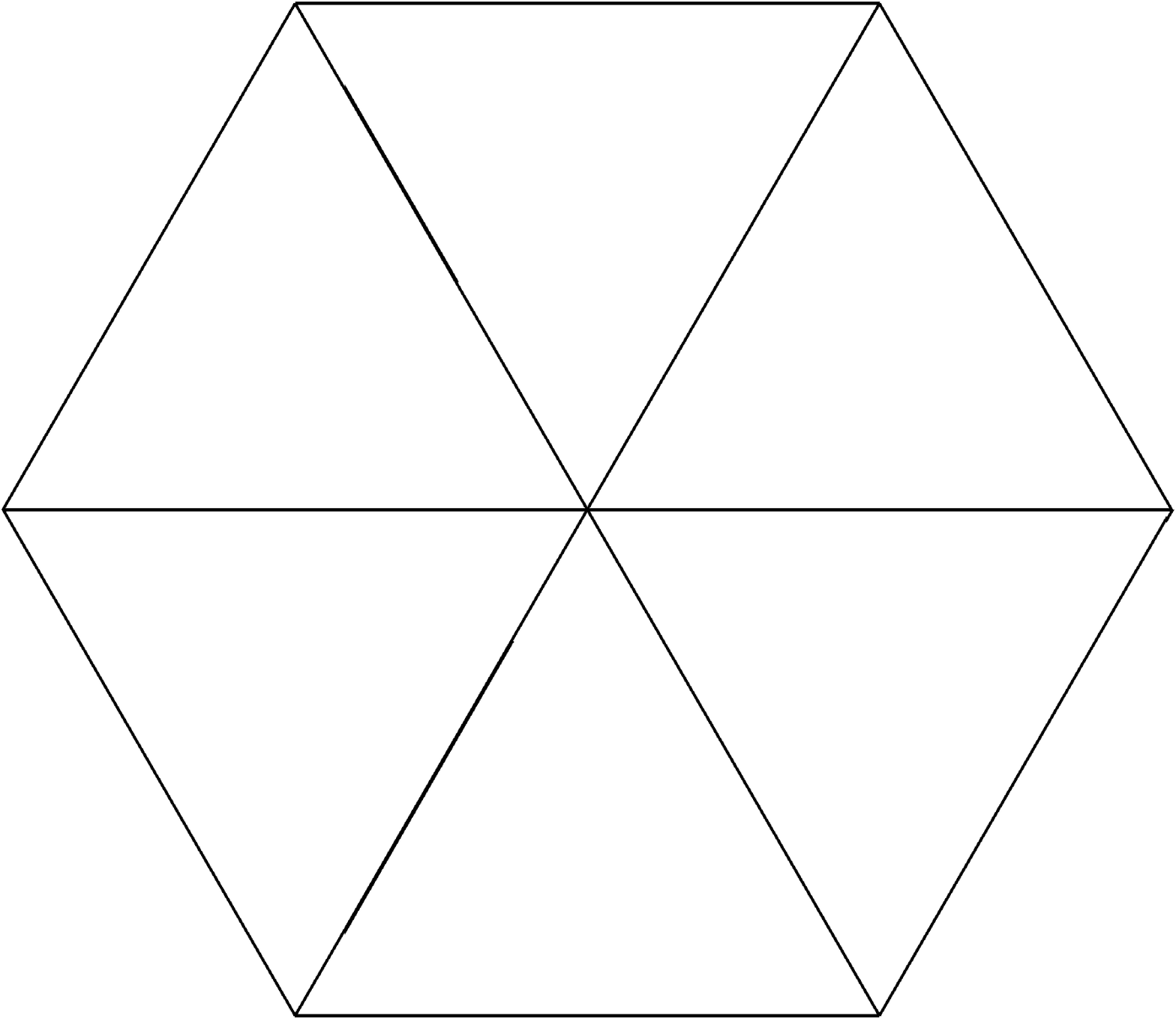, width = 0.2\textwidth}}
\subfigure[]{\label{fig:meshes_falk_prob_b} \epsfig{file = 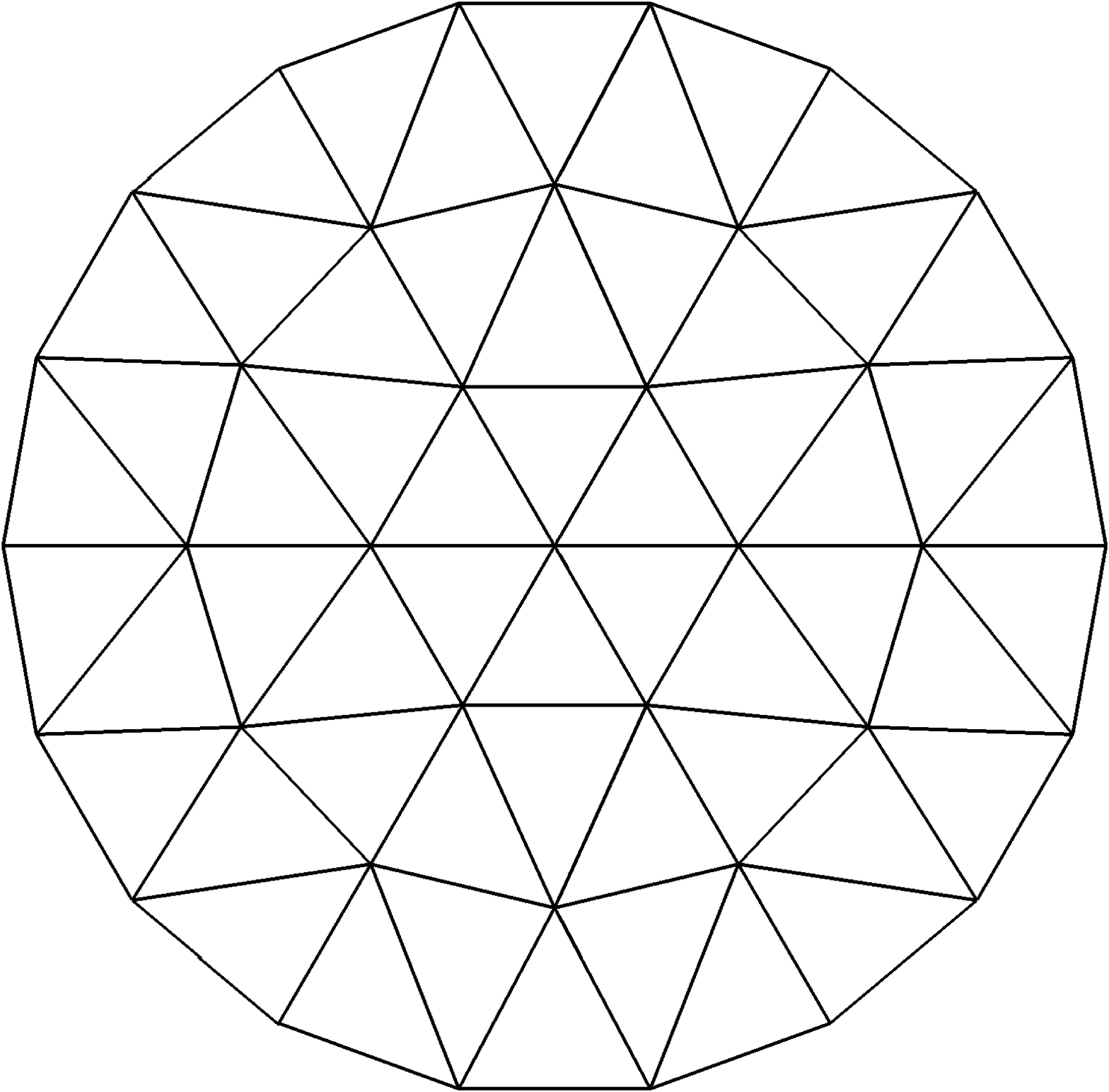, width = 0.2\textwidth}}
\subfigure[]{\label{fig:meshes_falk_prob_c} \epsfig{file = 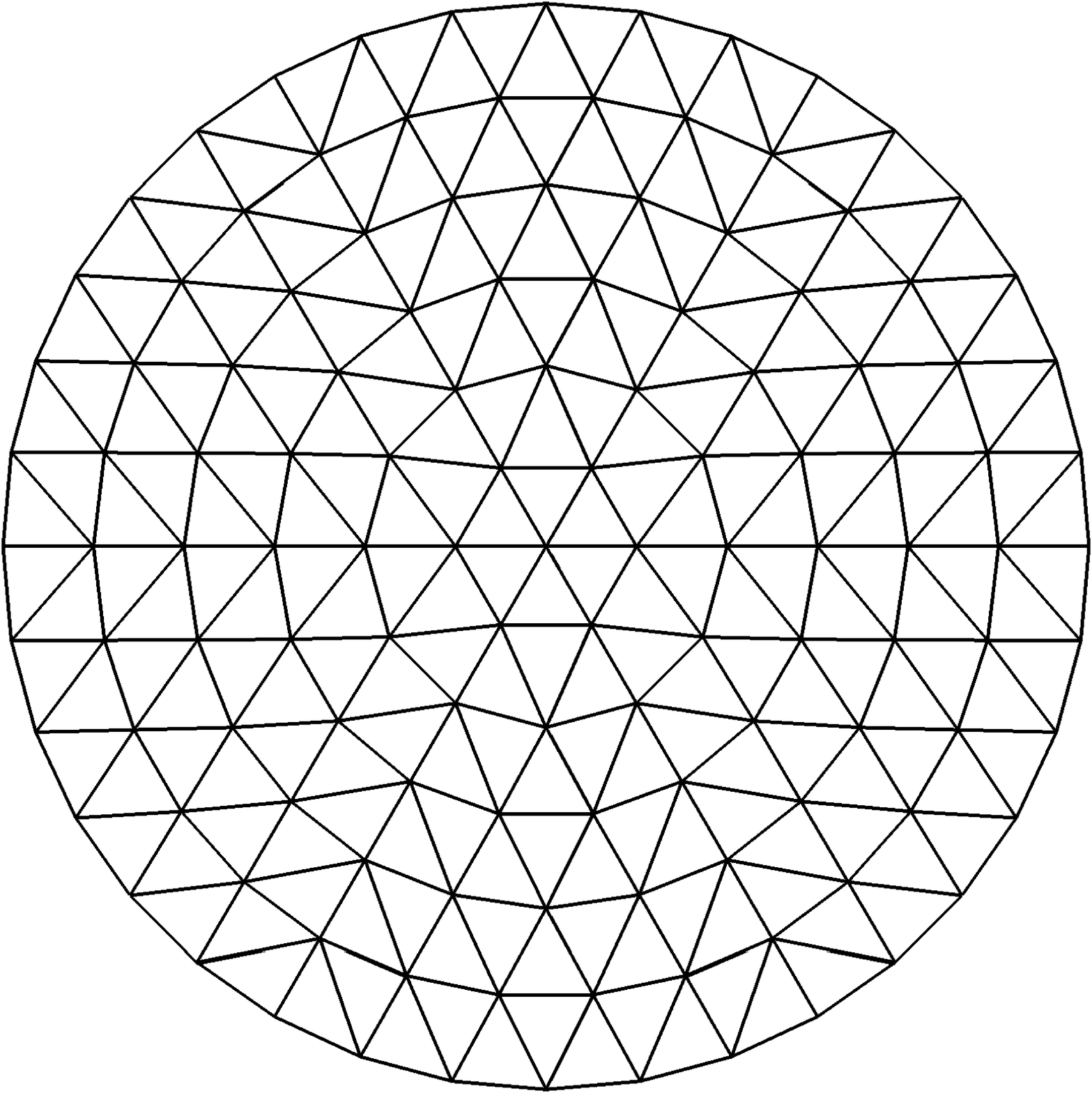, width = 0.2\textwidth}}
\subfigure[]{\label{fig:meshes_falk_prob_d} \epsfig{file = 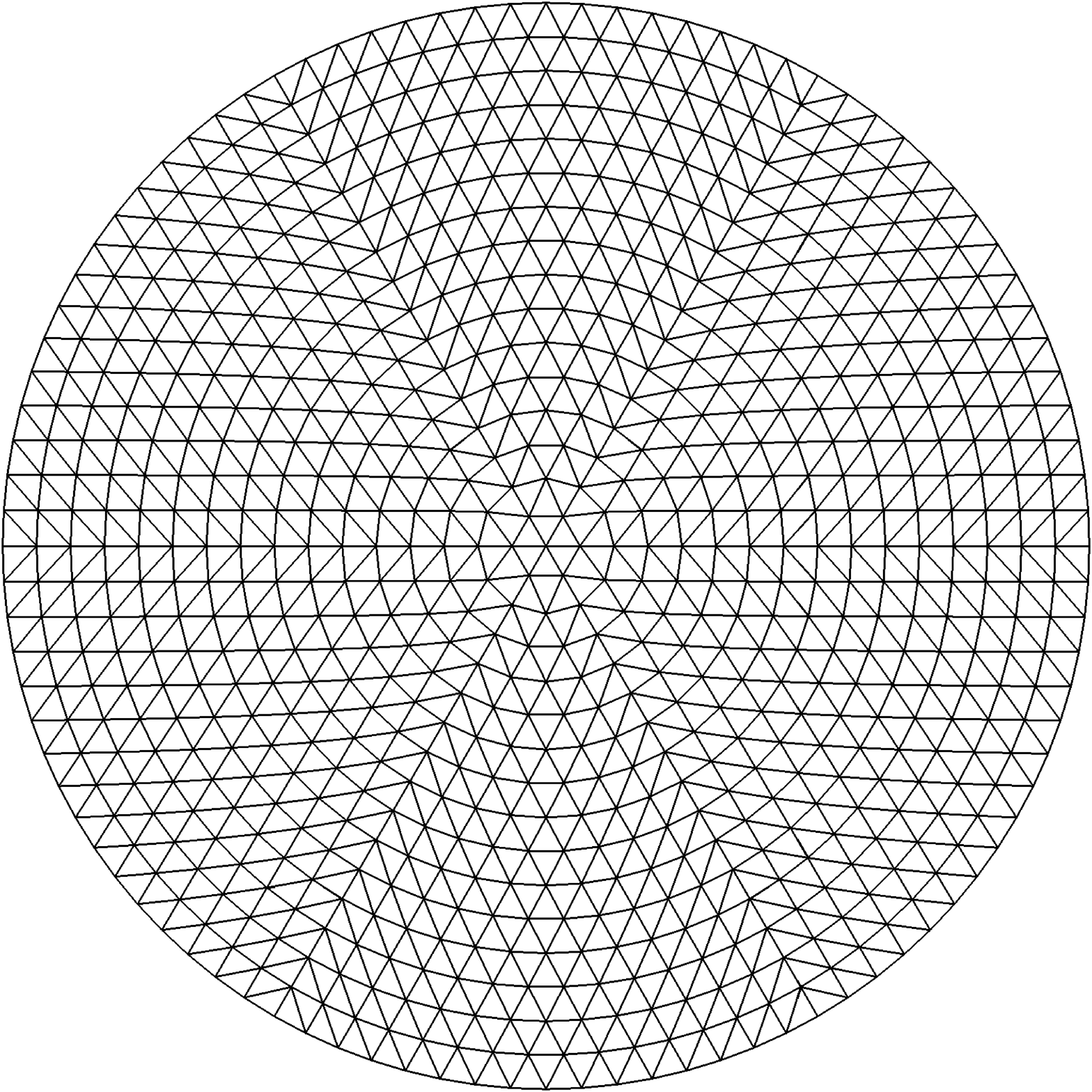, width = 0.2\textwidth}}
}
\caption{Integration meshes for the circular plate subjected to a uniform load problem.}
\label{fig:meshes_falk_prob}
\end{figure}

The integration meshes that are considered for this problem are shown in~\fref{fig:meshes_falk_prob}.

\begin{figure}[!tbhp]
\centering
\mbox{
\subfigure[]{\label{fig:norms_falk_prob_a} \epsfig{file = 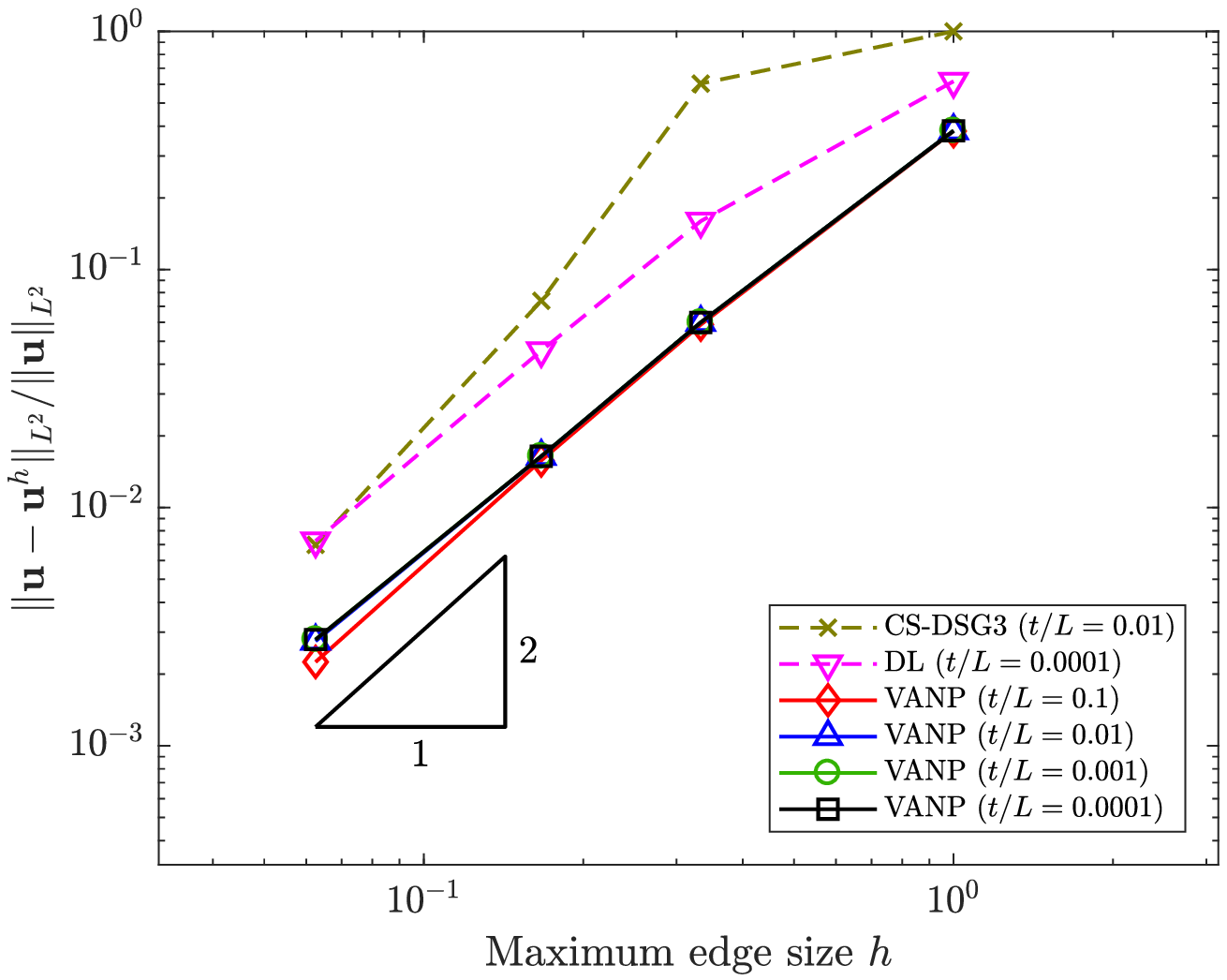, width = 0.5\textwidth}}
\subfigure[]{\label{fig:norms_falk_prob_b} \epsfig{file = 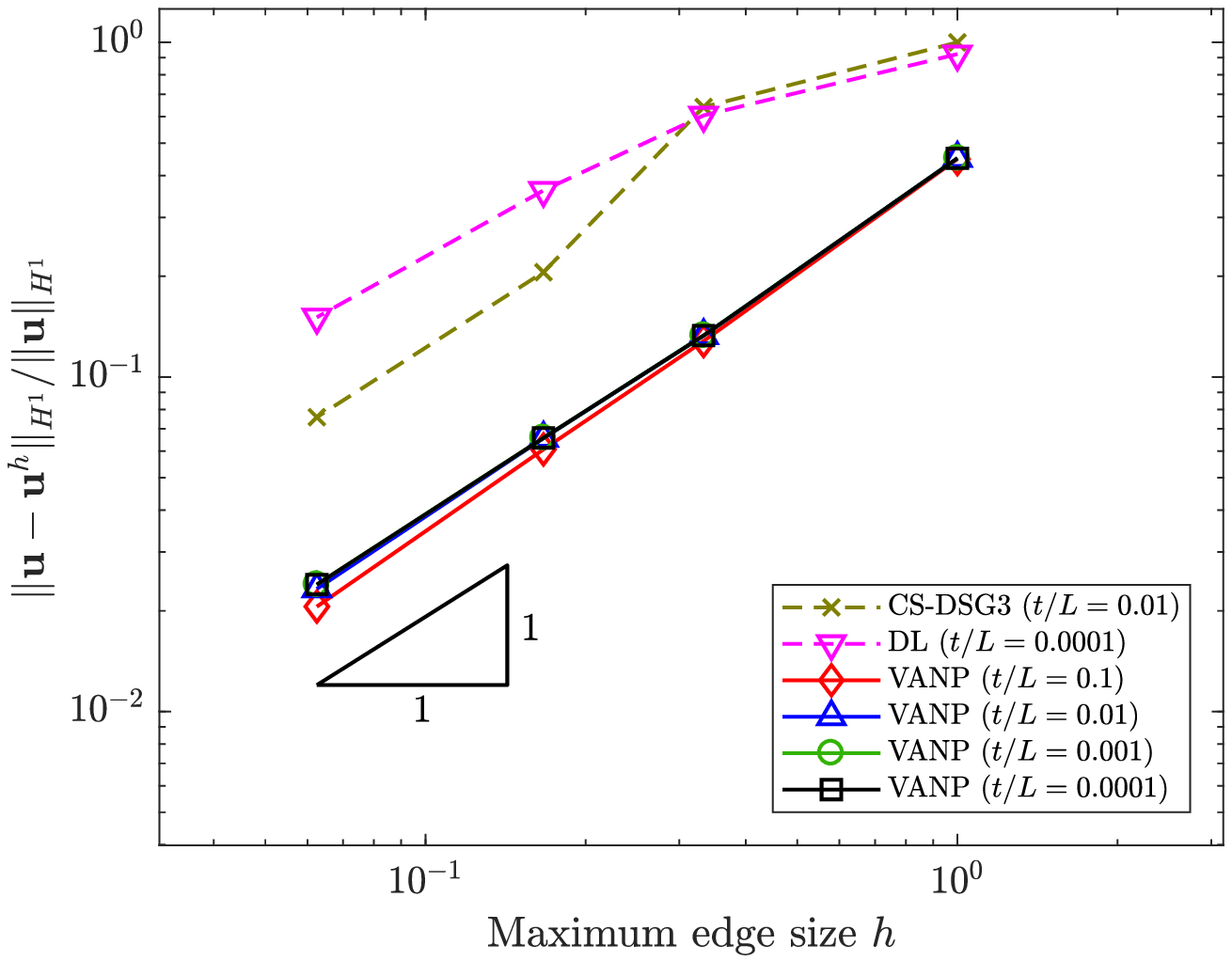, width = 0.5\textwidth}}
}
\caption{Rates of convergence for the circular plate subjected to a uniform load. (a) $L^2$-norm of the error and (b) $H^1$-seminorm of the error for several values of $t/L$. The \texttt{VANP}, \alejandro{DL and CS-DSG3} approaches deliver the optimal rates of convergence, but the accuracy of the \texttt{VANP} approach is superior to the accuracy of \alejandro{the DL and CS-DSG3 elements}.}
\label{fig:norms_falk_prob}
\end{figure}

We start by studying the convergence of the proposed \texttt{VANP} formulation as the integration mesh is refined. For comparison purposes, we also include the convergence results for the mixed triangular finite element of Dur\'an and Liberman~\cite{duran_mixed_1992}, which we denote by DL, \alejandro{and the three-node triangular element with cell-based smoothing for bending strain and discrete shear gap method for shear-locking (CS-DSG3) of Nguyen-Thoi et al.~\cite{nguyen_phung_thai_2012}}. The following normalized thicknesses are considered for the \texttt{VANP} approach: $t/L=\{0.1,\, 0.01,\, 0.001,\, 0.0001\}$. For the DL element, we only show the convergence curve for $t/L=0.0001$ since the curves for the other normalized thicknesses do not change significantly. \alejandro{For the CS-DSG3 element only the curve for $t/L=0.01$ is shown because this element did not perform well for thinner plates}. The convergence rates are shown in~\fref{fig:norms_falk_prob}, where it is observed that the optimal rates of convergence, 2 and 1, \alejandro{are delivered by the \texttt{VANP}, DL and CS-DSG3 approaches} in both the $L^2$-norm and the $H^1$-seminorn of the error, respectively. However, the accuracy of the \texttt{VANP} formulation is superior to the accuracy of the \alejandro{DL and CS-DSG3 elements}.

\begin{figure}[!tbhp]
\centering
\mbox{
\subfigure[]{\label{fig:norms_gamma_falk_prob_a} \epsfig{file = 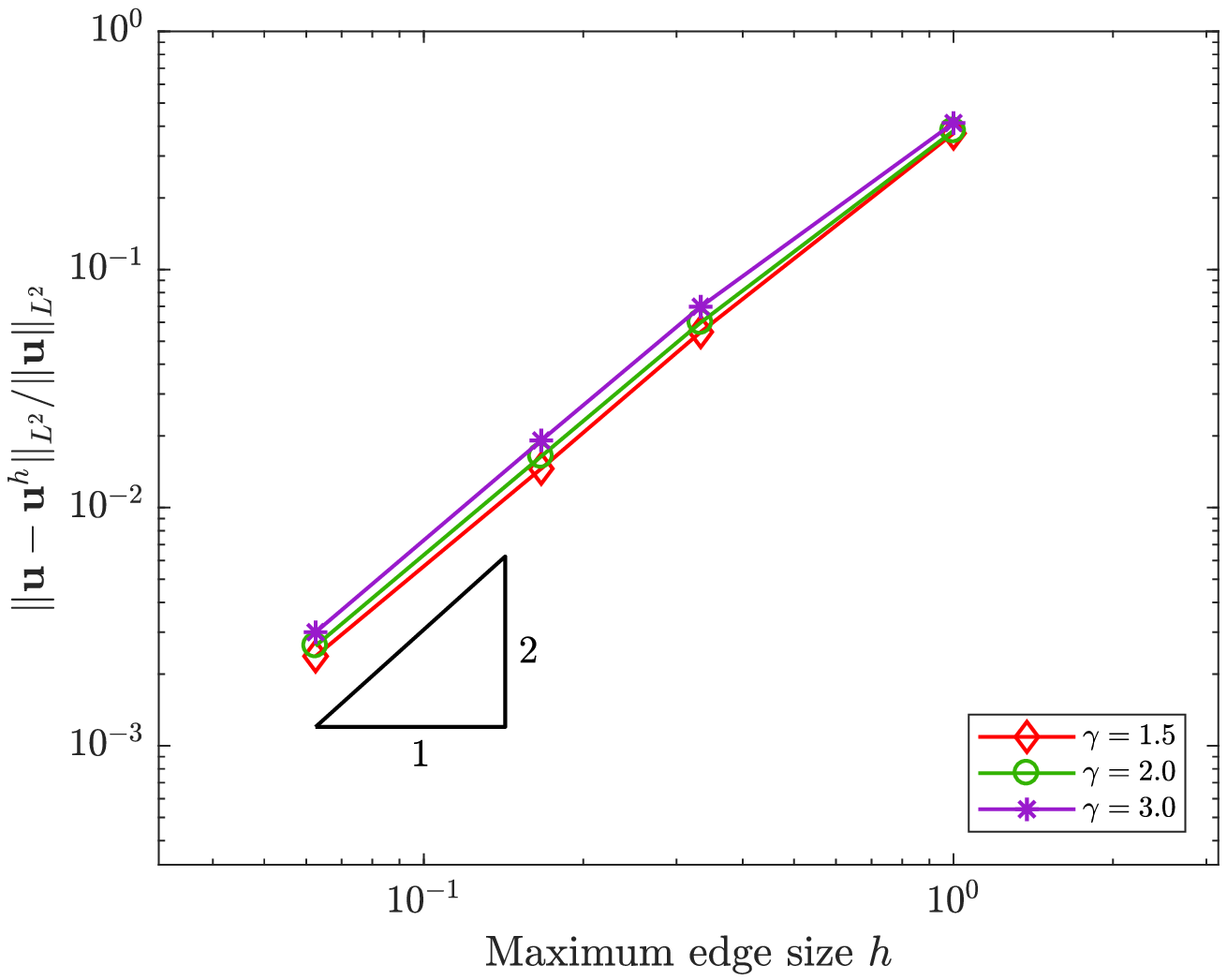, width = 0.5\textwidth}}
\subfigure[]{\label{fig:norms_gamma_falk_prob_b} \epsfig{file = 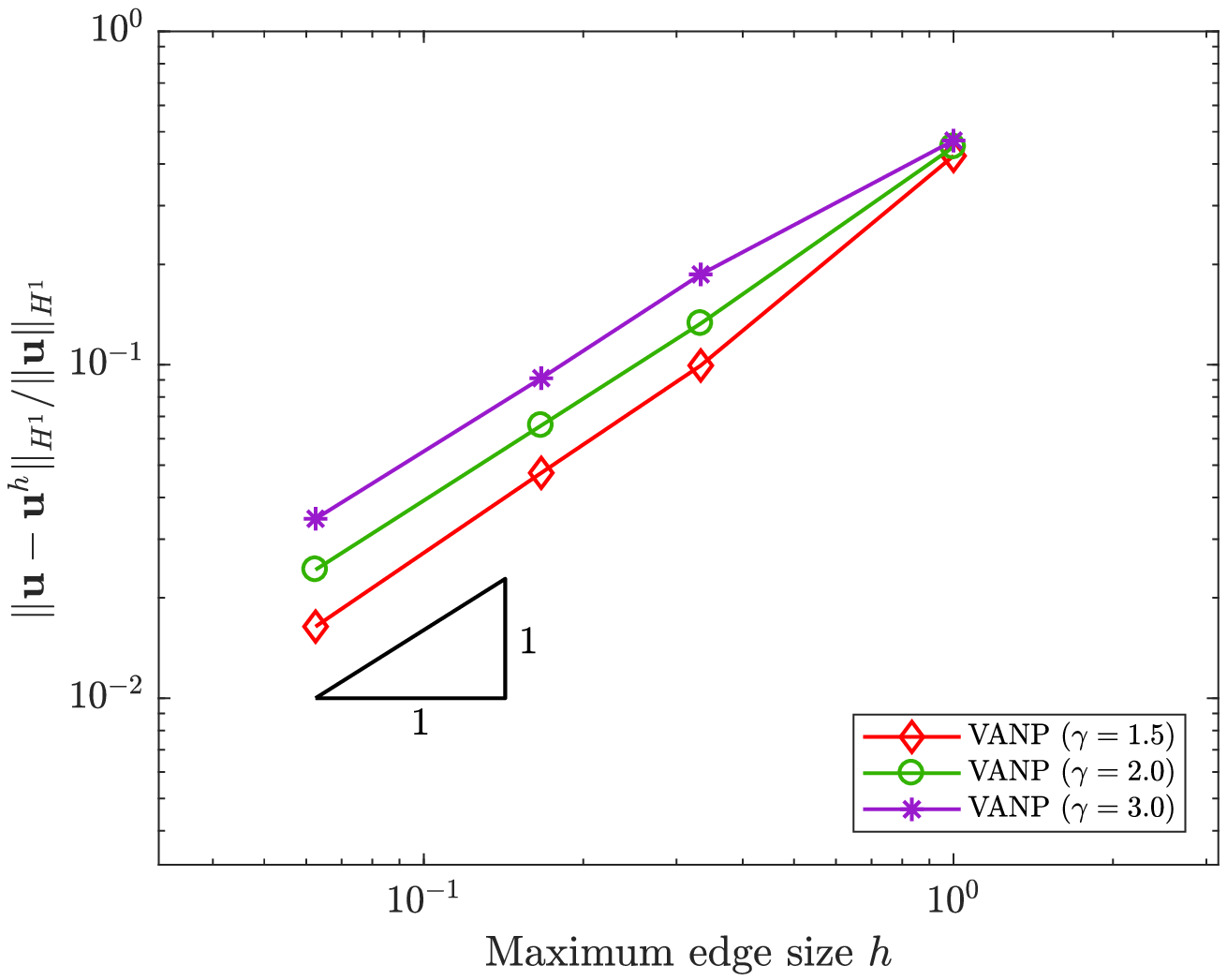, width = 0.5\textwidth}}
}
\caption{Influence of the maxent basis function support parameter ($\gamma$) on the \texttt{VANP} convergence rates. Three values for $\gamma$ are considered. Optimal convergence rates in the (a) $L^2$ norm and (b) the $H^1$ seminorm of the error are obtained for all these cases.}
\label{fig:norms_gamma_falk_prob}
\end{figure}

We also study the sensitivity of the convergence rates to the support parameter ($\gamma$) of the maxent basis functions. Three values are considered: $\gamma=\{1.5,\, 2.0,\, 3.0\}$, where the largest one results in the smaller support. For this test, the normalized thickness $t/L=0.0001$ is considered. The convergence rates are presented in~\fref{fig:norms_gamma_falk_prob}, where it is observed that the optimal rates of 2 and 1 are delivered by the \texttt{VANP} formulation in both the $L^2$-norm and the $H^1$-seminorn of the error, respectively, independently of the basis function support parameter. It is also observed that the \texttt{VANP} accuracy decreases as the support gets smaller, which is a reasonable behavior since as the support gets smaller the maxent basis function approaches the ``hat'' finite element basis function~\cite{arroyo:2006:LME}.

\subsection{Square plate subjected to a nonuniform load}
\label{sec:numexamples_chinosi}

\begin{figure}[!htbp]
  \centering
  \epsfig{file = 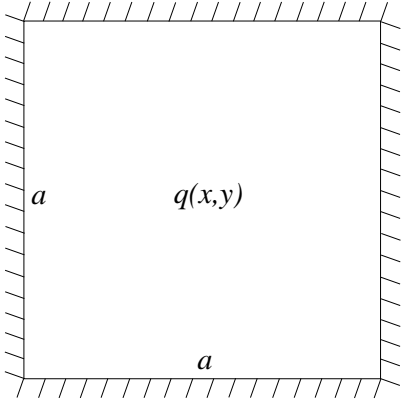, width = 0.35\textwidth}
  \caption{Square plate subjected to a nonuniform load.}
  \label{fig:chinosi_prob}
\end{figure}

In this example, we study the convergence properties of the \texttt{VANP} formulation in a more complicated setting, which includes nonuniform integration meshes and a nonuniform load. As shown in~\fref{fig:chinosi_prob}, the problem domain is a square plate that is clamped along its entire boundary. The side of the plate is taken as the characteristic length for defining the normalized thickness of the plate as $t/L$. In this problem, we set the side of the plate to $a=1$ in so that the characteristic length becomes $L=1$ in. The following elastic parameters are considered for the material of the plate: $E_\textrm{Y}=10.92\times 10^6$ psi and $\nu=0.3$. The nonuniform load is given by
\begin{eqnarray*}
q &=& \frac{E_\textrm{Y}}{12(1-\nu^2)}\left[12y(y-1)(5x^2-5x+1)(2y^2(y-1)^2+x(x-1)(5y^2-5y+1))\right.\\
 & & {} + \left. 12x(x-1)(5y^2-5y+1)(2x^2(x-1)^2+y(y-1)(5x^2-5x+1))\right],
\end{eqnarray*}
and the exact solution is~\cite{chinosi:NAMFERM:1995}:
\begin{eqnarray*}
r_x &=& -y^3 (y-1)^3 x^2 (x-1)^2 (2x-1), \quad r_y=-x^3 (x-1)^3 y^2 (y-1)^2 (2y-1),\\
w &=& \frac{1}{3}x^3 (x-1)^3 y^3 (y-1)^3 - \frac{2t^2}{5(1-\nu)}\left[y^3 (y-1)^3 x(x-1)(5x^2-5x+1)\right.\\
& & {} +\left.x^3 (x-1)^3 y(y-1)(5y^2-5y+1)\right].
\end{eqnarray*}

\begin{figure}[!tbhp]
\centering
\mbox{
\subfigure[]{\label{fig:meshes_chinosi_prob_a} \epsfig{file = 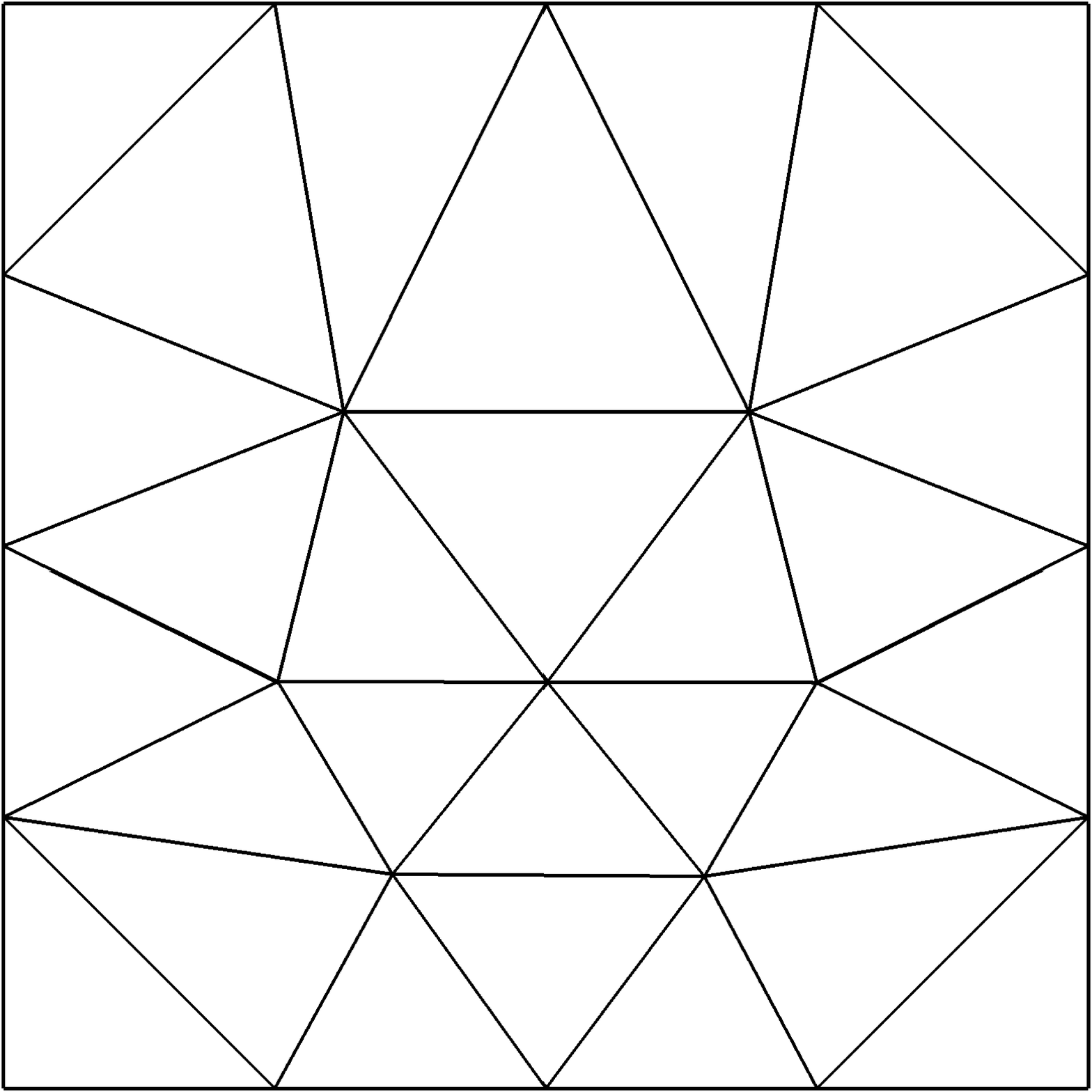, width = 0.2\textwidth}}
\subfigure[]{\label{fig:meshes_chinosi_prob_b} \epsfig{file = 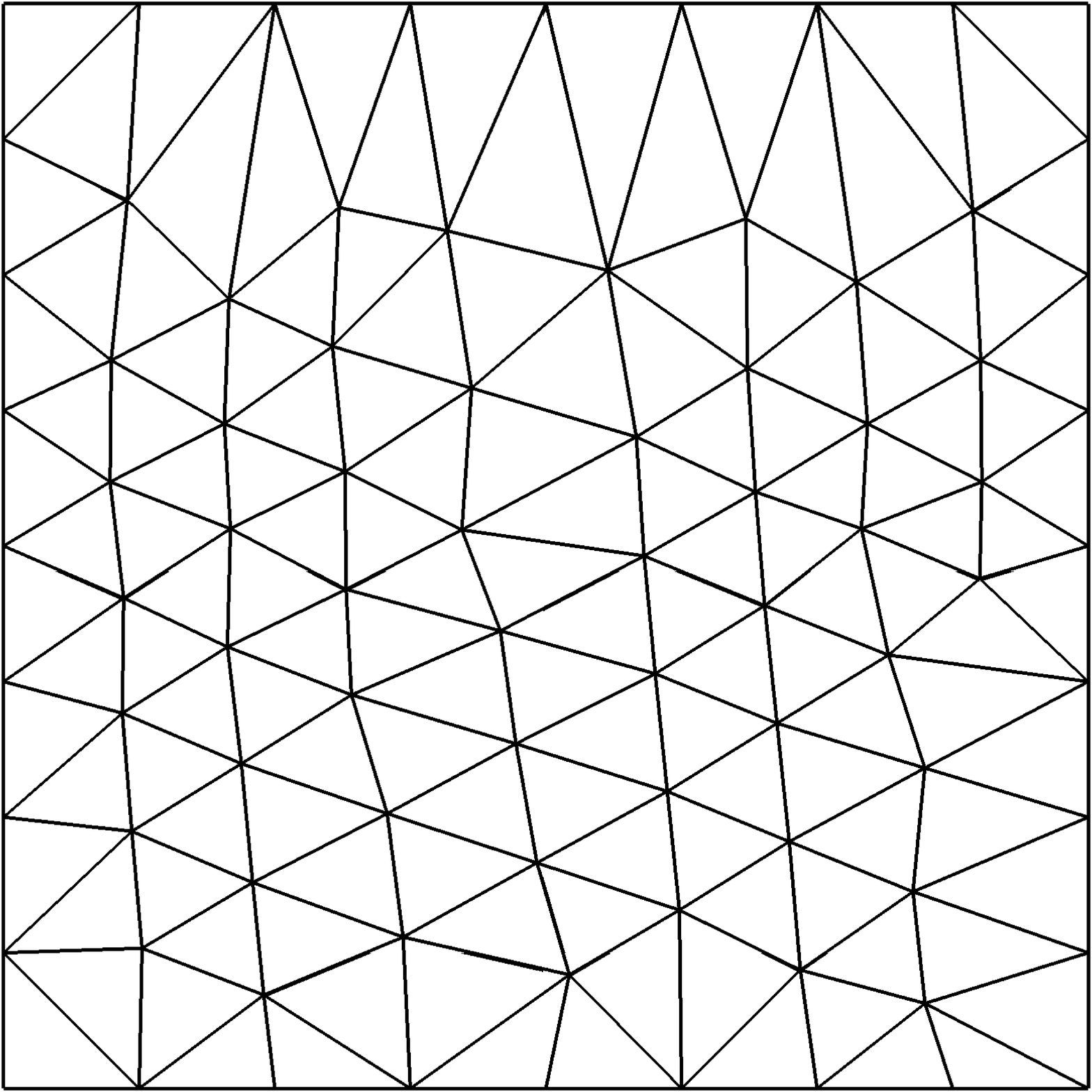, width = 0.2\textwidth}}
\subfigure[]{\label{fig:meshes_chinosi_prob_c} \epsfig{file = 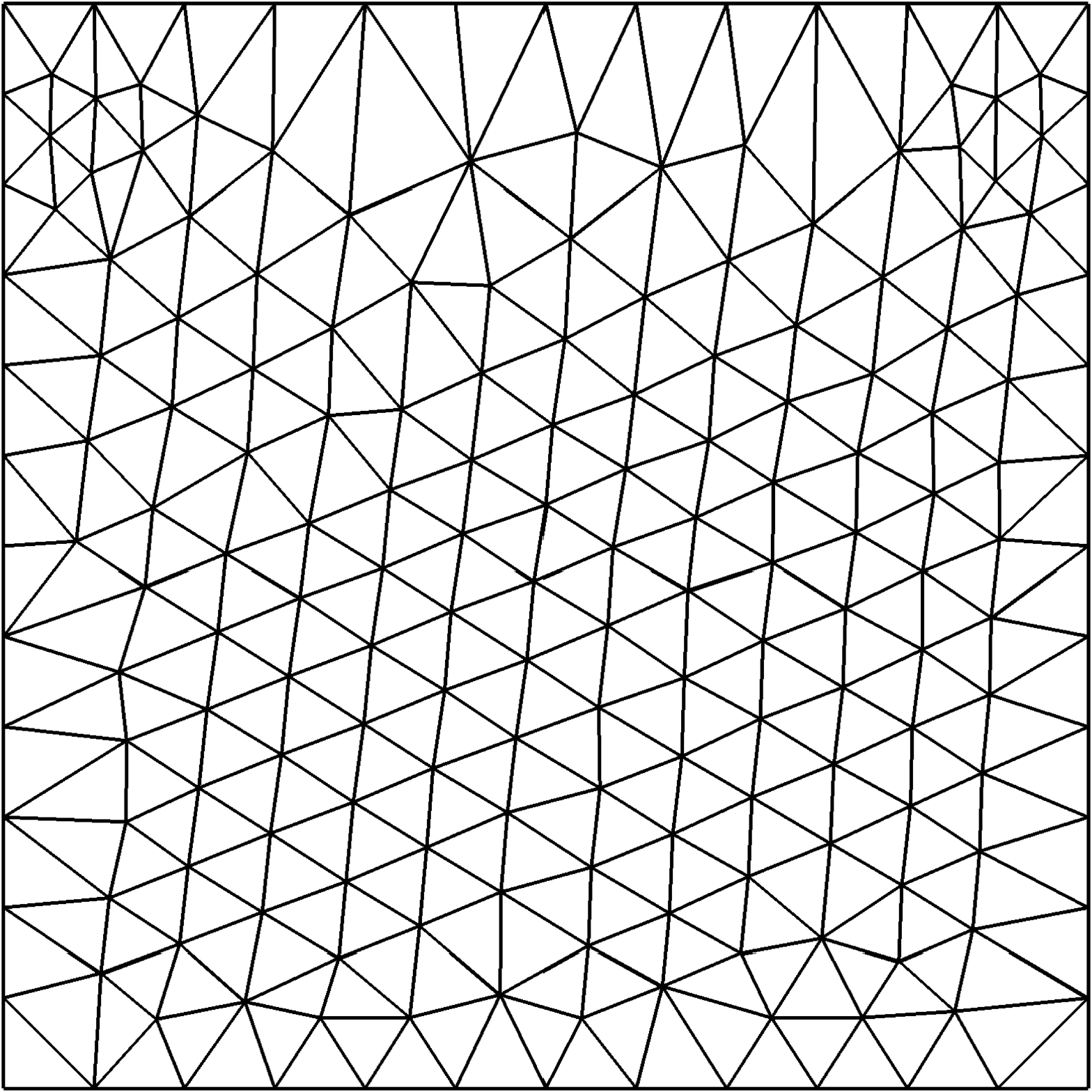, width = 0.2\textwidth}}
\subfigure[]{\label{fig:meshes_chinosi_prob_d} \epsfig{file = 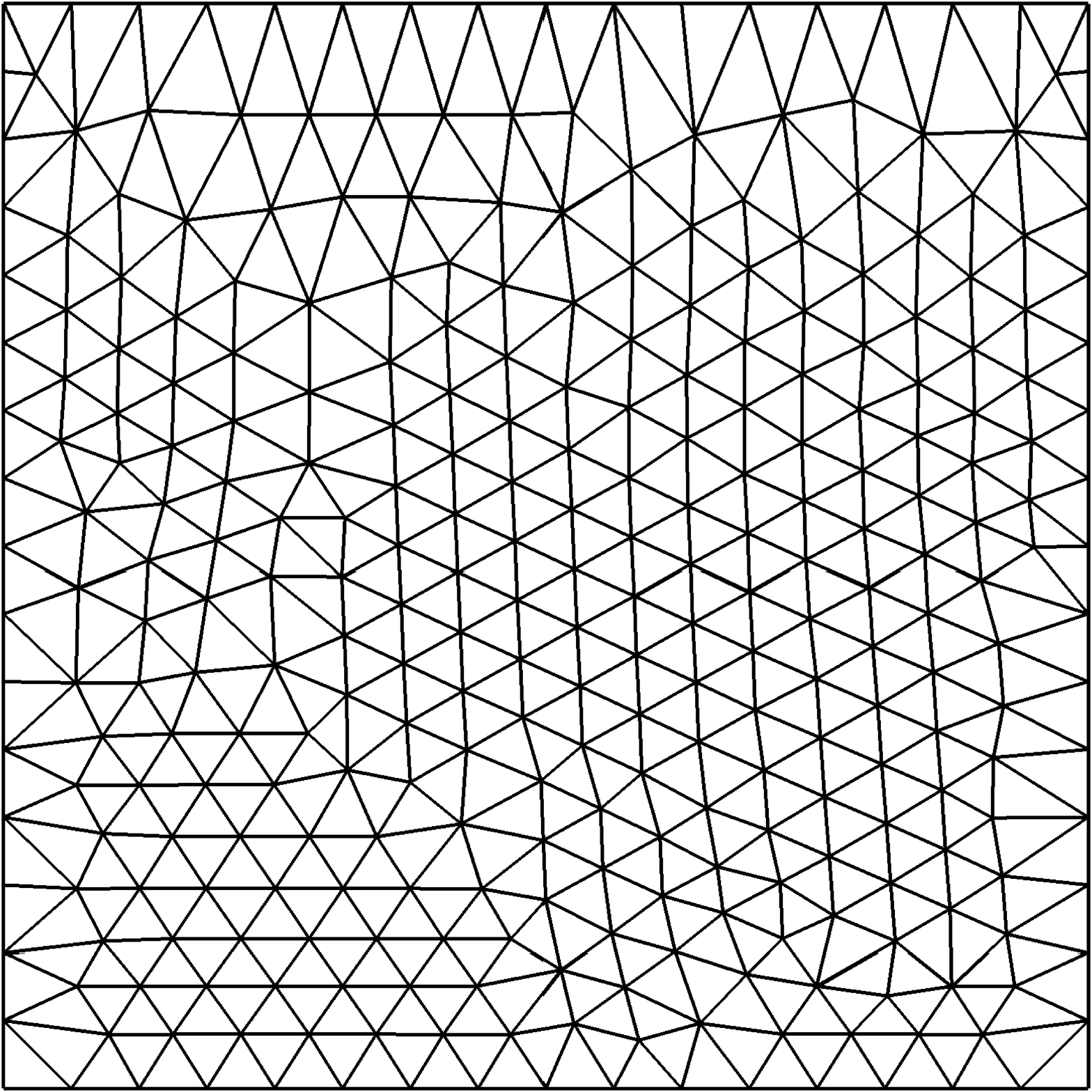, width = 0.2\textwidth}}
}
\mbox{
\subfigure[]{\label{fig:meshes_chinosi_prob_e} \epsfig{file = 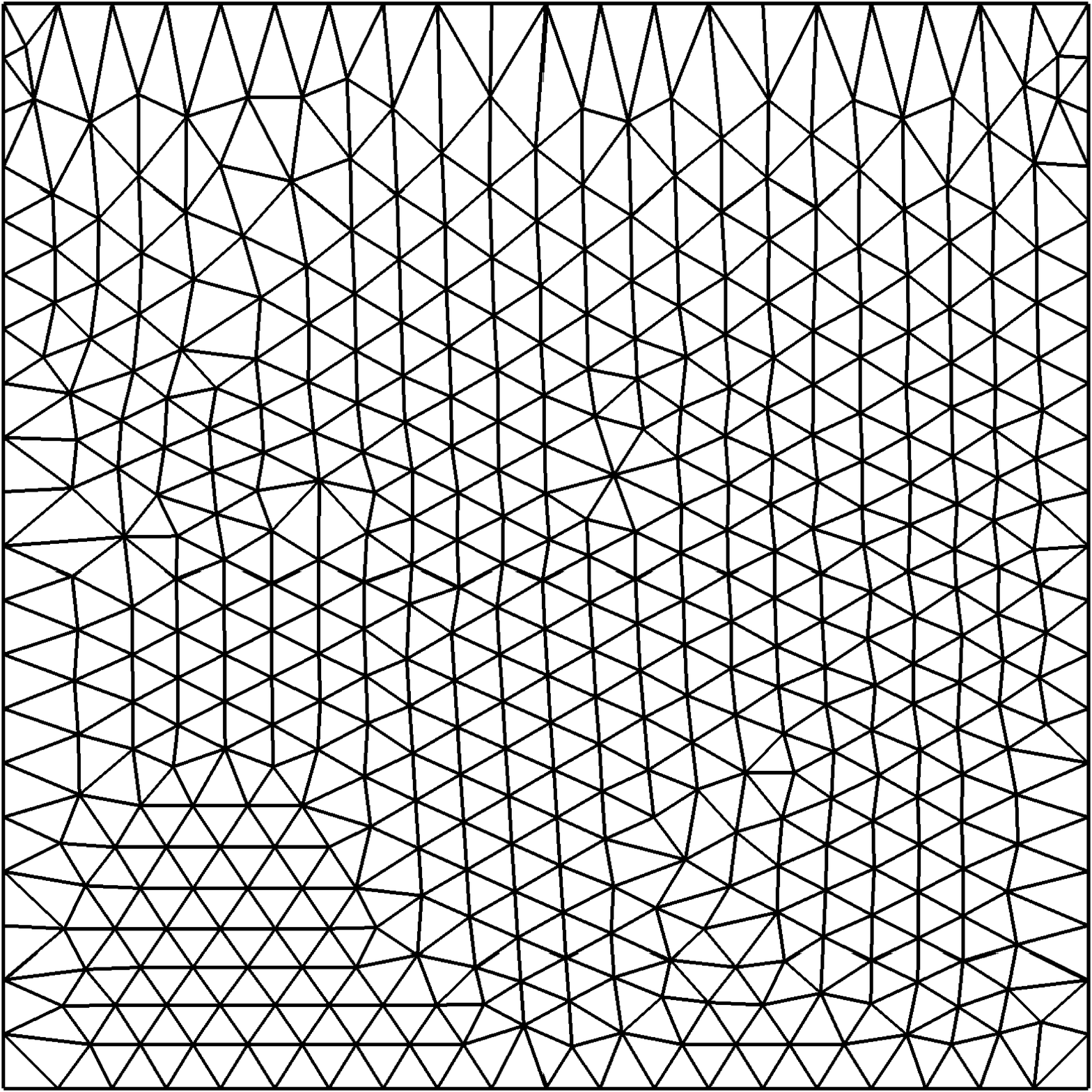, width = 0.2\textwidth}}
\subfigure[]{\label{fig:meshes_chinosi_prob_f} \epsfig{file = 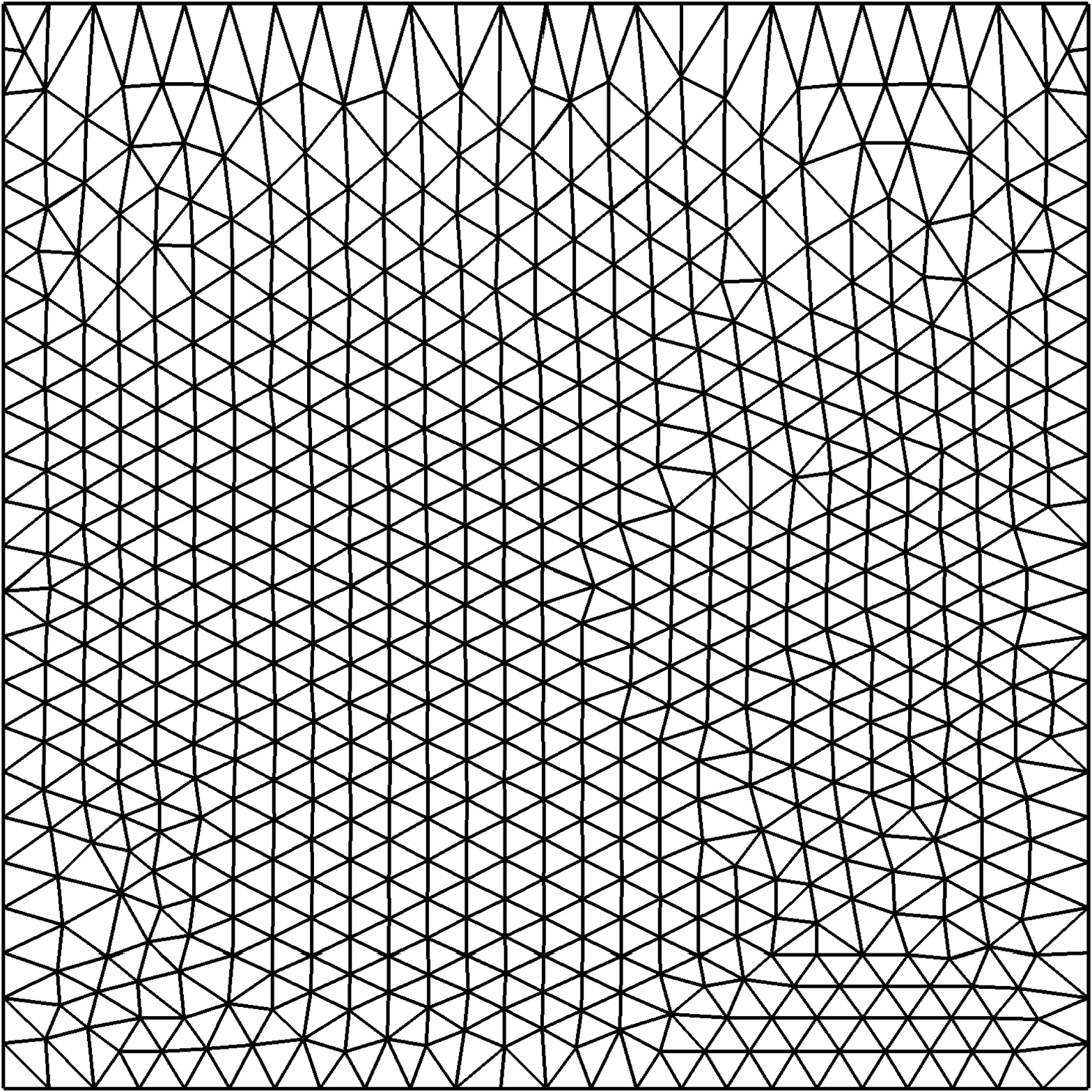, width = 0.2\textwidth}}
\subfigure[]{\label{fig:meshes_chinosi_prob_g} \epsfig{file = 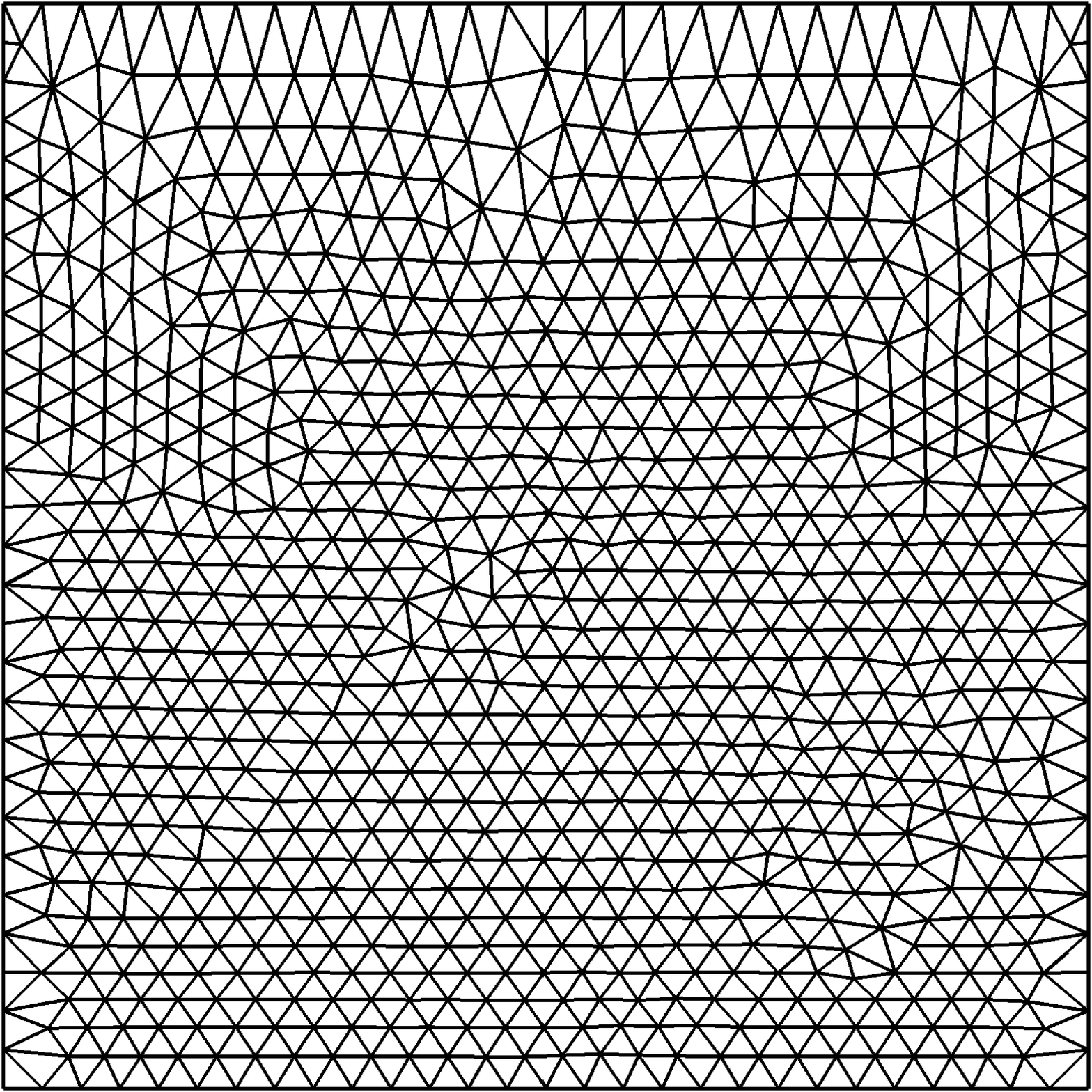, width = 0.2\textwidth}}
\subfigure[]{\label{fig:meshes_chinosi_prob_h} \epsfig{file = 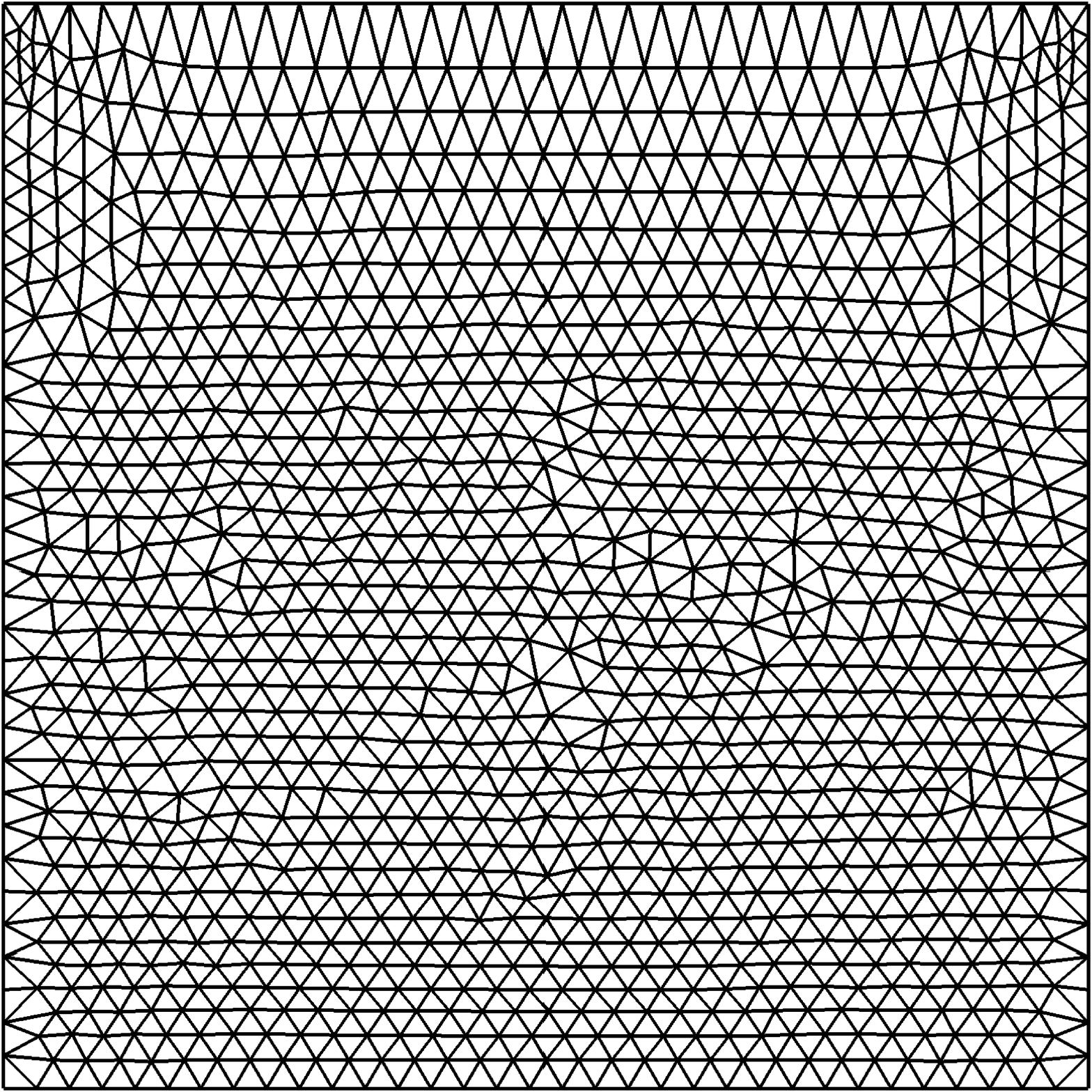, width = 0.2\textwidth}}}
\caption{Integration meshes for the square plate subjected to a nonuniform load problem.}
\label{fig:meshes_chinosi_prob}
\end{figure}

The integration meshes that are considered for this problem are shown in~\fref{fig:meshes_chinosi_prob}.

\begin{figure}[!tbhp]
\centering
\mbox{
\subfigure[]{\label{fig:norms_chinosi_prob_a} \epsfig{file = 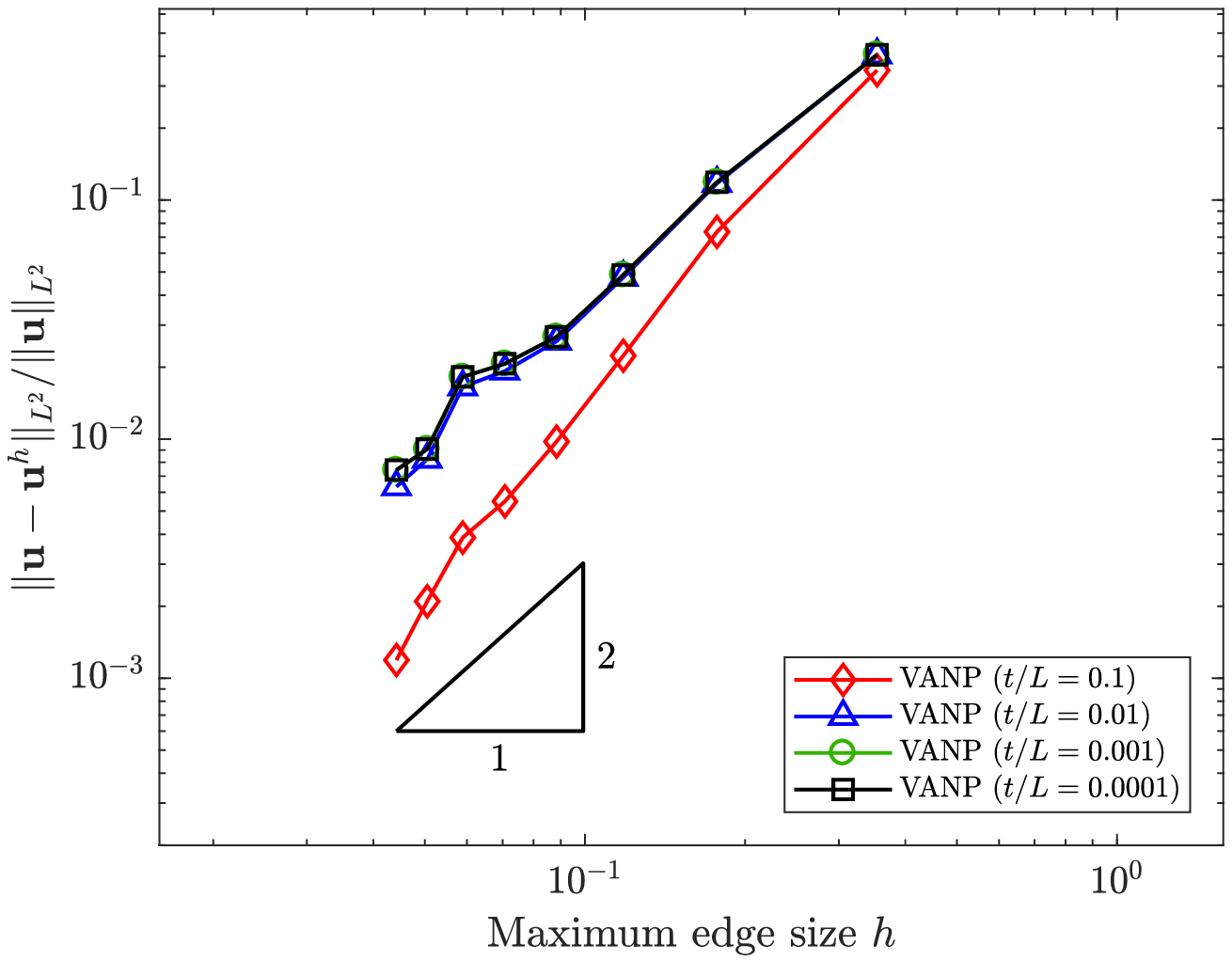, width = 0.5\textwidth}}
\subfigure[]{\label{fig:norms_chinosi_prob_b} \epsfig{file = 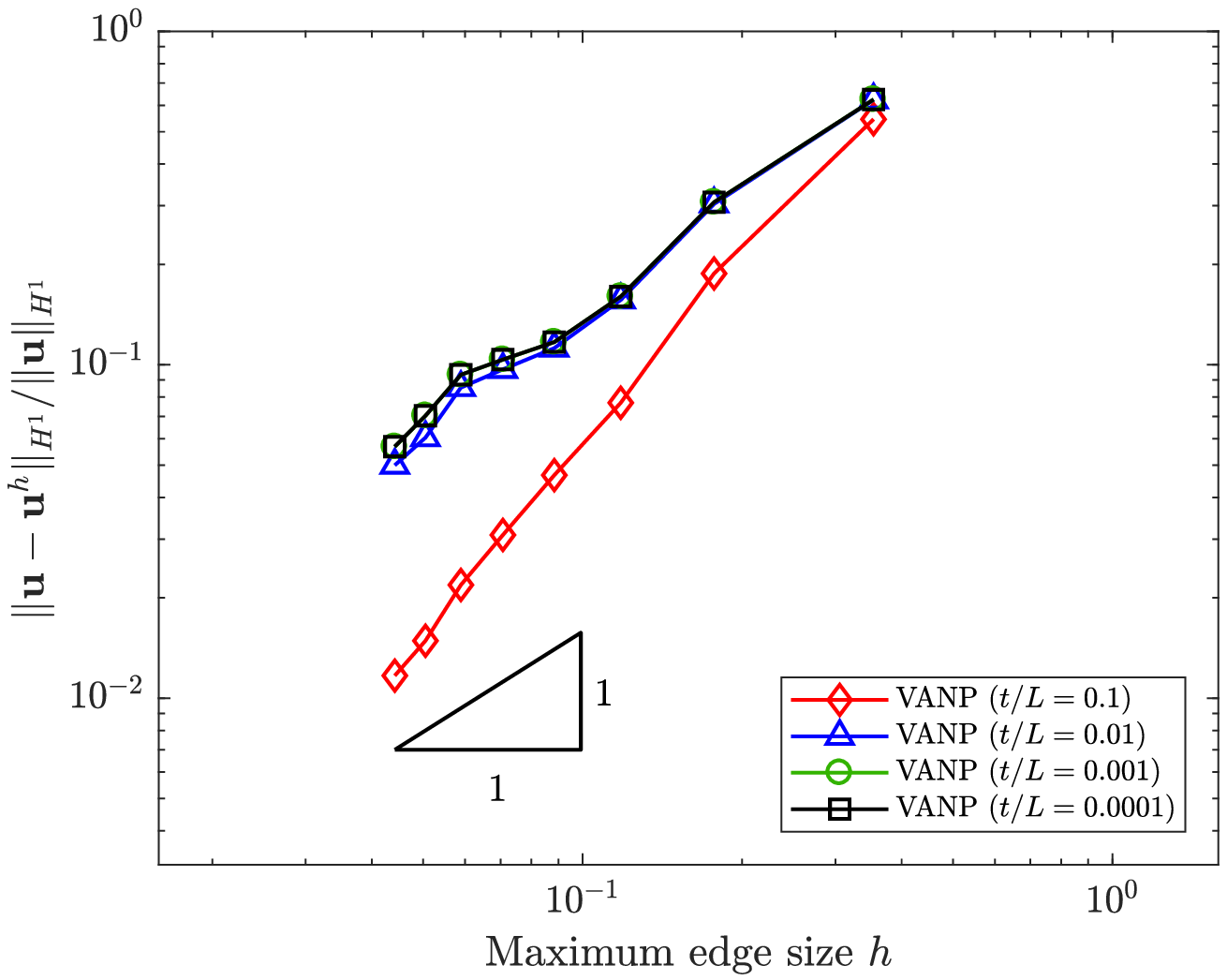, width = 0.5\textwidth}}
}
\caption{Rates of convergence for the square plate subjected to a nonuniform load. (a) $L^2$-norm of the error and (b) $H^1$-seminorm of the error for several values of $t/L$. The \texttt{VANP} method delivers the optimal rates of convergence for $t/L=\{0.01,\, 0.001,\, 0.0001\}$ and above the optimal for $t/L=0.1$.}
\label{fig:norms_chinosi_prob}
\end{figure}

The convergence of the \texttt{VANP} approach as the integration mesh is refined is studied for the following normalized thicknesses: $t/L=\{0.1,\, 0.01,\, 0.001,\, 0.0001\}$. The convergence rates are shown in~\fref{fig:norms_chinosi_prob}, where it is observed that the optimal rates of convergence, 2 and 1, are delivered by the \texttt{VANP} method in both the $L^2$-norm and the $H^1$-seminorn of the error, respectively, for the normalized thicknesses $t/L=\{0.01,\, 0.001,\, 0.0001\}$. On the other hand, the convergence rates for $t/L=0.1$ (the thicker plate case) are above the optimal.

\begin{figure}[!tbhp]
\centering
\mbox{
\subfigure[]{\label{fig:norms_integration_chinosi_prob_a} \epsfig{file = 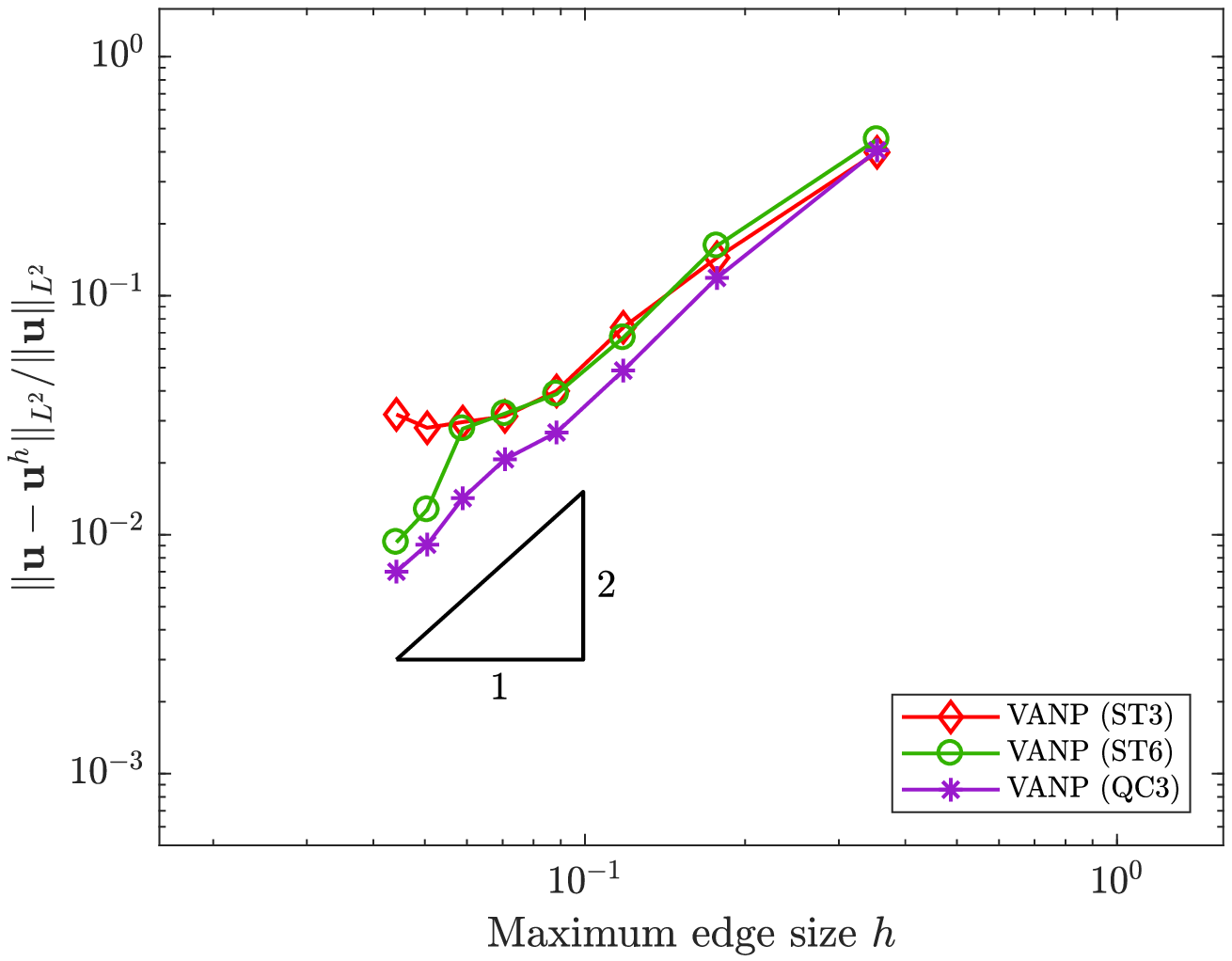, width = 0.5\textwidth}}
\subfigure[]{\label{fig:norms_integration_chinosi_prob_b} \epsfig{file = 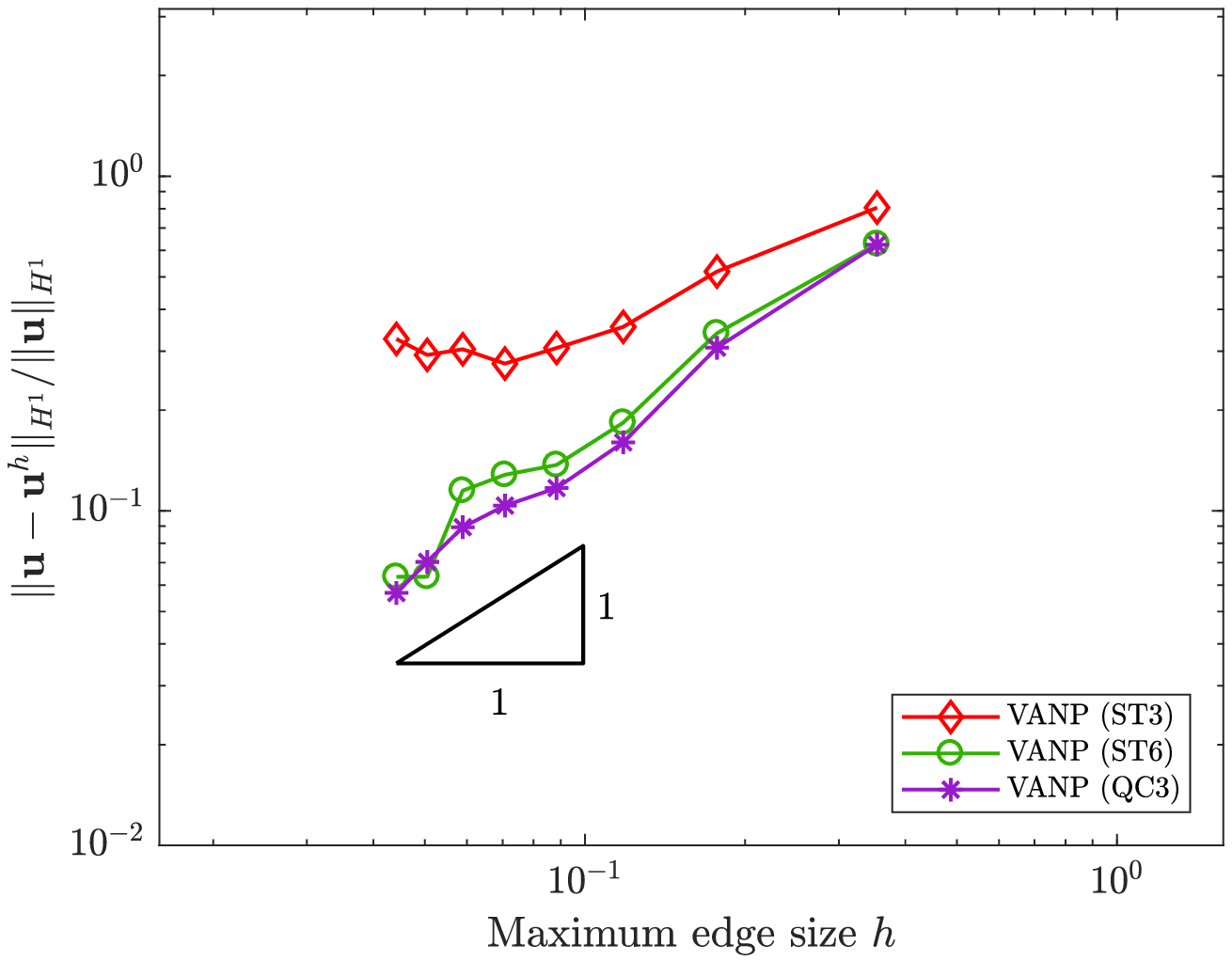, width = 0.5\textwidth}}
}
\caption{Influence of the numerical integration on the \texttt{VANP} convergence rates. (a) $L^2$-norm of the error and (b) $H^1$-seminorm of the error for the ST3, ST6 and QC3 integration schemes. The ST3 scheme fails to converge, the ST6 improves the convergence and the QC3 provides the optimal convergence.}
\label{fig:norms_integration_chinosi_prob}
\end{figure}

To illustrate the influence of the numerical integration in the accuracy of the \texttt{VANP} formulation, we compare the numerical solutions using three integration rules on the triangular cell: the 3-point standard Gauss rule (ST3), the 6-point standard Gauss rule (ST6) and the default \texttt{VANP}'s integration scheme (QC3) that was developed in Section~\ref{sec:numericalintegration}. For this test, the normalized thickness $t/L=0.0001$ is considered. \fref{fig:norms_integration_chinosi_prob} depicts the convergence curves for each of these integration schemes. As can be seen, the ST3 scheme fails to converge in both the $L^2$-norm and the $H^1$-seminorn of the error. Even though the ST6 scheme exhibits much better convergence properties than the ST3 scheme, the optimal performance is observed for the QC3 scheme.

\subsection{Parallelogram plate subjected to a uniform load}
\label{sec:numexamples_morley}

\begin{figure}[!tbhp]
\centering
\epsfig{file = 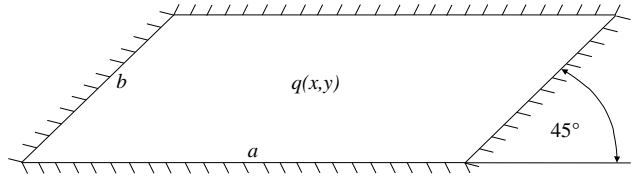, width = 0.55\textwidth}
\caption{Parallelogram plate subjected to a uniform load.}
\label{fig:morley_prob}
\end{figure}

This example is tailored to show the performance of the \texttt{VANP} formulation when distorted integration meshes are used. The problem consists in a parallelogram plate of unit thickness that is clamped along the entire boundary and subjected to a uniform load, as shown in~\fref{fig:morley_prob}. The problem parameters are set as follows: $a=200$ in, $b=100$ in, $q=100$ psi. The plates' material parameters are $E_\textrm{Y}=10.92\times 10^6$ psi and $\nu=0.3$. The analytical reference value for the maximum transverse displacement can be found in Ref.~\cite{mansfield:1989:BSP}.

\begin{figure}[!tbhp]
\centering
\mbox{
\subfigure[]{\label{fig:meshes_morley_prob_a} \epsfig{file = 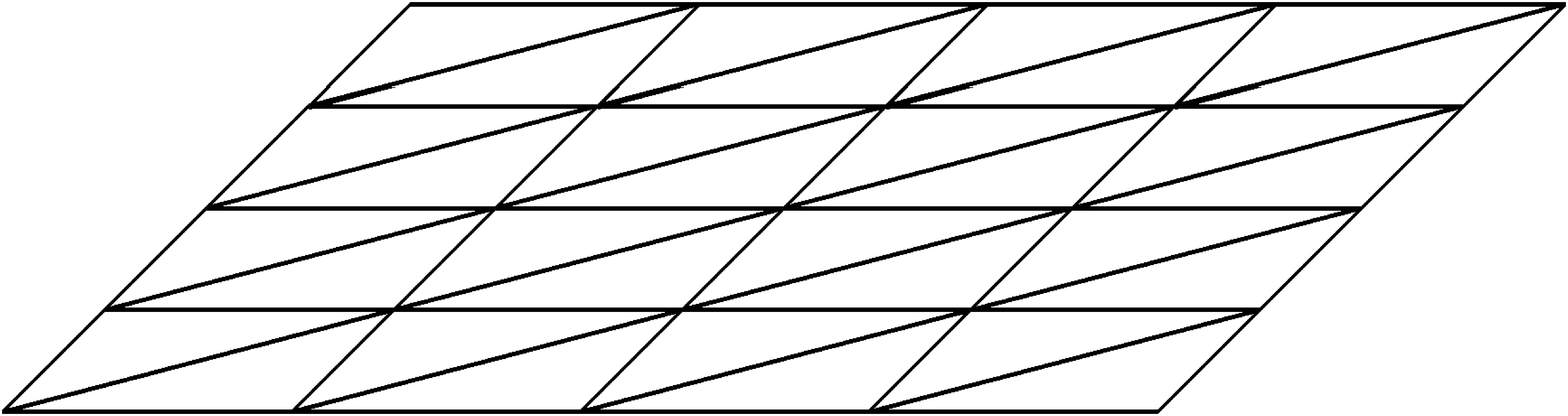, width = 0.5\textwidth}}
\subfigure[]{\label{fig:meshes_morley_prob_b} \epsfig{file = 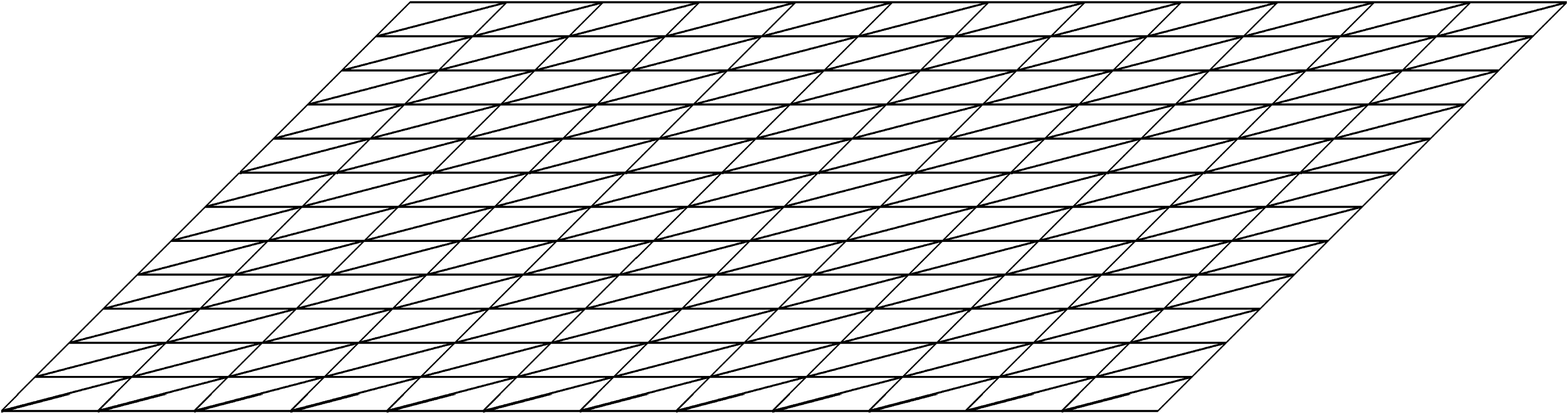, width = 0.5\textwidth}}
}
\mbox{
\subfigure[]{\label{fig:meshes_morley_prob_c} \epsfig{file = 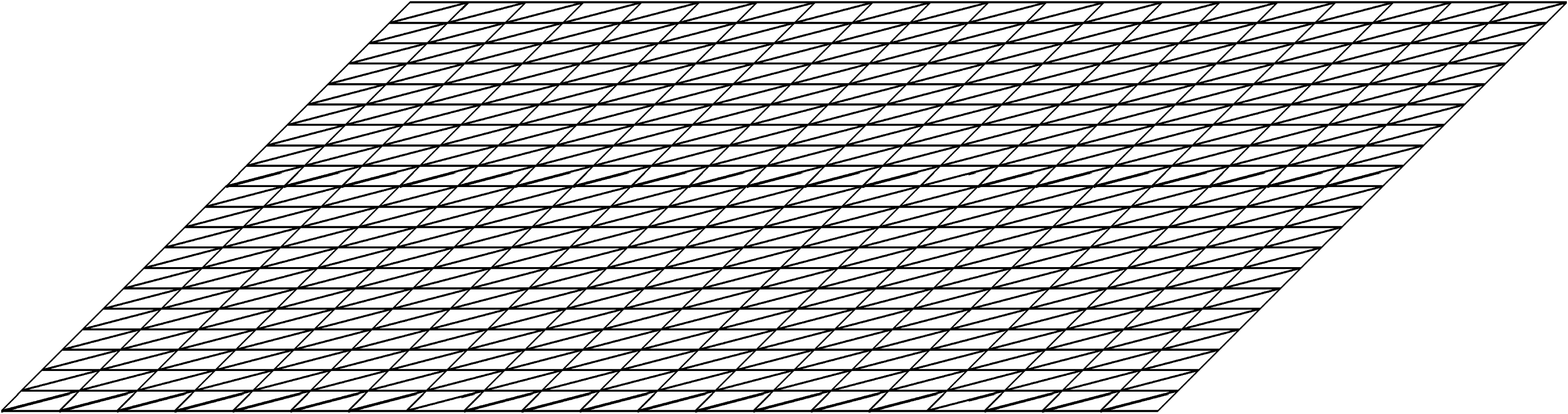, width = 0.5\textwidth}}
\subfigure[]{\label{fig:meshes_morley_prob_d} \epsfig{file = 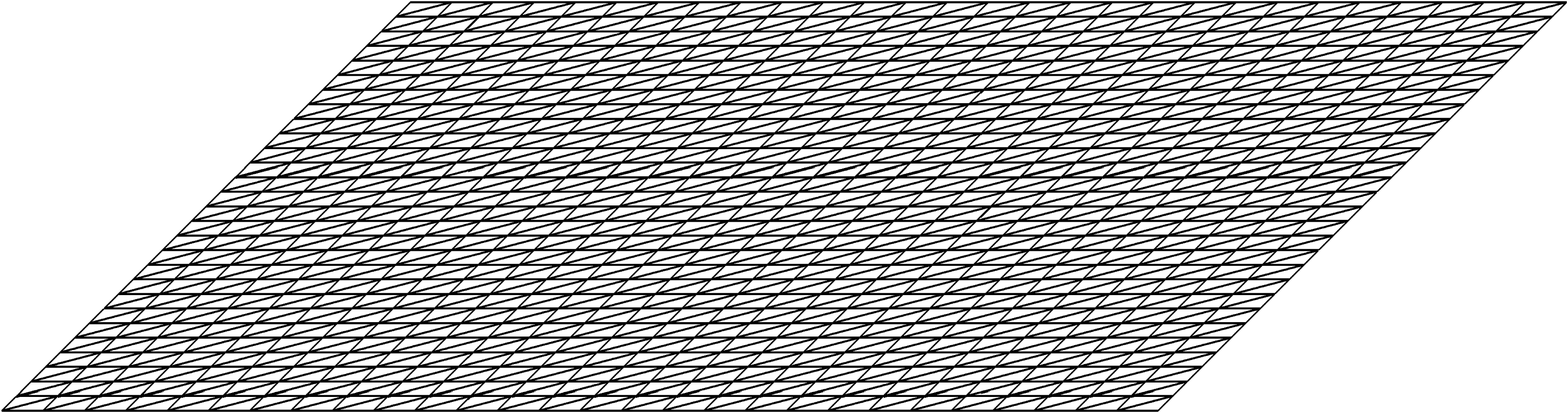, width = 0.5\textwidth}}
}
\subfigure[]{\label{fig:meshes_morley_prob_e} \epsfig{file = 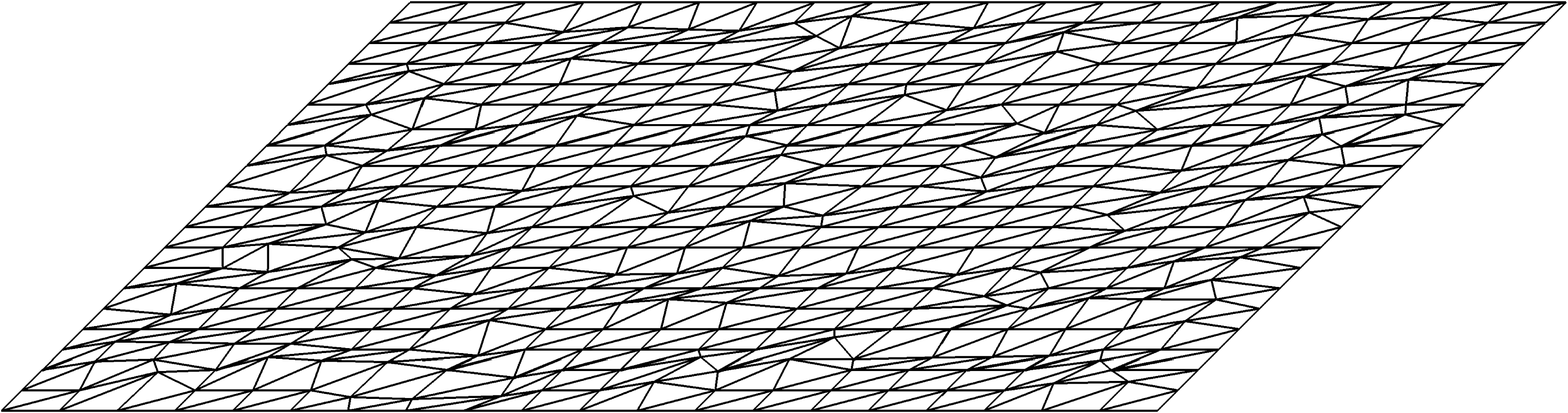, width = 0.5\textwidth}}
\caption{Integration meshes for the parallelogram plate problem.}
\label{fig:meshes_morley_prob}
\end{figure}

The integration meshes that are considered for this problem are depicted in~\fref{fig:meshes_morley_prob}.

\begin{figure}[!tbhp]
\centering
\mbox{
\subfigure[]{\label{fig:disp_morley_prob_a} \epsfig{file = 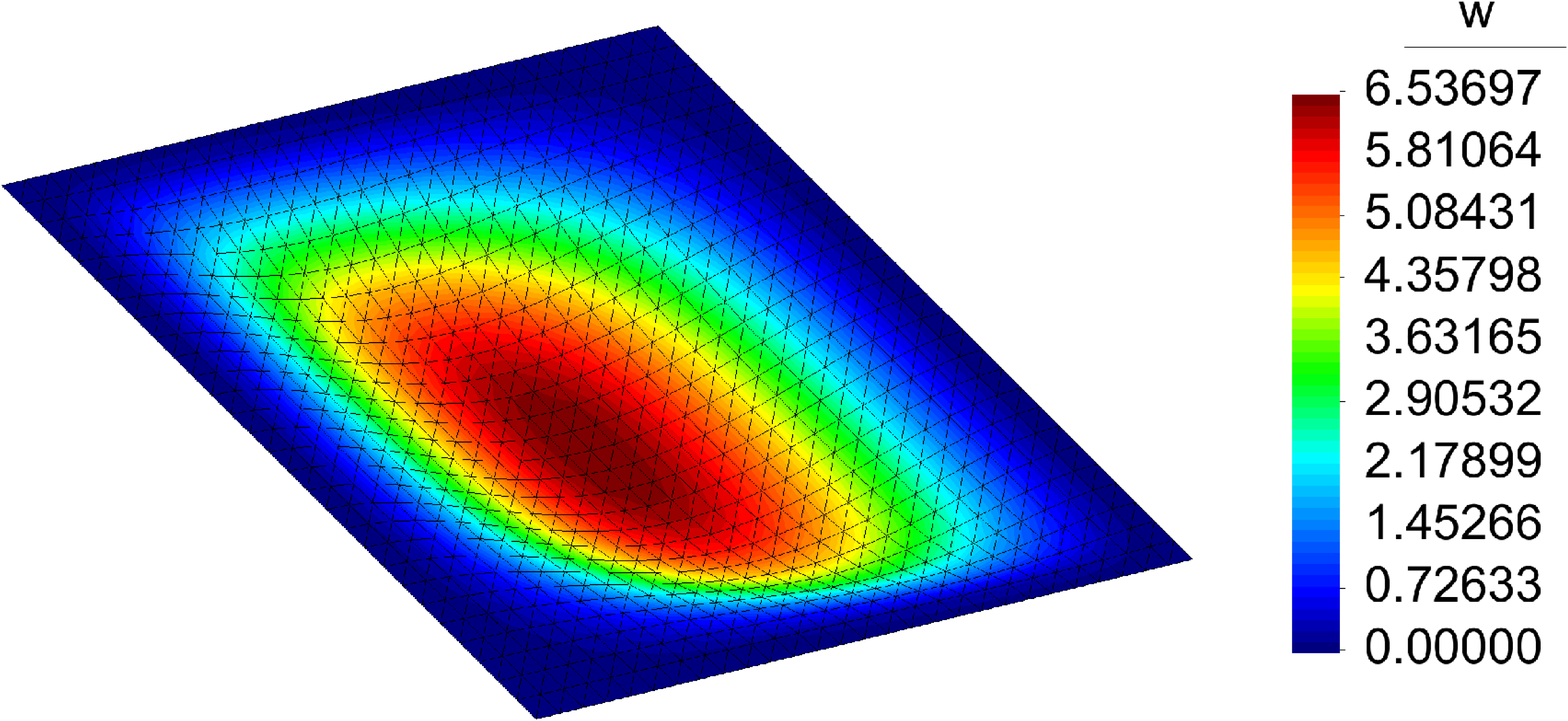, width = 0.5\textwidth}}
\subfigure[]{\label{fig:disp_morley_prob_b} \epsfig{file = 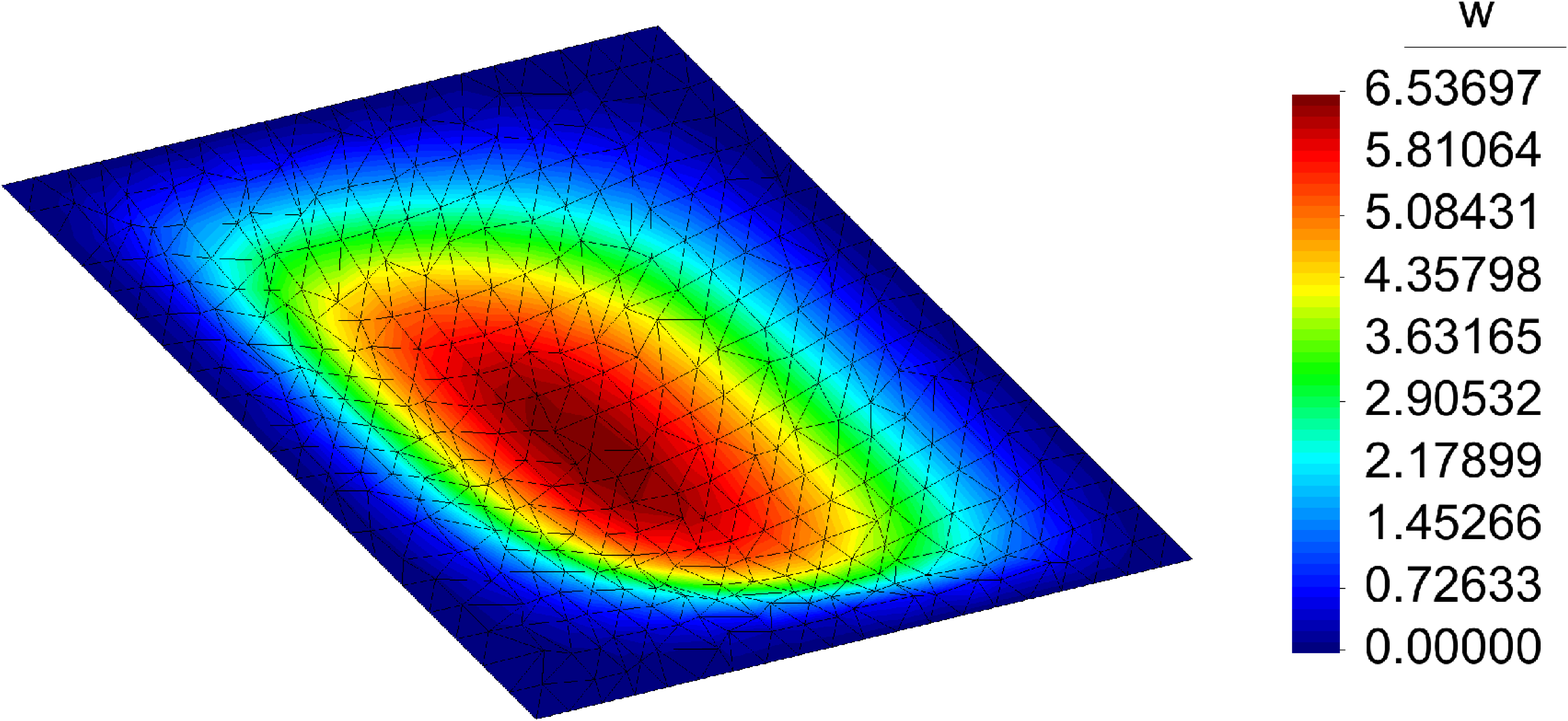, width = 0.5\textwidth}}
}
\caption{Transverse displacement solution for the parallelogram plate problem when the integration meshes (a) \ref{fig:meshes_morley_prob_d} and (b) \ref{fig:meshes_morley_prob_e} are used.}
\label{fig:disp_morley_prob}
\end{figure}

The transverse displacement field solution for the integration meshes shown in Figs.~\ref{fig:meshes_morley_prob_d} and~\ref{fig:meshes_morley_prob_e} are presented in~\fref{fig:disp_morley_prob}. Table~\ref{tab:morley45} summarizes the maximum transverse displacement (located at the center of the plate) that is obtained for each of the integration meshes considered. The table also provides the analytical reference solution. It is observed that the numerical solutions are close to the analytical reference solution for all the integration meshes considered.

\begin{table}[]
\centering
\caption{Transverse displacement for the parallelogram plate problem.}
\label{tab:morley45}
\begin{tabular}{|ccccc|c|}
\hline
\multicolumn{5}{|c|}{Integration mesh}        &                                    \\ \hline
\ref{fig:meshes_morley_prob_a}  & \ref{fig:meshes_morley_prob_b} & \ref{fig:meshes_morley_prob_c} & \ref{fig:meshes_morley_prob_d}  &\ref{fig:meshes_morley_prob_e} & \multicolumn{1}{c|}{Ref. solution} \\ \hline
6.51227 & 6.52950 & 6.53558 & 6.53697 & 6.52780 & \multicolumn{1}{c|}{6.52000} \\ \hline
\end{tabular}
\end{table}

\begin{figure}[!tbhp]
\centering
\epsfig{file = 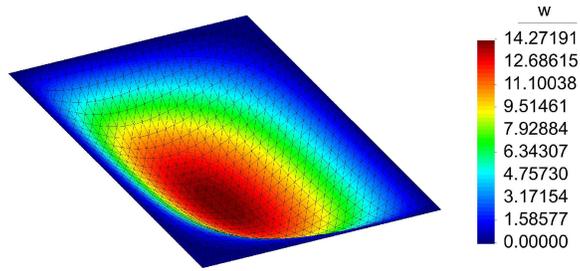, width = 0.5\textwidth}
\caption{Parallelogram plate subjected to a uniform load. The use of the 3-point standard Gauss rule (ST3) on the integration mesh \ref{fig:meshes_morley_prob_d} leads to an erroneous transverse deflection field solution.}
\label{fig:disp_morley_prob_std_integration}
\end{figure}

Finally, once again we show the importance of the QC3 integration scheme that was developed for the \texttt{VANP} formulation. \fref{fig:disp_morley_prob_std_integration} provides the transverse displacement field solution for the integration mesh shown in~\fref{fig:meshes_morley_prob_d} when the 3-point standard Gauss rule (ST3) is used. A comparison between the results shown in~\fref{fig:disp_morley_prob_a} and~\fref{fig:disp_morley_prob_std_integration} reveals that the ST3 integration scheme leads to an erroneous transverse deflection field.

\subsection{Performance of the scaled transverse shear stress solution}
\label{sec:numexamples_shearstress}

\begin{figure}[!tbhp]
\centering
\mbox{
\subfigure[]{\label{fig:meshes_chinosi_prob_shearstress_a} \epsfig{file = 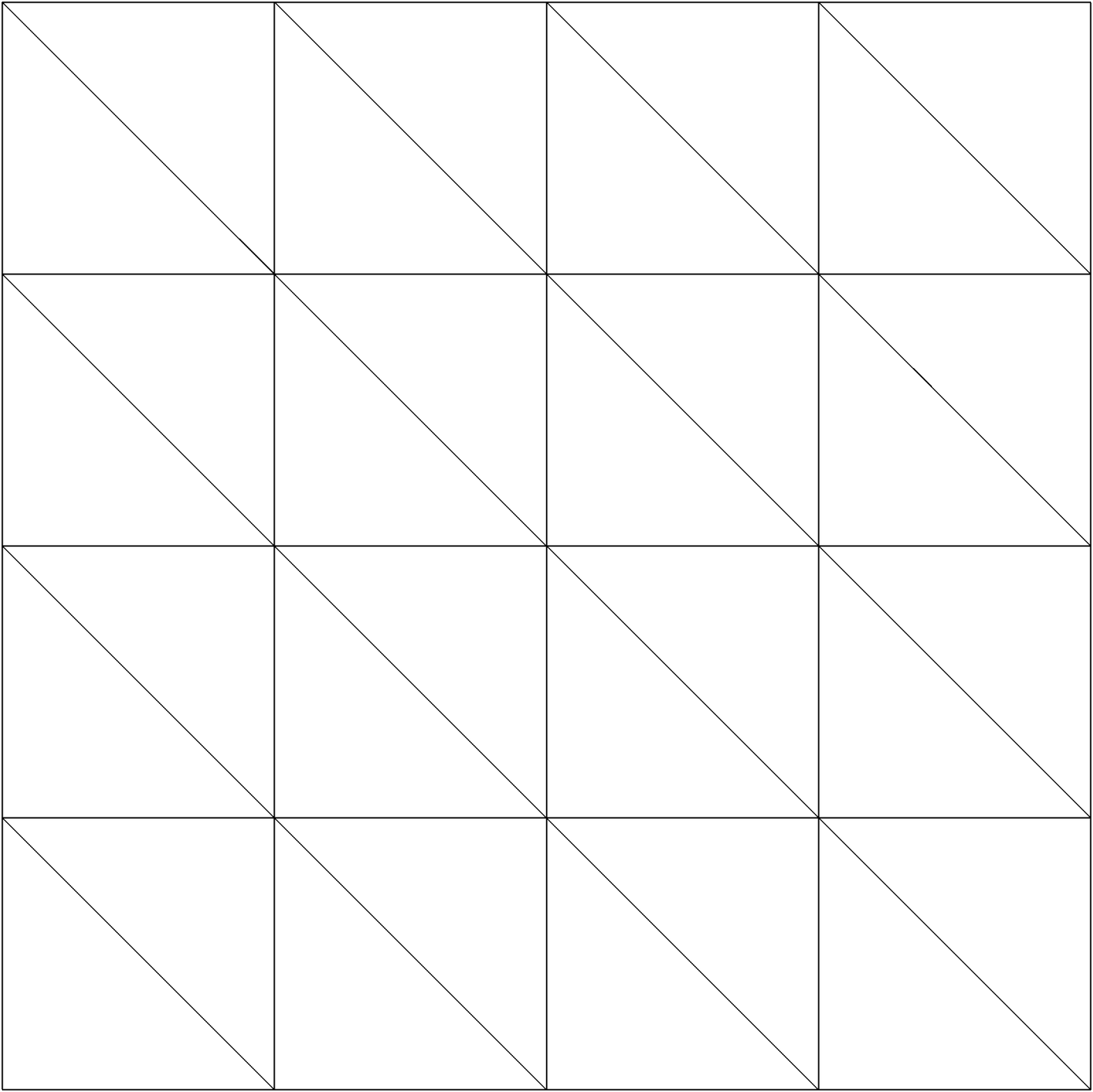, width = 0.2\textwidth}}
\subfigure[]{\label{fig:meshes_chinosi_prob_shearstress_b} \epsfig{file = 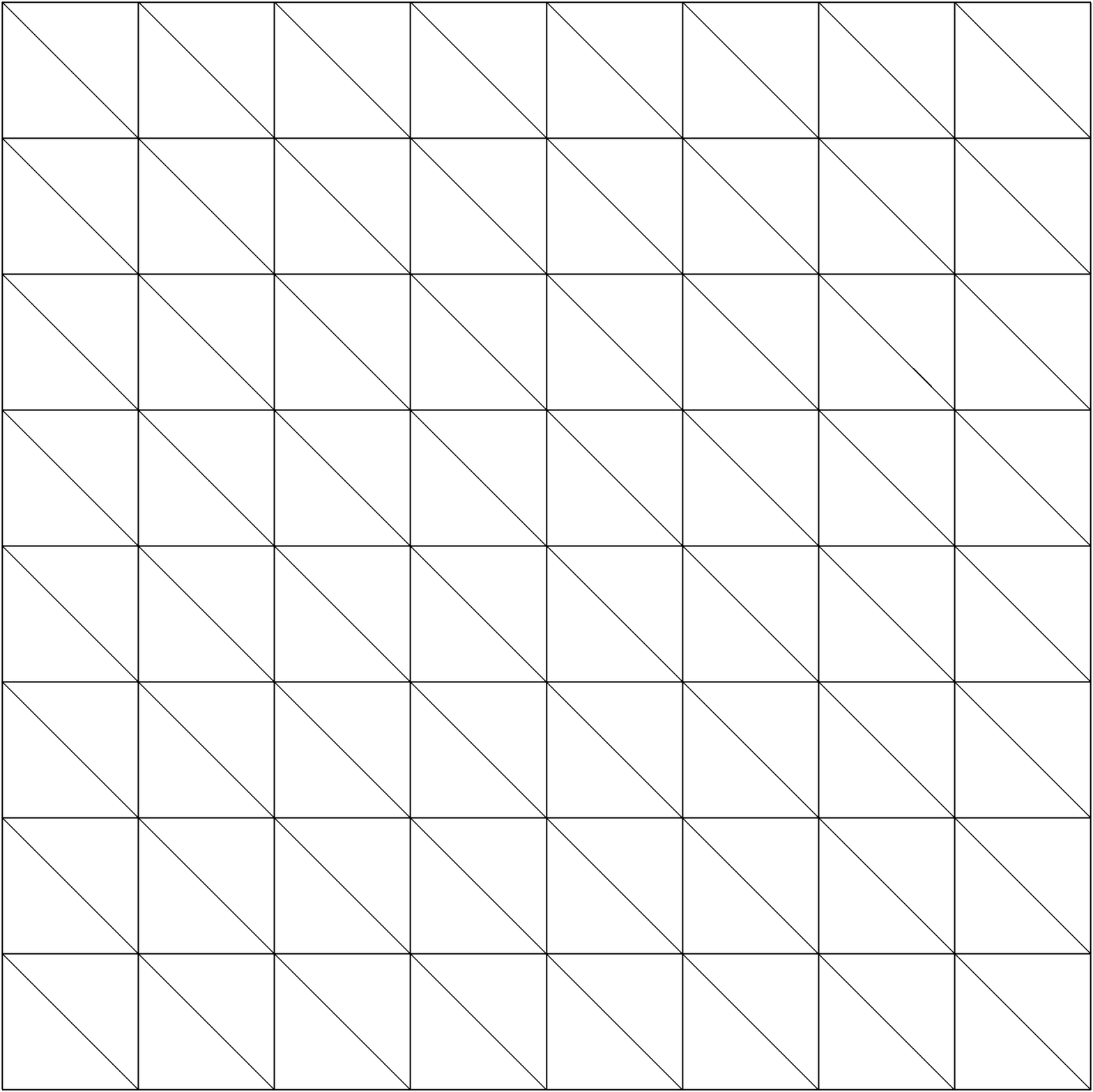, width = 0.2\textwidth}}
\subfigure[]{\label{fig:meshes_chinosi_prob_shearstress_c} \epsfig{file = 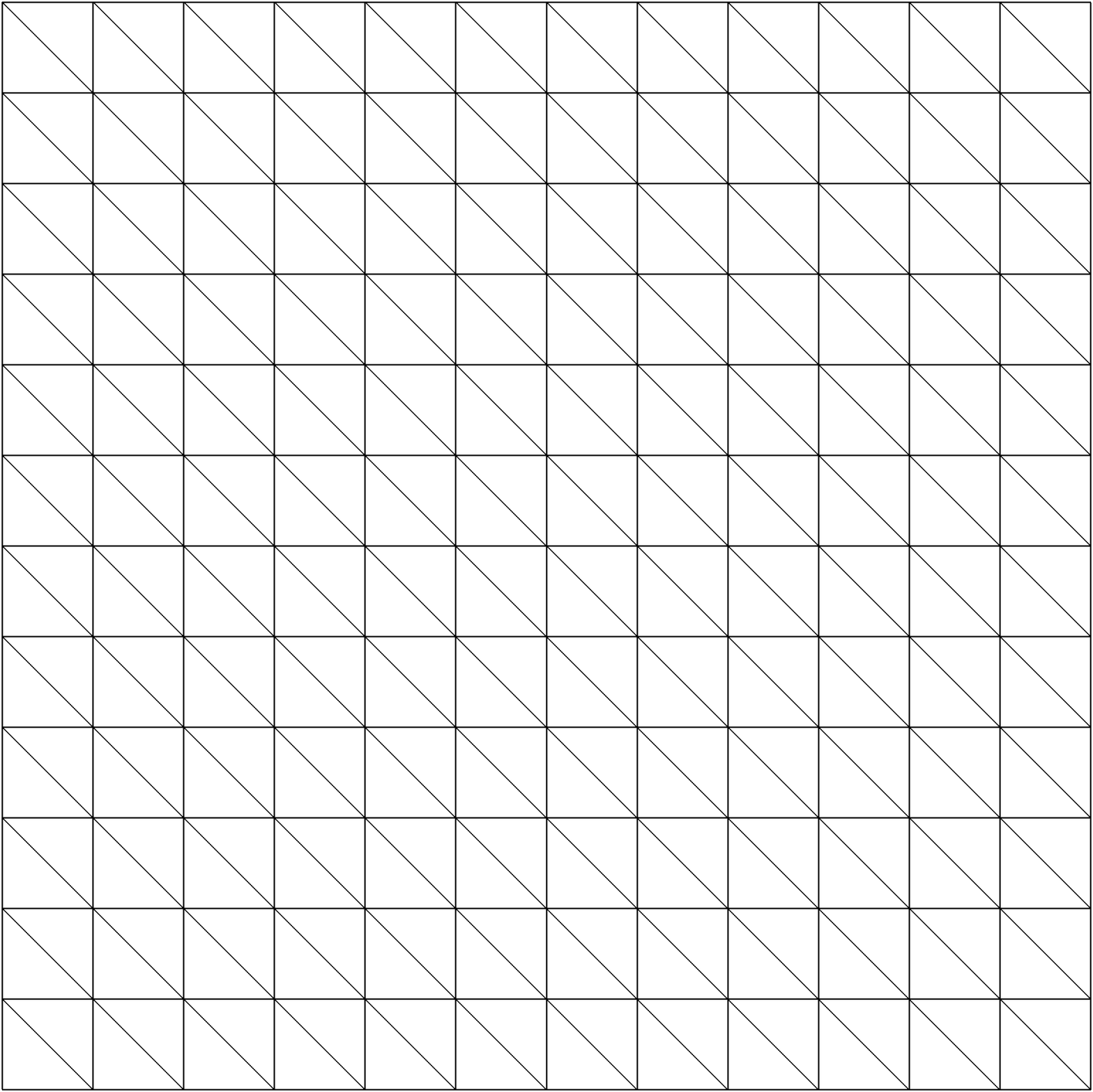, width = 0.2\textwidth}}
\subfigure[]{\label{fig:meshes_chinosi_prob_shearstress_d} \epsfig{file = 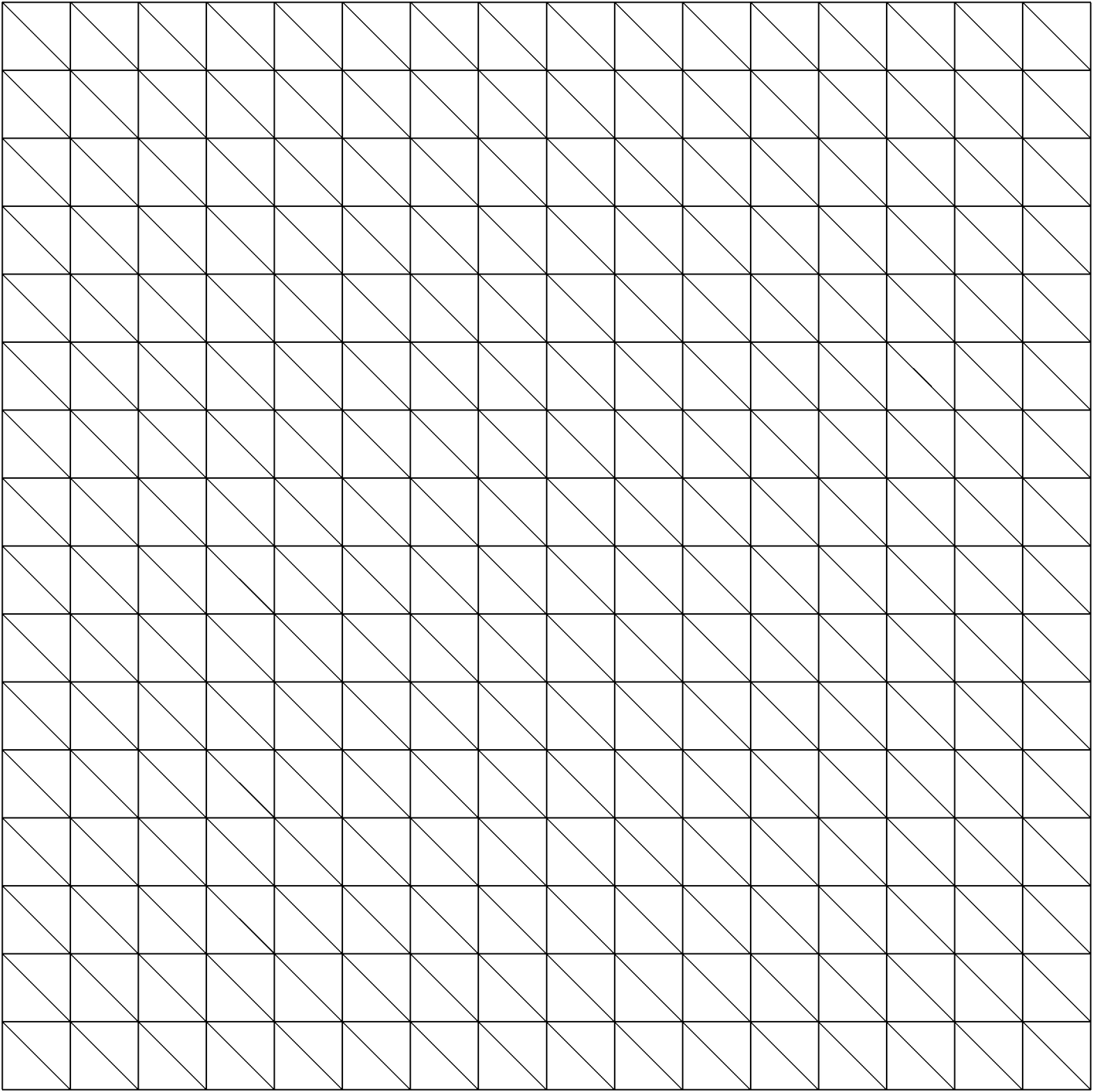, width = 0.2\textwidth}}
}
\mbox{
\subfigure[]{\label{fig:meshes_chinosi_prob_shearstress_e} \epsfig{file = 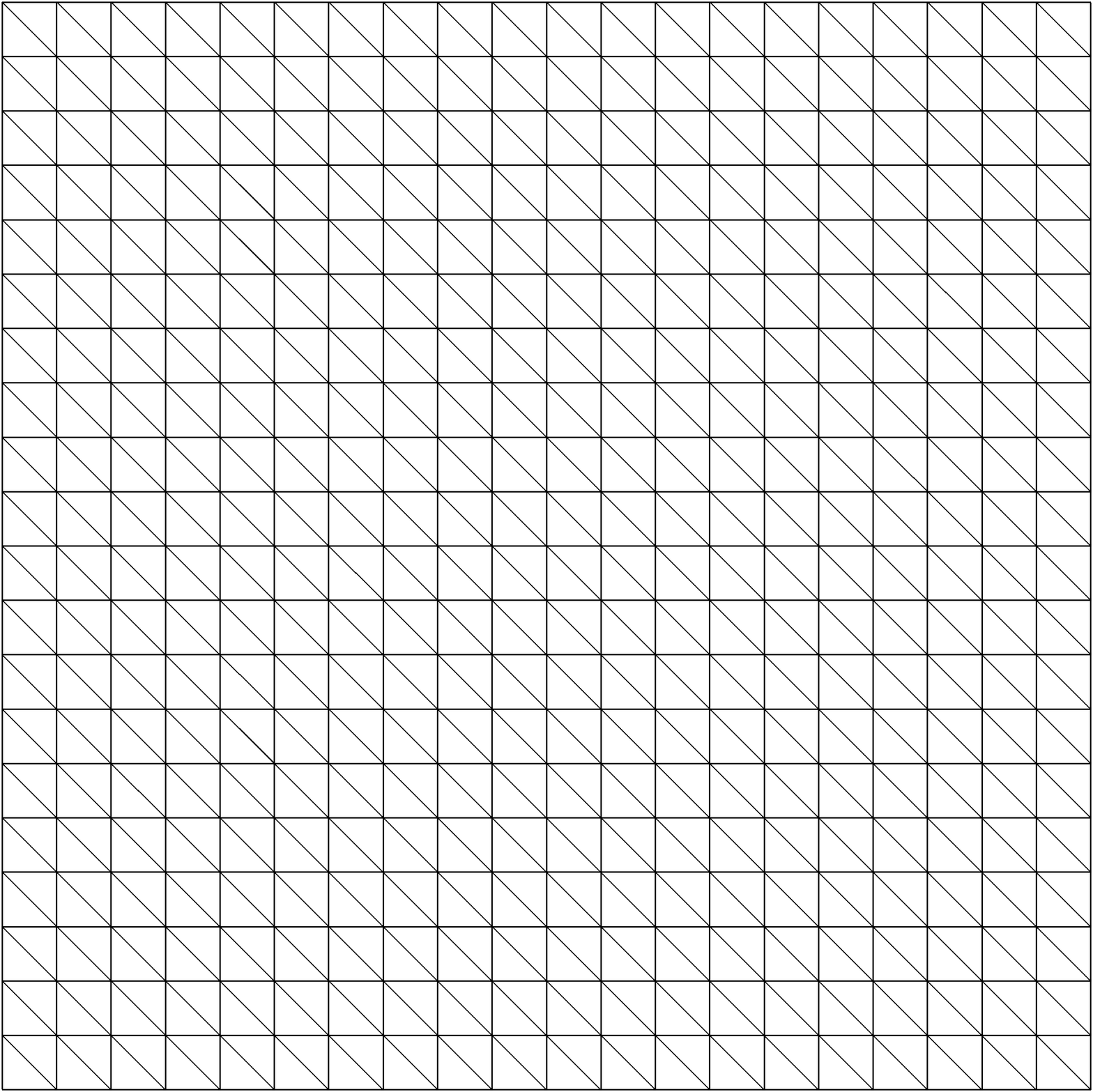, width = 0.2\textwidth}}
\subfigure[]{\label{fig:meshes_chinosi_prob_shearstress_f} \epsfig{file = 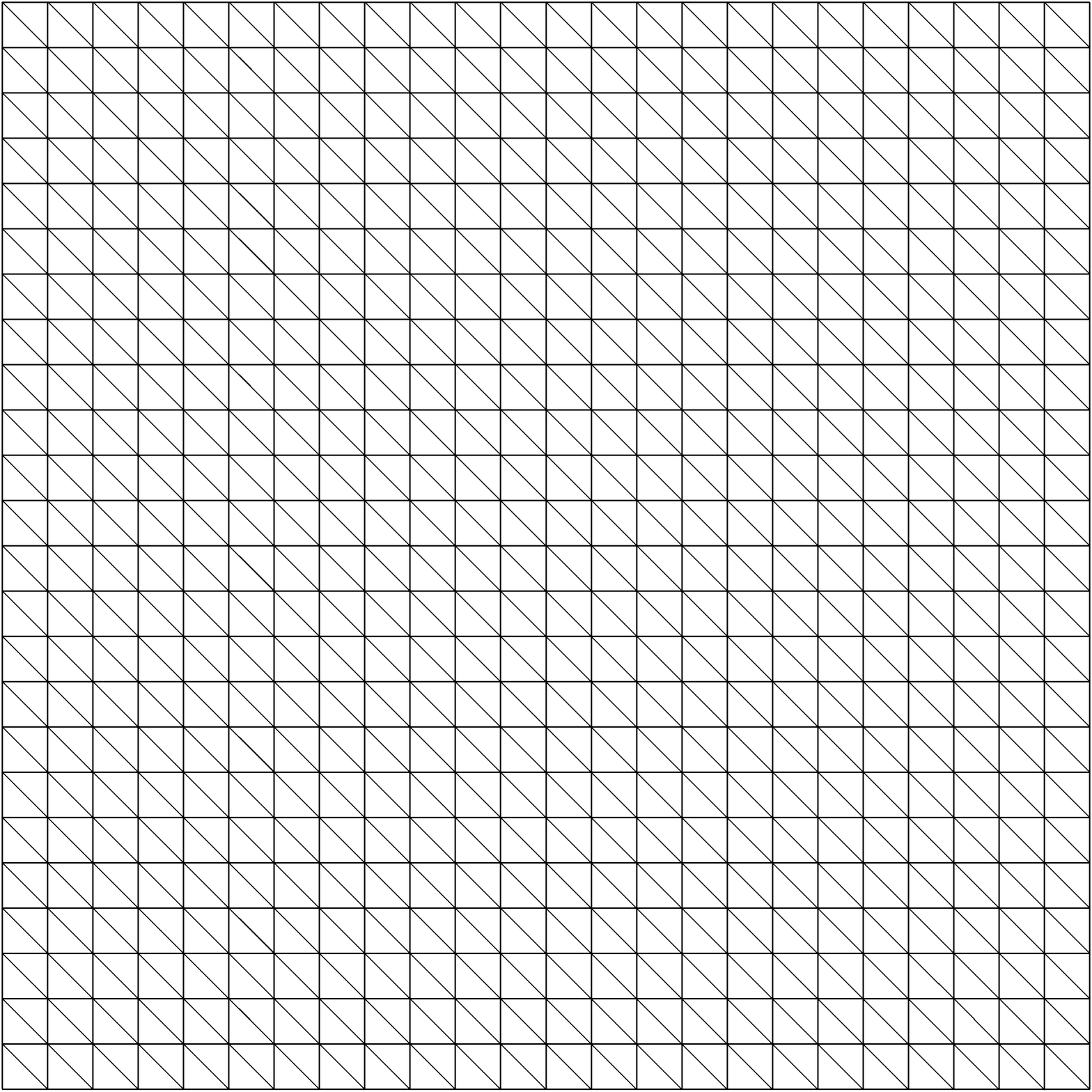, width = 0.2\textwidth}}
\subfigure[]{\label{fig:meshes_chinosi_prob_shearstress_g} \epsfig{file = 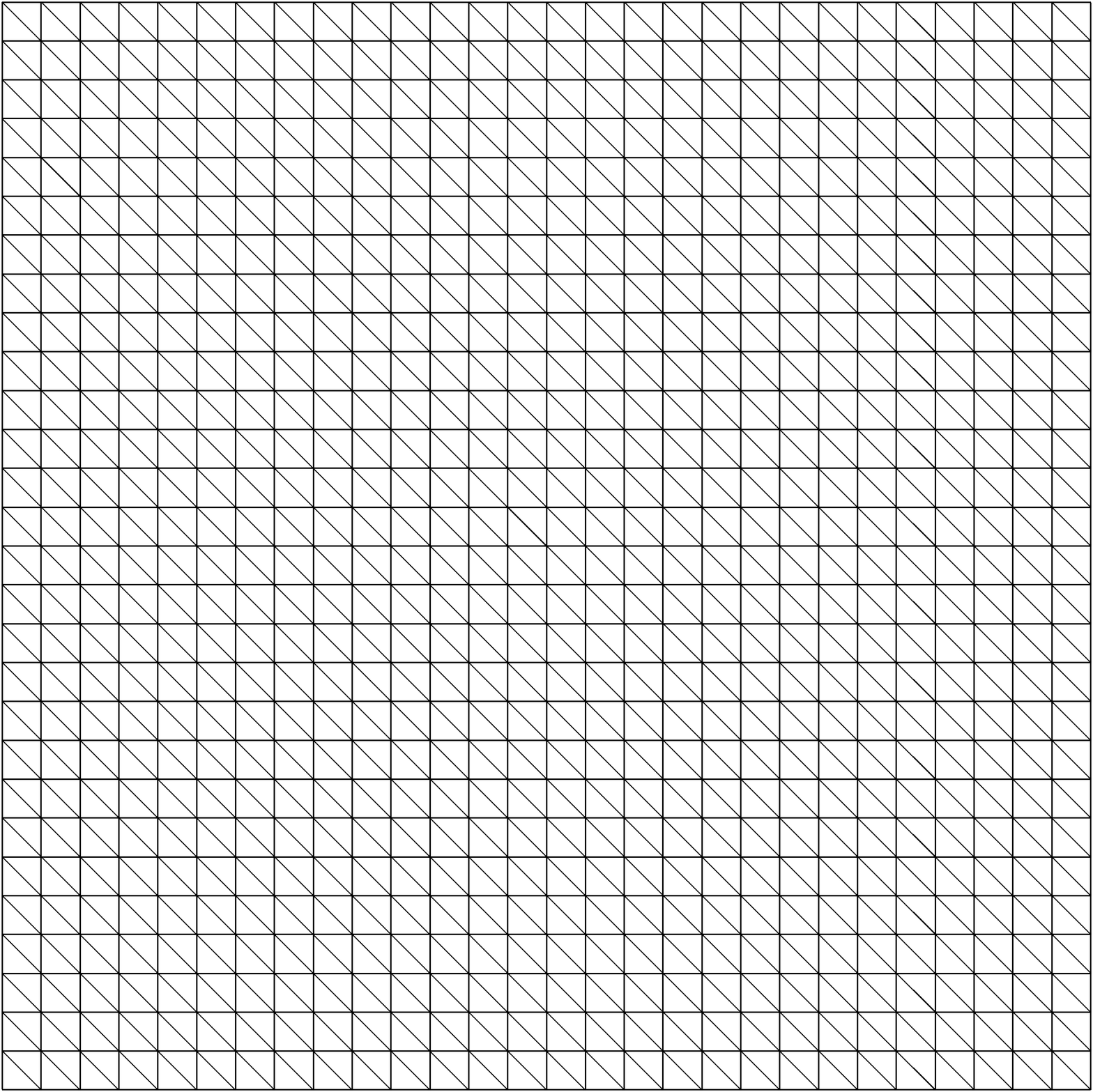, width = 0.2\textwidth}}
\subfigure[]{\label{fig:meshes_chinosi_prob_shearstress_h} \epsfig{file = 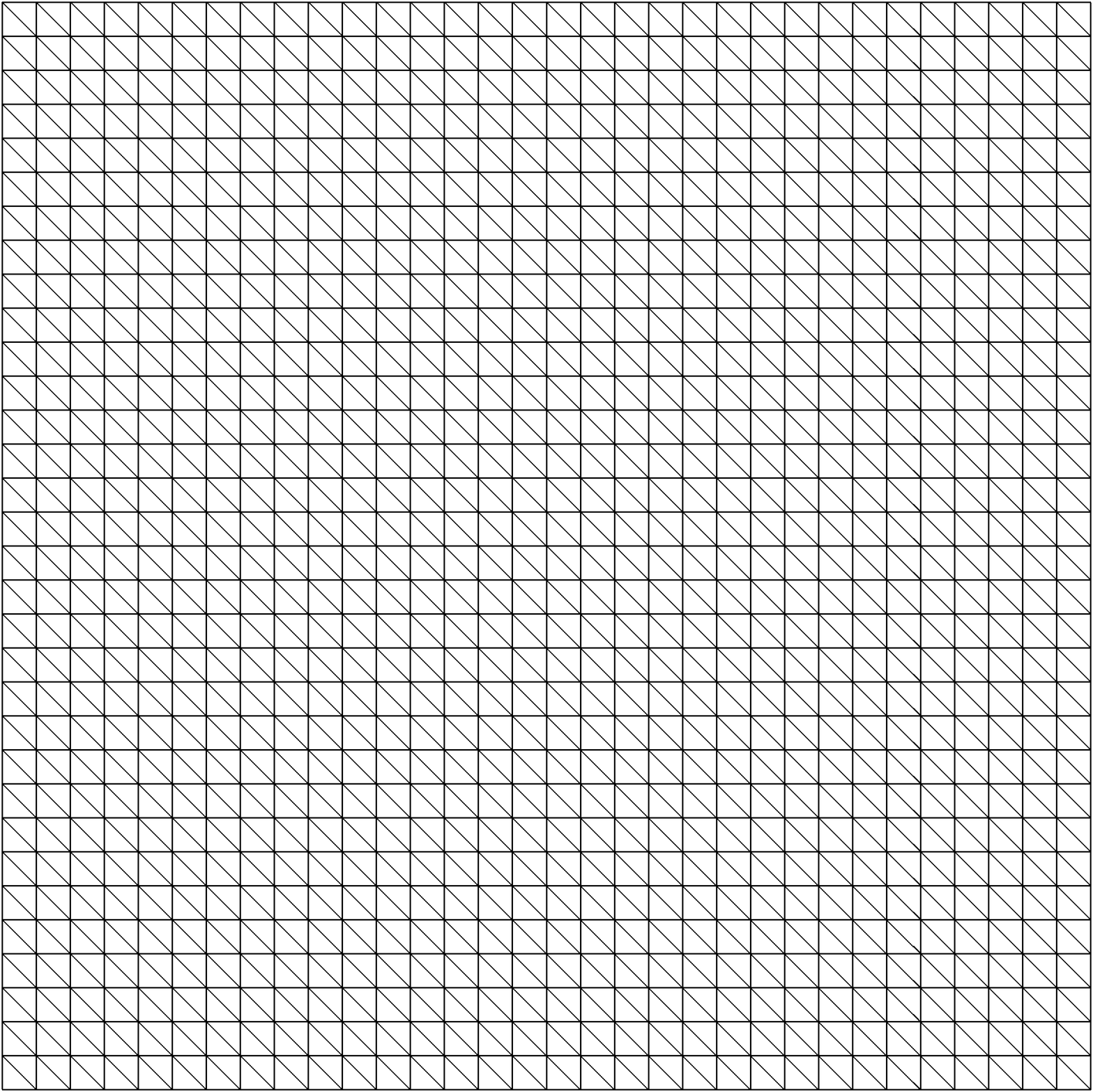, width = 0.2\textwidth}}}
\caption{Structured integration meshes to assess the performance of the transverse shear stress predictions in the \texttt{VANP} formulation.}
\label{fig:meshes_chinosi_prob_shearstress}
\end{figure}

\alejandro{The performance of the recovered scaled transverse shear stress predictions is now assessed. The problem already presented in~\sref{sec:numexamples_chinosi} is used to this aim. In addition to the unstructured integration meshes shown in~\fref{fig:meshes_chinosi_prob}, we consider the set of structured integration meshes depicted in~\fref{fig:meshes_chinosi_prob_shearstress}. It is recalled that the scaled transverse shear stress variable is directly recovered at the nodes using~\eref{eq:nodal_shearstress_lumped} after the primitive variables are computed.}

\begin{figure}[!tbhp]
\centering
\mbox{
\subfigure[]{\label{fig:norms_chinosi_prob_shearstress_a} \epsfig{file = 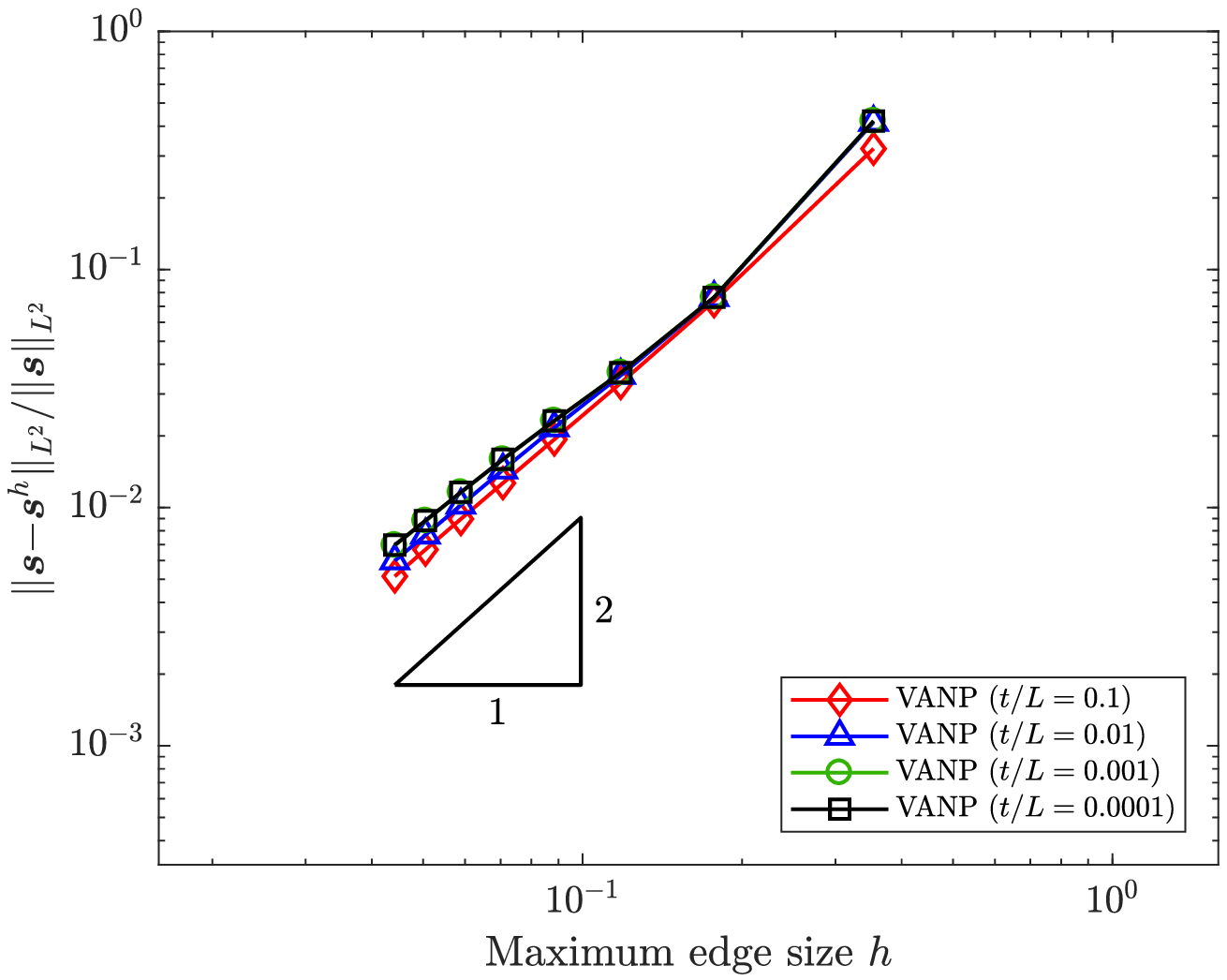, width = 0.5\textwidth}}
\subfigure[]{\label{fig:norms_chinosi_prob_shearstress_b} \epsfig{file = 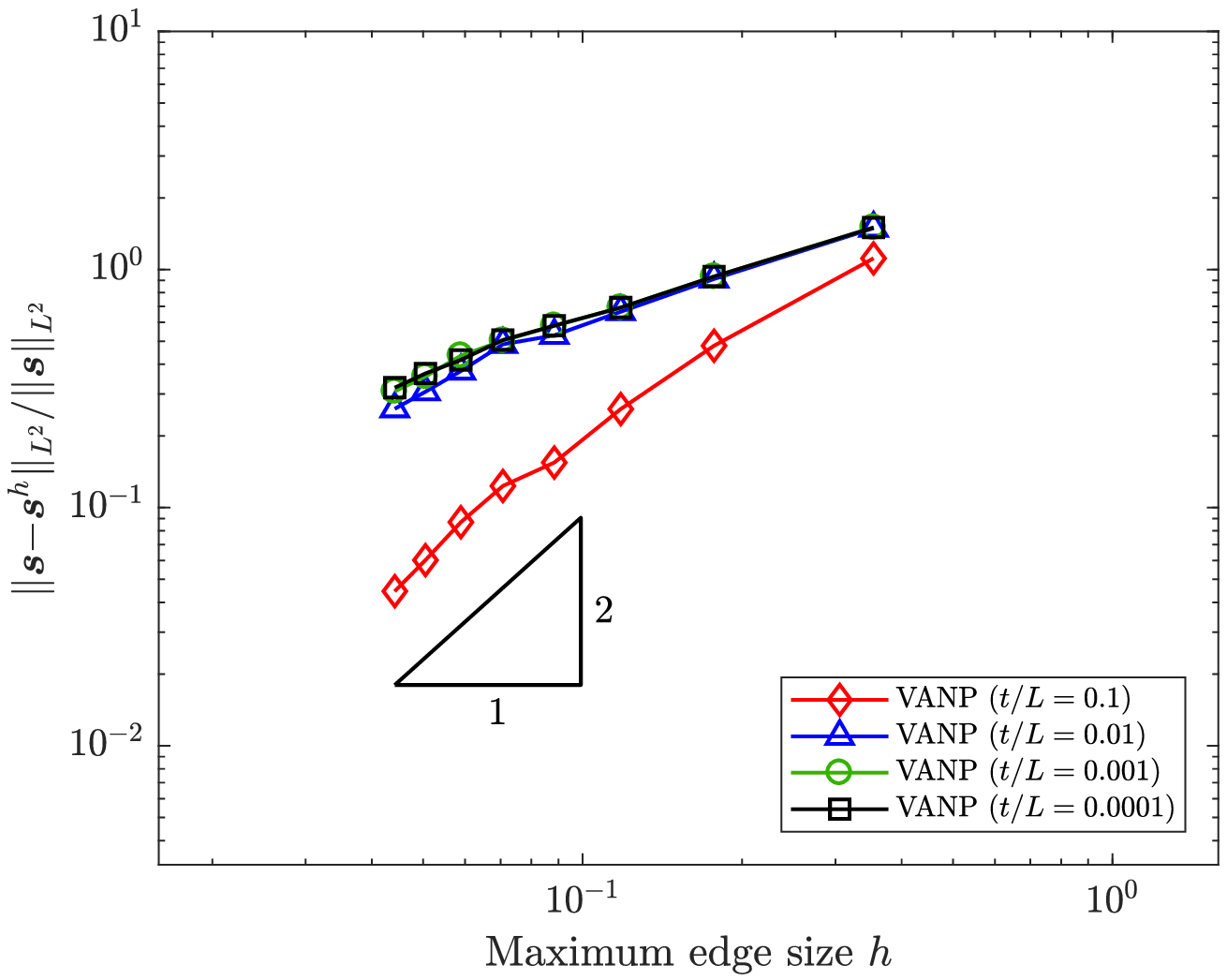, width = 0.5\textwidth}}
}
\caption{$L^2$-norm of the nodal error of the scaled transverse shear stress solution using (a) structured  and (b) unstructured integration meshes in the \texttt{VANP} formulation. For the structured integration meshes, optimal rates of convergence are delivered irrespectively of the plate thickness. The unstructured integration meshes deliver the optimal rate of convergence only for the plate with normalized thickness $t/L=0.1$, whereas the reminder plates converge at about half of the optimal rate.}
\label{fig:norms_chinosi_prob_shearstress}
\end{figure}

\alejandro{
\fref{fig:norms_chinosi_prob_shearstress} presents the $L^2$-norm of the nodal error of the scaled transverse shear stress solution. The optimal convergence rate is delivered by the \texttt{VANP} formulation irrespectively of the plate thickness when the structured integration meshes are used (\fref{fig:norms_chinosi_prob_shearstress_a}), whereas (with the exception of the thicker plate considered) the convergence rate deteriorates to a rate of about half of its optimal value when the unstructured integration meshes are used (\fref{fig:norms_chinosi_prob_shearstress_b}) resulting in degraded accuracy. Notwithstanding this deteriorated performance, the predicted solutions for the primitive variables is excellent and optimally convergent as it was shown in~\sref{sec:numexamples_chinosi}. We also stress that this poor convergence behavior of the scaled transverse shear stress variable should not be a concern as its uniform convergence in the $L^2$-norm is in general very difficult to achieve~\cite{chapelle:2011:FEASH}.}

\section{Concluding Remarks}
\label{sec:conclusions}

In this paper, a volume-averaged nodal projection (\texttt{VANP}) method for the solution of Reissner-Mindlin plate problems using primitive variables (i.e., rotations and transverse displacement) was presented. The proposed approach relies on the construction of a projection operator that permits the computation of the shear strain in terms of the primitive variables without presenting shear-locking issues in the limit of the thin-plate theory. The \texttt{VANP} method uses linear maximum-entropy approximations and bubble-like enrichment of the rotation degrees of freedom is added for stability purposes. A special integration scheme on triangular meshes was developed to fix integration errors in the computation of the meshfree stiffness matrices. The assessing of the \texttt{VANP} formulation through several benchmark problems, which included \alejandro{a zero shear deformation patch test}, a circular plate subjected to a uniform load, a square plate subjected to a nonuniform load and a parallelogram plate subjected to a uniform load, confirmed the accuracy and optimal convergence of the \texttt{VANP} approach for a wide range of plate thicknesses \alejandro{without experiencing shear-locking issues}.

\alejandro{Further improvement of the numerical integration of the stiffness matrix and force vector is being explored by developing a nodal integration technique. From a mathematical standpoint, the construction of error estimates for the \texttt{VANP} approach would help in understanding its optimal performance. The extension of the \texttt{VANP} approach to the von K\'arm\'an theory for nonlinear plates is worth being developed and explored. These topics will be addressed in subsequent works.}


\section*{Acknowledgement}
AOB acknowledges the research support of the Chilean National Fund for Scientific and Technological Development (FONDECYT) through grant CONICYT/FONDECYT No. 1181192.

\section*{References}
\bibliographystyle{elsarticle-num}
\bibliography{meshfree,books,near-incomp,RM-plates}

\end{document}